\documentclass[12pt]{article}
\usepackage{amsfonts}
\usepackage{mathrsfs}
\usepackage{amsmath}
\usepackage{latexsym}
\usepackage{amssymb}
\usepackage{amsthm}
\usepackage{amscd}

\title{On Pointed Hopf Algebras with Weyl Groups of Exceptional Type}
\author{  \small
Shouchuan Zhang $^{a, b}$,  Yao-Zhong Zhang $^b$,    \ \ Peng Wang $^a$,
\ \ Jing Cheng $^a$,\ \  Hui Yang$^a$ \\
\small $a$. Department  of Mathematics,
Hunan University\\  \small   Changsha  410082, \ P.R. China \\
\small $b$. School of Mathematics and Physics, The University of Queensland\\
\small Brisbane 4072, Australia\\
 }
\date{}

\begin{document}
\newtheorem{Proposition}{\quad Proposition}[section]
\newtheorem{Theorem}{\quad Theorem}
\newtheorem{Definition}[Proposition]{\quad Definition}
\newtheorem{Corollary}[Proposition]{\quad Corollary}
\newtheorem{Lemma}[Proposition]{\quad Lemma}
\newtheorem{Example}[Proposition]{\quad Example}
\newtheorem{Remark}[Proposition]{\quad Remark}

\maketitle \addtocounter{section}{-1}

\numberwithin{equation}{section}

\date{}

\begin {abstract} All  $-1$-type pointed Hopf algebras and central
quantum linear  spaces with Weyl groups of exceptional type are
found. It is proved that  every non $-1$-type pointed Hopf algebra
with real $G(H)$ is infinite dimensional and every central quantum
linear space over finite group is finite dimensional. It is proved
that except a few cases Nichols algebras of reducible Yetter-Drinfeld modules over
 Weyl groups of exceptional type are infinite dimensional.


keywords: Quiver, Hopf algebra, Weyl group.
\end {abstract}

\section{Introduction}\label {s0}

This article is to contribute to the classification of
finite-dimensional complex pointed Hopf algebras $H$ with Weyl
groups of exceptional type. The classification of finite dimensional
pointed Hopf algebra with finite abelian groups has been completed ( see
\cite{AS98, AS02, AS00, AS05, He06}). Papers \cite { AG03, Gr00,
AZ07, Fa07, AF06, AF07} considered some  non-abelian cases, for
example, symmetric group, dihedral group, alternating group and
 the Mathieu simple groups. It was shown in \cite {HS} that every
Nichols algebra of reducible Yetter-Drinfeld module over
non-commutative finite  simple group and symmetric group is infinite
dimensional.

In this paper we find all  $-1$-type pointed Hopf algebras and quantum linear  spaces with Weyl
groups of exceptional type. We  show that  every non $-1$-type
pointed Hopf algebra is infinite dimensional and every quantum linear space
is finite dimensional.
It is  desirable to do this in view of the importance of Weyl groups in the theories of
Lie groups, Lie algebras and algebraic groups.
 We first give the
relation between the bi-one Nichols algebra $\mathfrak{B} ({\mathcal O}_s, \rho
)$ introduced in \cite {Gr00, AZ07, AHS08, AFZ08} and the arrow Nichols algebra introduced in
\cite {CR97, CR02, ZZC04, ZCZ08}.
\cite {ZWCY08a, ZWCY08b} applied the software GAP to  compute
 the representatives of
conjugacy classes, centralizers of these representatives  and
character tables of these centralizers in Weyl groups of exceptional
type. Using the results in  \cite {ZWCY08a, ZWCY08b} and the
classification theorem of quiver Hopf algebras and Nichols algebras
in \cite [Theorem 1]{ZCZ08} we find all $-1$-type pointed Hopf
algebras and quantum linear spaces with Weyl groups of exceptional
type. We  prove that Nichols algebras of reducible Yetter-Drinfeld
modules over Weyl groups of exceptional type are infinite
dimensional except a few cases by applying  \cite [Theorem 8.2,
8.6]{HS}.

This paper is organized as follows. In section 1 it is shown that
bi-one arrow Nichols algebras and $\mathfrak{B}({\mathcal O}_s,
\rho)$ introduced in \cite {DPR,  Gr00, AZ07, AHS08, AFZ08} are the
same up to isomorphisms. In  section 2 it is proved that  every non
$-1$-type pointed Hopf algebra with real $G(H)$ is infinite
dimensional. In  section 3 it is shown that every central quantum
linear space
 is  finite dimensional with an  arrow {\rm PBW} basis. In  section 4 the programs
to  compute the representatives of
conjugacy classes, centralizers of these representatives  and character tables
of these centralizers
in Weyl groups
of exceptional type are given. In  section 5 all $-1$- type bi-one Nichols algebras over Weyl groups of  exceptional type
 up to  graded pull-push
{\rm YD} Hopf algebra isomorphisms  are listed in tables. In  section 6 all $-1$- type bi-one Nichols algebras over Weyl groups of  exceptional type
 up to  graded pull-push
{\rm YD} Hopf algebra isomorphisms  are listed in tables. In  section 7
 all central  quantum linear spaces over Weyl groups of exceptional type are found.
In  section 8 it is proved that except a few cases Nichols algebras of reducible
Yetter-Drinfeld modules over
 Weyl groups of exceptional type are infinite dimensional.

\section*{Preliminaries And Conventions}

Throughout this paper let  $k$ be  the complex field; $G$ be a finite group;
$\hat{{ G}}$
denote the set of all isomorphic classes of irreducible
representations of group $G$; $G^{s}$ denote the centralizer of $s$; $Z(G)$ denote the center of
$G$. For $h \in G$ and an isomorphism $\phi $ from $G$ to $G'$,
define a map $\phi _h$ from $G$ to $G'$ by sending $x$ to $\phi
(h^{-1}xh)$ for any $x\in G$. Let $s^G$ or ${\mathcal O}_s$ denote the conjugacy class containing $s$ in $G$.
The Weyl groups of $E_6,$ $ E_7$, $E_8$, $F_4$ and $G_2$ are called
{\it Weyl groups of  exceptional type}. Let {\rm deg }$\rho$ denote the dimension of the representation space
$V$ for a representation  $(V, \rho).$

Let ${\mathbb N}$ and ${\mathbb Z}$ denote the sets of  all positive integers and
all integers, respectively. For a set $X$,
we denote by $|X|$ the  number of elements in $X$.  If $X = \oplus _{i\in I} X_{(i)}$ as vector
spaces, then we denote by $\iota _i$ the natural injection from
$X_{(i)}$ to $X$ and by $\pi _i$ the corresponding projection from $X$
to $X_{(i)}$. We will use $\mu$ to denote the multiplication  of an
algebra and use $\Delta$ to denote the comultiplication of a
coalgebra. For a (left or right) module and a (left or right)
comodule, denote by $\alpha ^-$, $\alpha ^+$, $\delta ^-$ and
$\delta ^+$ the left module, right module, left comodule and right
comodule structure maps, respectively. The Sweedler's sigma
notations for coalgebras and comodules are $\Delta (x) = \sum
x_{1}\otimes x_{2}$, $\delta ^- (x)= \sum x_{(-1)} \otimes
x_{(0)}$, $\delta ^+ (x)= \sum x_{(0)} \otimes x_{(1)}$.

A quiver $Q=(Q_0,Q_1,s,t)$ is an oriented graph, where  $Q_0$ and
$Q_1$ are the sets of vertices and arrows, respectively; $s$ and $t$
are two maps from  $Q_1$ to $Q_0$. For any arrow $a \in Q_1$, $s(a)$
and $t(a)$ are called {\it its start vertex and end vertex}, respectively,
and $a$ is called an {\it arrow} from $s(a)$ to $t(a)$. For any $n\geq 0$,
an $n$-path or a path of length $n$ in the quiver $Q$ is an ordered
sequence of arrows $p=a_na_{n-1}\cdots a_1$ with $t(a_i)=s(a_{i+1})$
for all $1\leq i\leq n-1$. Note that a 0-path is exactly a vertex
and a 1-path is exactly an arrow. In this case, we define
$s(p)=s(a_1)$, the start vertex of $p$, and $t(p)=t(a_n)$, the end
vertex of $p$. For a 0-path $x$, we have $s(x)=t(x)=x$. Let $Q_n$ be
the set of $n$-paths. Let $^yQ_n^x$ denote the set of all $n$-paths
from $x$ to $y$, $x, y\in Q_0$. That is, $^yQ_n^x=\{p\in Q_n\mid
s(p)=x, t(p)=y\}$.

A quiver $Q$ is {\it finite} if $Q_0$ and $Q_1$ are finite sets. A
quiver $Q$ is {\it locally finite} if $^yQ_1^x$ is a finite set for
any $x, y\in Q_0$.

Let ${\mathcal K}(G)$ denote the set of
conjugate classes in $G$. A formal sum $r=\sum_{C\in {\mathcal
K}(G)}r_CC$  of conjugate classes of $G$  with cardinal number
coefficients is called a {\it ramification} (or {\it ramification
data} ) of $G$, i.e.  for any $C\in{\mathcal K}(G)$, \  $r_C$ is a
cardinal number. In particular, a formal sum $r=\sum_{C\in {\mathcal
K}(G)}r_CC$  of conjugate classes of $G$ with non-negative integer
coefficients is a ramification of $G$.

 For any ramification $r$ and  $C \in {\mathcal K}(G)$, since $r_C$ is
 a cardinal number,
we can choose a set $I_C(r)$ such that its cardinal number is $r_C$
without loss of generality.
 Let ${\mathcal K}_r(G):=\{C\in{\mathcal
K}(G)\mid r_C\not=0\}=\{C\in{\mathcal K}(G)\mid
I_C(r)\not=\emptyset\}$.  If there exists a ramification $r$ of $G$
such that the cardinal number of $^yQ_1^x$ is equal to $r_C$ for any
$x, y\in G$ with $x^{-1}y \in C\in {\mathcal K}(G)$, then $Q$ is
called a {\it Hopf quiver with respect to the ramification data
$r$}. In this case, there is a bijection from $I_C(r)$ to $^yQ_1^x$,
and hence we write  ${\ }^yQ_1^x=\{a_{y,x}^{(i)}\mid i\in I_C(r)\}$
for any $x, y\in G$ with $x^{-1}y \in C\in {\mathcal K}(G)$.

$(G, r, \overrightarrow \rho, u)$ is called a {\it ramification system
with irreducible representations}  (or {\rm RSR } in short ), if $r$
is a ramification of $G$; $u$ is a map from ${\mathcal K}(G)$ to $G$
with $u(C)\in C$ for any $C\in {\mathcal K}(G)$;  $I_{C} (r, u)$ and
$J_C(i)$ are  sets   with $ \mid\!J_{C} (i )\!\mid $  = $ {\rm deg}
(\rho _C^{(i)})$ and $I_C(r) = \{(i, j) \mid i \in I_C(r, u), j\in
J_C(i)\}$  for any $C\in {\mathcal K}_r(G)$, $i\in I_C(r, u)$;
$\overrightarrow \rho=\{\rho_C^{(i)} \}_ { i\in I_C(r,u ),
C\in{\mathcal K}_r(G)} \ \in \prod _ { C\in{\mathcal K}_r(G)}
(\widehat{    {G^{u(C)}}}) ^{\mid I_{C} (r, u) \mid }$ with
$\rho_C^{(i)} \in \widehat{ {G^{u(C)}}} $ for any
 $ i \in I_C(r, u), C\in {\mathcal
K}_r(G)$.
In this paper we always assume that $I_C(r, u)$ is a finite set for any $C\in
{\mathcal K}_r (G).$
 Furthermore, if $\rho _C^{(i)}$ is a one dimensional representation for any $C\in
{\mathcal K}_r(G)$, then $(G, r, \overrightarrow \rho, u)$ is called a {\it ramification system
with characters }  (or {\rm RSC } $(G, r, \overrightarrow \rho, u)$ in short ) (see \cite [Definition 1.8]{ZZC04}). In this case,
 $ a_{y, x}^{(i, j)}$ is written as $a_{y, x}^{(i)}$ in short  since $J_C(i)$ has only one element.

For ${\rm RSR}(G, r, \overrightarrow \rho, u)$, let $\chi
_{C}^{(i)}$ denote the character of $\rho _C^{(i)}$ for any $i \in
I_C(r, u)$, $C \in {\mathcal K}_r(C)$. If ramification $r = r_CC$
and $I_C(r,u) = \{ i  \}$ then we say that ${\rm RSR}(G, r, \overrightarrow \rho, u)$ is   bi-one, written  as
${\rm RSR } (G, {\mathcal O}_{s}, \rho )$ with $s= u(C)$ and $\rho = \rho _C^{(i)}$ in short,
since $r$ only has one conjugacy class $C$ and $\mid\!I_C(r,u)\!\mid
=1$. Quiver Hopf algebras,  Nichols algebras and
Yetter-Drinfeld modules, corresponding to a bi-one ${\rm RSR}(G, r,
\overrightarrow \rho, u)$, are said to be    bi-one.

If  $(G, r, \overrightarrow \rho, u)$ is an ${\rm RSR}$,  then it is clear that ${\rm RSR } (G, {\mathcal O}_{u(C)},
\rho _C^{(i)})$  is   bi-one  for any $C\in {\mathcal K}$ and
$i \in I_C(r, u)$, which is called a {\it bi-one } sub-${\rm RSR}$ of ${\rm RSR}(G, r, \overrightarrow \rho, u)$,

If  $\phi: A\rightarrow A'$ is  an algebra homomorphism and  $(M,
\alpha ^-)$ is a left $A'$-module, then $M$ becomes a left
$A$-module with the $A$-action given by $a \cdot x =\phi (a) \cdot x
$ for any $a\in A$, $x\in M$, called a {\it pullback } $A$-module through
$\phi$, written as  $_{\phi}M$.  Dually, if  $\phi: C\rightarrow C'$
is a coalgebra homomorphism and  $(M, \delta ^- )$ is  a left
$C$-comodule, then $M$ is a left $C'$-comodule with the
$C'$-comodule structure given by $ {\delta'}^-:=(\phi\otimes{\rm
id})\delta^-$, called  a {\it push-out } $C'$-comodule through $\phi$,
written as  $^{\phi}M$.

If $B$ is a Hopf algebra and $M$ is a $B$-Hopf bimodule, then  we
say that $(B, M)$ is a  Hopf bimodule. For any two Hopf bimodules
$(B,M)$ and $(B', M')$, if $\phi$ is a Hopf algebra homomorphism
from $B$ to $B'$
 and $\psi$ is simultaneously a $B$-bimodule homomorphism from
$M$ to $_\phi M'{}_\phi $ and a $B'$-bicomodule homomorphism from
$^\phi M ^\phi $ to $M'$, then $(\phi, \psi)$ is called a {\it pull-push }
Hopf bimodule homomorphism. Similarly, we say that $(B, M)$ and $(B,
X)$ are a Yetter-Drinfeld ({\rm YD}) module
 and  {\rm YD} Hopf algebra, respectively, if $M$ is a
{\rm YD} $B$-module and $X$ is a braided Hopf algebra in
{\rm YD} category $^B_B {\mathcal YD}$.
 For any two
{\rm YD} modules $(B,M)$ and $(B', M')$, if $\phi$ is a Hopf algebra
homomorphism from $B$ to $B'$,  and $\psi$ is simultaneously a left
$B$-module homomorphism from  $M$ to $_\phi M' $ and a left
$B'$-comodule homomorphism from $^\phi M  $ to $M'$, then $(\phi,
\psi)$ is called a pull-push {\rm YD} module homomorphism. For any
two {\rm YD} Hopf algebras $(B,X)$ and $(B', X')$, if $\phi$ is a
Hopf algebra homomorphism from $B$ to $B'$,   $\psi$ is
simultaneously a left $B$-module homomorphism from  $X$ to $_\phi X'
$ and a left $B'$-comodule homomorphism from $^\phi X  $ to $X'$,
meantime, $\psi$ also is algebra and coalgebra homomorphism from $X$
to $X'$, then $(\phi, \psi)$ is called a pull-push {\rm YD} Hopf
algebra homomorphism (see \cite [the remark after Th.4]{ZZC04}).

For  $s\in G$ and  $(\rho, V) \in  \widehat {G^s}$, here is a
precise description of the {\rm YD} module $M({\mathcal O}_s,
\rho)$, introduced in \cite {Gr00, AZ07}. Let $t_1 = s$, \dots,
$t_{m}$ be a numeration of ${\mathcal O}_s$, which is a conjugacy
class containing $s$,  and let $g_i\in G$ such that $g_i \rhd s :=
g_i s g_i^{-1} = t_i$ for all $1\le i \le m$. Then $M({\mathcal
O}_s, \rho) = \oplus_{1\le i \le m}g_i\otimes V$. Let $g_iv :=
g_i\otimes v \in M({\mathcal O}_s,\rho)$, $1\le i \le m$, $v\in V$.
If $v\in V$ and $1\le i \le m $, then the action of $h\in G$ and the
coaction are given by
\begin {eqnarray} \label {e0.11}
\delta(g_iv) = t_i\otimes g_iv, \qquad h\cdot (g_iv) =
g_j(\gamma\cdot v), \end {eqnarray}
 where $hg_i = g_j\gamma$, for
some $1\le j \le m$ and $\gamma\in G^s$. The explicit formula for
the braiding is then given by
\begin{equation} \label{yd-braiding}
c(g_iv\otimes g_jw) = t_i\cdot(g_jw)\otimes g_iv =
g_{j'}(\gamma\cdot v)\otimes g_iv\end{equation} for any $1\le i,j\le
m$, $v,w\in V$, where $t_ig_j = g_{j'}\gamma$ for unique $j'$, $1\le
j' \le m$ and $\gamma \in G^s$. Let $\mathfrak{B} ({\mathcal O}_s,
\rho )$ denote $\mathfrak{B} (M ({\mathcal O}_s, \rho ))$. $M({\mathcal O}_s,
\rho )$ is a simple {\rm YD} module (see \cite [Section 1.2 ] {AZ07}).
Furthermore, if $\chi$ is the character of $\rho$, then we also
denote $\mathfrak{B} ({\mathcal O}_s, \rho )$ by $\mathfrak{B}
({\mathcal O}_s, \chi )$.


\section {Relation between bi-one arrow Nichols algebras and $\mathfrak{B}({\mathcal O}_s, \rho)$}
\label {s1}

In this section it is shown that
bi-one arrow Nichols algebras and
$\mathfrak{B}({\mathcal O}_s, \rho)$ introduced in \cite {Gr00, AZ07, AHS08, AFZ08} are the
same  up to isomorphisms.

For any ${\rm RSR} (G, r, \overrightarrow \rho, u)$, we can
construct an arrow Nichols algebra
 $\mathfrak{B} (kQ_1^1, ad (G, r, \overrightarrow{\rho},
$ $u))$ ( see \cite [Pro. 2.4] {ZCZ08}), written as $\mathfrak{B} (G,
r, \overrightarrow{\rho}, $ $ u)$ in short.
Let us recall the  precise description of arrow {\rm YD} module.
For an ${\rm RSR}(G, r, \overrightarrow \rho, u)$ and a $kG$-Hopf
bimodule $(kQ_1^c, G, r, \overrightarrow {\rho}, u)$ with the module
operations $\alpha^-$ and $\alpha^+$, define a new left $kG$-action
on $kQ_1$ by
$$g\rhd x:=g\cdot x\cdot g^{-1},\ g\in G, x\in kQ_1,$$
where $g\cdot x=\alpha^-(g\otimes x)$ and $x\cdot
g=\alpha^+(x\otimes g)$ for any $g\in G$ and $x\in kQ_1$. With this
left $kG$-action and the original left (arrow) $kG$-coaction
$\delta^-$, $kQ_1$ is a  Yetter-Drinfeld   $kG$-module. Let  $Q_1^1:=\{a\in Q_1 \mid
s(a)=1\}$, the set of all arrows with starting vertex $1$. It is
clear that $kQ_1^1$ is a  Yetter-Drinfeld   $kG$-submodule of $kQ_1$, denoted by
$(kQ_1^1, ad(G, r, \overrightarrow {\rho}, u))$, called the arrow {\rm YD} module.

\begin {Lemma} \label {1.1} For any $s\in G$ and  $\rho \in  \widehat
{G^s}$, there exists a bi-one arrow Nichols algebra $\mathfrak{B} (G,
r, \overrightarrow{\rho}, u)$ such that
$$\mathfrak{B}({\mathcal O}_s, \rho)
\cong \mathfrak{B} (G, r, \overrightarrow{\rho}, u)$$ as graded
braided Hopf algebras in $^{kG}_{kG} \! {\mathcal YD}$.

\end {Lemma}
{\bf Proof.}  Assume that $V$ is the representation space of $\rho$
 with $\rho (g) (v)= g\cdot v$ for any $g\in G, v\in V$. Let $C = {\mathcal O_s}$,
$r = r_C C$, $r_C ={\rm deg } \rho$, $u(C)
= s$, $I_C(r, u)= \{1 \}$ and  $ (v)\rho _C^{(1)}(h)= \rho (h^{-1})(v)$
for any $h\in G$, $v\in V$. We get a bi-one arrow Nichols algebra $
\mathfrak{B}( G, r, \overrightarrow{\rho}, u)$.

We now  only need to show that $M ({\mathcal O}_s, \rho) \cong
(kQ_1^1, ad( G, r, \overrightarrow{\rho}, u)) $ in $^{kG}_{kG} \!
{\mathcal YD}$. We recall the notation in \cite [Proposition
1.2]{ZCZ08}. Assume  $J_C(1) = \{1, 2, \cdots, n\}$ and $X_C^{(1)} =
V$ with basis $\{x_C^{(1, j)} \mid j =1, 2, \cdots, n\}$ without
loss of  generality. Let $v_j $ denote $x_C^{(1,j)}$ for convenience. In
fact,
 the  left and right coset decompositions of $G^s$ in $G$ are
\begin {eqnarray} \label {e1.1.1}
G =\bigcup_{i=1}^m g_{i} G^s    \ \ \hbox {and } \ \    G
&=&\bigcup_{i=1} ^m G^sg_{i} ^{-1} \ \ ,
\end {eqnarray} respectively.

 Let $\psi$ be a map from $M({\mathcal O}_s,
\rho)$ to $(kQ_1^1, {\rm ad } (G, r, \overrightarrow{\rho}, u))$ by
sending $ g_i v_j$ to $a _{t_i, 1}^{(1, j)}$ for any $1\le i\le m,
1\le j \le n$. Since the dimension is $mn$, $\psi$ is a
bijective. See
\begin {eqnarray*}
\delta ^- (\psi (g_i v_j)) &=& \delta ^- (a _{t_i, 1}^{(1,j)})\\
&=& t_i \otimes a_{t_i, 1}^{(1, j)} = (id \otimes \psi) \delta ^-
(g_iv_j). \end {eqnarray*} Thus $\psi$ is a $kG$-comodule
homomorphism. For any $h\in G$, assume  $hg_i = g_{i'} \gamma $ with
$\gamma \in G^s$. Thus $g_i ^{-1} h^{-1} = \gamma ^{-1}
g_{i'}^{-1}$, i.e. $\zeta _i (h^{-1})= \gamma ^{-1}$, where
$\zeta_i$ was defined  in \cite [(0.3)]{ZZC04}. Since $\gamma \cdot
x^{(1,j)} \in V$, there exist $k_{C, h^{-1}}^{(1, j, p)}\in k $,
$1\le p \le n$, such that $\gamma \cdot x^{(1,j)} =\sum _{p=1}^n
k_{C, h^{-1}} ^{(1, j, p)} x^{(1, p)}$. Therefore
\begin {eqnarray}
\label {e1.11}x^{(1,j)}\cdot \zeta _i (h^{-1}) &=& \gamma \cdot
x^{(1,j)} \ \ ( \hbox {by definition of } \rho _C^{(1)} ) \nonumber
\\ &=&\sum _{p=1}^n k_{C, h^{-1}} ^{(1, j, p)} x^{(1,
p)}.\end {eqnarray}

See
\begin {eqnarray*}
 \psi (h \cdot g_iv_j) &=& \psi ( g _{i'}(\gamma v_j)) \\
&=& \psi (g_{i'} (\sum _{p=1}^n k_{C, h^{-1}} ^{(1, j, p)} v_p) )\\
 &= &\sum _{p=1}^n k_{C, h^{-1}} ^{(1, j, p)} a _{t_{i'}, 1} ^{(1,
p)}
\end {eqnarray*}
and
\begin {eqnarray*}
h \rhd (\psi (g_i v_j)) &=& h \rhd (a _{t_i, 1}^{(1,j)})\\
&=& a_{ht_i, h}^{(1, j)} \cdot h^{-1} \\
&=& \sum _{p =1}^n k _{C, h^{-1}} ^{(1, j, p)} a_{t_{i'}, 1}^{(1,
p)} \ \ (\hbox {by \cite [Pro.1.2]{ZCZ08} and }  (\ref {e1.11})).
\end {eqnarray*} Therefore
$\psi $ is a $kG$-module homomorphism. $\Box$

 Therefore  we also say that
$\mathfrak{B}({\mathcal O}_s, \rho)$ is  a bi-one Nichols Hopf
algebra.

\begin {Remark} \label {1.2'} The representation $\rho$ in   $\mathfrak{B}({\mathcal O}_s, \rho)$   introduced in \cite {Gr00, AZ07} and
$\rho _C^{(i)}$ in {\rm RSR} are different.  $\rho (g)$ acts on its representation
space from  the left and $\rho _C^{(i)} (g)$
acts on its representation space from  the right.
\end {Remark}


$s \in G$ is  {\it  real } if $s$ and $s^{-1}$ are in the same
conjugacy class. If every element in $G$ is real, then $G$ is real. Obviously, Weyl groups are real.

\begin {Lemma}\label {1.2} Assume that $s\in G$ is real and $\chi$ is the character of
$\rho \in \widehat {G^s}$. If $\chi (s)\neq-{\rm deg} (\rho)$ or the
order of $s$ is odd, then ${\rm dim}
\mathfrak{B}(\mathcal{O}_{s},\rho)=\infty$.
\end {Lemma}

{\bf Proof}.  If the order of $s$ is odd, it follows  from \cite [Lemma
2.2]{AZ07} and \cite [Lemma 1.3] {AF07}.  Now assume
that $\chi (s)\neq-{\rm deg} (\rho)$. Since $\rho(s)=q_{ss}{\rm id}$,
$\chi(s)=q_{ss}({\rm deg }(\rho))$. Therefore $q_{ss}\neq-1$ and
${\rm dim}\mathfrak{B}(\mathcal{O}_{s},\rho)=\infty$ by \cite [Lemma
2.2]{AZ07} and \cite [Lemma 1.3] {AF07}. $\Box$.

\begin {Lemma} \label {1.8}
$(kG, \mathfrak{B}(\mathcal{O}_{s},\rho)) \cong
(kG', \mathfrak{B}(\mathcal{O}_{s'},\rho'))$  as graded pull-push
{\rm YD} Hopf algebras if and only if there exist
$h\in G$ and a group isomorphism
 $\phi$ from  $G$ to $G'$ such that $\phi (h^{-1} sh) = s' $ and $\rho ' \phi _h \cong  \rho$,
  where $\phi _h (g) = \phi (h^{-1}gh)$ for any $g\in G.$
\end {Lemma}
{\bf Proof.} Let $C$ and $ C' $ be  conjugacy classes of $G$ and $G'$, respectively;
$r= r_C C$ and $r'= r_{C'} C'$ be  ramifications of $G$ and $G'$, respectively. Applying
Lemma \ref {1.1}, we only need show that $(kG, \mathfrak{B}(G, r,
\overrightarrow{\rho}, u)) \cong
 (kG', \mathfrak{B}(G, r',
\overrightarrow{\rho'}, u'))$ as   graded pull-push {\rm YD}
Hopf algebras if and only if there exist $h\in G$ and a
group automorphism
  group isomorphism $\phi$ from  $G$ to $G'$ such that $\phi (h^{-1} u(C)h) = u'(C') $ and $\rho ' {}^{(i')}_{C'}\phi _h
 \cong
 \rho_{C}^{(i)}$.
Applying \cite [Theorem 4]{ZCZ08}, we only need show that ${\rm
RSR}(G, r, \overrightarrow{\rho}, u) \cong
 {\rm RSR}(G, r',
\overrightarrow{\rho'}, u')$  if and only if there exist $h\in G$
and a group isomorphism
 $\phi$ from  $G$ to $G'$ such that $\phi (h^{-1} u(C)h) = u'(C') $ and $\rho ' {}^{(i')}_{C'}\phi _h
 \cong
 \rho_{C}^{(i)}$. This is clear. $\Box$

If we define  a relation on group $G$ as follows: $x \sim y$ if and only if there exists a group automorphism $\phi $
of $G$ such that $\phi (x)$ and $y$ are
contained in the same conjugacy  class, then this is an equivalent relation. Let set $\{s_i \mid i\in \Omega \}$ denote all
representatives of the equivalent
classes, which is called  the
 {\it  representative system of conjugacy classes } of $G$ under isomorphism relations, or
   the representative system of iso-conjugacy classes of $G$ in short.

\begin {Proposition} \label {1.9} Let  $\{s_i \mid i\in \Omega \}\subseteq
G$ be the representative system of iso-conjugacy classes of $G$. Then
 $\{   {\mathfrak B} (  {\mathcal O}_{s_i}, \rho) \mid i \in
\Omega, \rho \in \widehat {G^{s_i}} \}$ are all representatives of
the bi-one Nichols algebra over $G$, up to graded pull-push {\rm YD}
Hopf algebra isomorphisms.
\end {Proposition}
{\bf Proof.} If ${\mathfrak B} (  {\mathcal O}_{s}, \rho)$ is a
bi-one Nichols Hopf algebra over $G$, then there exist $i\in \Omega
$,   $\phi \in {\rm Aut } (G)$ and  $h\in G$  such that $\phi _h
(s) = s_i$. Let $\rho ' = \rho ( \phi _{h})^{-1} $. By Lemma \ref
{1.8}, $( kG, {\mathfrak B} (  {\mathcal O}_{s}, \rho)) \cong ( kG, {\mathfrak B}
( {\mathcal O}_{s_i}, \rho'))$ as graded pull-push {\rm YD} Hopf
algebras.

It follows from Lemma \ref {1.8} that $(kG, \mathfrak{B}({\mathcal
O}_{s_i}, \rho))$ and $ (kG,  \mathfrak{B}({\mathcal O}_{s_j}, \rho'))$ are
not  graded pull-push {\rm YD} Hopf algebra isomorphisms when $i
\not= j$ and $i, j \in \Omega$. $\Box$


\section {Diagram}\label {s2}
In this section it is proved that  every non $-1$-type pointed Hopf
algebra with real $G(H)$ is infinite dimensional.

If $H$ is a graded Hopf algebra, then there exists the diagram of
$H$, written ${\rm diag} (H)$, (see \cite [Section 3.1]{ZZC04} and
\cite {Ra85}). If $H$ is a pointed Hopf algebra, then  the coradical
filtration Hopf algebra ${\rm gr} (H)$ is a graded Hopf algebra. So
${\rm gr} (H)$ has the diagram, written ${\rm diag} _{\rm filt}
(H)$, called the {\it filter diagram } of $H$. ${\rm diag} _{\rm
filt} (H)$ is written as ${\rm diag} (H)$ in short when it does not
cause confusion (see
  \cite [Introduction]{AS98} ).

A graded coalgebra $C = \oplus _{n=0} ^{\infty} C_{(n)}$ is
  strictly graded if  $C_{(0)} =k$ and $C_{(1)}= P(C)$ (see \cite
  [P232]{Sw69}).

\begin {Proposition}\label {1.3}
If  $H = \oplus _{n=0} ^{\infty} H_{(n)}$ is a  graded Hopf algebra
and $R := {\rm diag} (H)$ is  strictly graded as coalgebras, then $H
\cong {\rm gr} H$ as graded Hopf algebras.
\end {Proposition}
{\bf Proof.} By \cite [Lemma 2.5] {AS98}, $H$ is coradically
graded, i.e. $H_{m} = \oplus _{i=0}^{m} H_{(i)}$ for $m=0, 1, 2,
\cdots, $ where $H_0 \subseteq H_1\subseteq H_2 \subseteq H_3 \subseteq \cdots $
is the coradical filtration  of $H$. Define  a map $\psi $ from $H$ to ${\rm gr } H$ by
sending $a $ to $a + H_{m-1}$ for any $a \in H_{(m)}$ and $m = 0, 1,
2, 3, \cdots .$ Note $H_{-1} := 0.$ Obviously, $\psi $ is bijective.
If  $a\in H_{(m)}$, then there exist $a_s ^{(j)},$ $ b_s ^{(j)} \in
H_{(j)} $ for $0\le j \le m $,  $1\le s \le n_j$, such that  $\Delta
(a) = \sum _{i=0}^{m} \sum _{s=1}^{n_i} a_s^{(i)} \otimes
b_s^{(m-i)}$. See
\begin {eqnarray*}
(\psi \otimes \psi)\Delta (a) &=& (\psi \otimes \psi) \sum
_{i =0}^{m} \sum _{s=1}^{n_i}a_s^{(i)}  \otimes b_s^{(m-i)}\\
&=& \sum _{i =0}^{m}\sum _{s=1}^{n_i} (a_s^{(i)} +H_{i-1}) \otimes
(b_s^{(m-i)}+ H_{m-i-1})\\
&=& \Delta (a + H_{m-1}) \ \  \\
&& ( \hbox {by the definition of comultiplication of } {\rm gr } H
\hbox { in  \cite
[P229]{Sw69}  } )\\
&=& \Delta \psi (a).
\end {eqnarray*}
Thus $\psi$ is a coalgebra homomorphism. Similarly, $\psi $ is a
algebra homomorphism. $\Box$

 Consequently, every  pointed Hopf algebra of  type one ( since its diagram is Nichols algebra,
 see \cite [Section 2]{ZCZ08}) is isomorphic to its filtration Hopf
 algebra as graded Hopf algebras.

\begin {Lemma}\label {1.4}
If  $R$ is a graded braided Hopf algebra in $^{kG}_{kG} {\mathcal
  YD}$ and is strictly graded as coalgebra gradations, then the
  subalgebra $\bar R$ generated by $R_{(1)}$ as algebras is a Nichols algebra.
  Furthermore, $\bar R$ generated by $R_{(1)}$ as algebras
 in $R$ is a Nichols algebra when $R$ is the filter diagram  of a pointed Hopf algebra $H$.

\end {Lemma}
{\bf Proof.} We show the first claim by the following steps. Let
\begin {eqnarray*}
  x &=&  x^{(1)}x^{(2)} \cdots x^{(n)} \ \ \hbox {and} \ \ y= x ^{(1)}, \ z
  =x^{(2)}\cdots x^{(n)}
  \end {eqnarray*}
   with $x \in R$, $x^{(i)} \in R_{(1)}$ for $i = 1, 2, \cdots, n$.

{\rm (i)} $\bar R$ is $kG$-submodule
 of $R$. In fact $h\cdot x = h \cdot  x^{(1)}x^{(2)}\cdots x^{(n)}$
 =  $(h\cdot x^{(1)})(h\cdot x^{(2)})\cdots ( h\cdot  x^{(n)})$ $\in
  R_{(1) }  R_{(1) } \cdots  R_{(1) } \subseteq \bar R$ for any $h\in G.$

  {\rm (ii)} $\bar R$ is $kG$-subcomodule
 of $R$. We use induction on $n$ to show $\delta ^- (x) \in kG \otimes
 \bar R$. When $n=1$, it is clear. Assume $n>1.$ $\delta ^- (x) = \delta ^- (yz) = \sum
 y _{(-1)}  z _{(-1)} \otimes y_{(0)}
 z_{(0)}$ $\in kG \otimes \bar R$.

{\rm  (iii)} $\bar R$ is a subcoalgebra of $R$. We use induction on $n$ to
 show $\Delta (x) \in \bar R\otimes \bar R.$ When $n=1$ it is
 clear. Assume $n>1.$
  \begin {eqnarray*}\Delta _R (x) &=& \Delta _R(yz)\\
 &=& \sum_{(z)} y z_{1} \otimes
 z_2 +\sum_{(z), (y)} y_{(-1)}\cdot z_{1} \otimes y_{(0)} z_2,
\end {eqnarray*} which implies $\Delta _R(x) \in \bar R \otimes \bar R.$

 For the second claim, since  $R$ is  strictly
 graded as coalgebra gradations (see \cite  [Lemma 2.3 and Lemma  2.4]{AS98}),
   $\bar R$ is a Nichols algebra by the first claim. $\Box$

\begin {Remark} \label {1.5'} By \cite [Cor.2.3 ]{AS02}  $\bar R\cong \mathfrak B ( {\rm diag _{filt} } (H) _{(1)}    )$  as graded braided Hopf algebra
in $^{kG}_{kG} {\mathcal YD}$, where $\bar R$ is the same as in Lemma \ref {1.4}. There exists an ${\rm RSR} (G, r,
\overrightarrow{\rho}, u)$
such that $  {\mathfrak B}(G, r, \overrightarrow{\rho}, u)         \cong {\mathfrak B} ( {\rm diag _{filt} } (H) _{(1)}    )$  as graded braided Hopf algebra
in $^{kG}_{kG} {\mathcal YD}$,  by \cite [Pro. 2.4] {ZCZ08}. We call ${\mathfrak B} ( {\rm diag _{filt} } (H) _{(1)}    ) $ and ${\rm RSR} (G, r,
\overrightarrow{\rho}, u)$ the Nichols algebra and ${\rm RSR}$  of $H$, respectively.
\end {Remark}

\begin {Definition}\label {1.5}
{\rm (i)} ${\rm RSR} (G, r, \overrightarrow{\rho}, u)$ is of
 $-1$-type, if $u(C) $ is real and  the  order of $u(C)$ is  even
 with
  $\chi _C^{(i)}(u(C)) = -\chi _C^{(i)}(1)$ (i.e.
  $\chi _C^{(i)} (u(C)) = - {\rm deg }\rho _C^{(i)})$ for any  $C \in {\mathcal K}_r (G)$ and
  any $i\in I_C(r, u)$.

 {\rm  (ii)} Nichols algebra $R$ over group $G$ is of $-1$-type if there exists $-1$-type ${\rm RSR} (G, r, \overrightarrow{\rho},
  u)$ such that $R \cong \mathfrak {B}(G, r, \overrightarrow{\rho}, u)$
  as  graded pull-push {\rm YD} Hopf algebras.

 {\rm  (iii)} Pointed Hopf algebra $H$ with group $G= G(H)$ is of $-1$-type if the Nichols algebra of $H$ is of $-1$-type.
\end {Definition}

\begin {Proposition} \label {1.6}
{\rm (i)} If ${\rm RSR} (G, r, \overrightarrow{\rho}, u)\cong {\rm RSR}
(G', r', \overrightarrow{\rho'}, u')$  and ${\rm RSR} (G, r,
\overrightarrow{\rho}, u)$ is of
 $-1$-type, then so is ${\rm RSR} (G', r', \overrightarrow{\rho'},
 u')$.

 {\rm  (ii)} If   $(kG, R) \cong (kG', R')$ as   graded pull-push  {\rm YD} Hopf algebras
 and $R$ is
  of $-1$-type,  then so is $R'$,  where $R$  and $R'$  are
   Nichols algebras  over  group algebras $kG$ and $kG'$, respectively.

  {\rm (iii)} If pointed Hopf algebras $H$ and  $H'$
  are isomorphic as Hopf algebras and $H$ is of  $-1$-type, then so is
  $H'$.
\end {Proposition}
{\bf Proof.} {\rm (i)} There exist a group isomorphism $\phi: G\rightarrow
G'$,  an element $h_C\in G$ such that
$\phi(h^{-1}_Cu(C)h_C)=u'(\phi(C))$ for any $C\in{\mathcal K}(G)$
and a bijective map $\phi_C: I_C(r, u)\rightarrow I_{\phi (C)}(r',
u')$ such that $ \rho_C^{(i)}
 {\cong}  \rho'{}_{\phi(C)}^{(\phi_C(i))} \phi _{h_c}$  for any $i\in I_C(r, u)$.
 Therefore
 \begin {eqnarray*}\chi'{} _{\phi (C)}^{(\phi _C(i))} (u'(\phi (C))) &=&
 \chi'{} _{\phi (C)}^{(\phi _C(i))} (\phi (h_C^{-1}u(C)h_C))\\
&=&\chi _C^{(i)} (u(C))  \ \ (\hbox {by the isomorphism })\\
&=& - \chi _C^{(i)} (1) \ \ (\hbox {by the definition of } -1\hbox {-type})\\
&=& -\chi'{} _{\phi (C)}^{(\phi _C(i))} (\phi_{h_C} (1)) = -\chi'{}
_{\phi (C)}^{(\phi _C(i))} (1),\\ \end {eqnarray*}
 which proves the claim.

{\rm (ii)} By \cite [Pro.2.4]{ZCZ08}, there exist two ${\rm RSR} (G, r,
\overrightarrow{\rho}, u)$ and ${\rm RSR} (G', r',
\overrightarrow{\rho'}, u')$ such that $R \cong \mathfrak {B} (G, r,
\overrightarrow{\rho}, u)$ and $R' \cong \mathfrak {B} (G', r',
\overrightarrow{\rho'}, u')$ as graded {\rm YD} Hopf algebras. Thus
${\rm RSR} (G, r, \overrightarrow{\rho}, u)\cong {\rm RSR} (G', r',
\overrightarrow{\rho'}, u')$ by \cite [Theorem 4]{ZCZ08}. It follows
from Definition \ref {1.5} and Part (i) that ${\rm RSR} (G', r',
\overrightarrow{\rho'}, u')$ is of $-1$-type.

{\rm (iii)} It is clear that ${\rm gr }H \cong  {\rm gr }H'$ as graded
Hopf algebras. Thus $(kG, R) \cong (kG',  R')$ as graded pull-push {\rm YD }
Hopf algebras by \cite [Lemma 3.1] {ZCZ08}, where $kG $ and $kG'$ are the coradicals of $H$ and $H'$,  respectively;
 $R= {\rm diag } H$ and $R'= {\rm
diag } H'$. Let $\bar R$ and $\bar {
R'}$ denote the subalgebras generated by $R_{(1)}$ and $R'_{(1)}$ as
algebras in $R$ and $R'$, respectively. It is clear that $( kG, \bar R)
\cong (kG',  \bar {R'})$ as  graded pull-push {\rm YD } Hopf algebras. It follows from Part
{\rm (ii)} that $H'$ is of $-1$-type.
$\Box$

In fact, the proof of Part {\rm (iii)} above shows that if two pointed Hopf algebras are isomorphic, then their Nichols
 algebras are
graded pull-push isomorphic. Similarly, we can prove that their ${\rm RSR}'s$ are isomorphic.

\begin {Proposition} \label {1.7}
If $H$ is a pointed Hopf algebra with real $G= G(H)$ and is not of
$-1$-type, then $H$ is infinite dimensional.
\end {Proposition}
{\bf Proof.} Let $R$ be the (filter) diagram of $H$. By Lemma \ref
{1.4}, $\bar R$ generated by $R_{(1)}$ as algebras
 in $R$ is a Nichols algebra.
By \cite [Pro.2.4 (ii)]{ZCZ08}, there exists an ${\rm RSR} (G, r,
\overrightarrow{\rho}, u)$ such that $\bar R \cong \mathfrak{B} (G,
r, \overrightarrow{\rho }, u)$ is graded pull-push {\rm YD} Hopf
algebra isomorphism. By assumption, there exist $C \in {\mathcal
K}_r(G)$ and $i \in I_C(r, u)$ such that $\chi _C^{(i)} (u(C)) \not=
- {\rm deg} (\rho _C^{(i)})$ or the order of $u(C)$ is odd. It
follows from Lemma \ref {1.2} that the bi-one Nichols algebra $
\mathfrak{B} (G, r', \overrightarrow{\rho'}, u')$ is infinite
dimensional,  where ramification $r' = r'_C C$, $\rho' {} ^{(i)}_C =
\rho _C^{(i)}$, $u'(C) = u(C)$, $I_C(r', u') \subseteq I_C(r, u)$
with $\mid \!I_C(r',u')\!\mid =1 $. Let $Q'$ be a sub-quiver of $Q$ with $Q'_0 = Q_0 $ and $Q'_1:= \{a_{y, x}^{(i, j)}
\mid x^{-1}y \in C, j \in J_C(i)\}.$  Since $ ( k (Q')_1^1, ad (G,
r', \overrightarrow{\rho'}, u'))$ is a braided subspace of  $ ( k
Q_1^1, ad (G, r, \overrightarrow{\rho}, u))$, we have  $ {\rm dim }
\mathfrak{B} ( G, r, \overrightarrow{\rho}, u) = \infty$  and $H$
is  infinite dimensional. $\Box$

 ${\rm RSR}( G, r,
\overrightarrow{\rho}, u)$ is said to be of   {\it infinite type } if $\mathfrak
{B}( G, r, \overrightarrow{\rho}, u)$ is  infinite dimensional. Otherwise, it is said to be of
  {\it finite type }.
For
any ${\rm RSR}( G, r, \overrightarrow{\rho}, u)$, according to the
proof above, if there exist $C \in {\mathcal K}_r(G)$ and $i\in
I_C(r, u)$ such that ${\rm dim}\mathfrak{B} ( {\mathcal O}_{u(C)},
\rho _C^{(i)}) = \infty$, then ${\rm dim} \mathfrak{B} ( G, r,
\overrightarrow{\rho}, u) = \infty$. In this case ${\rm RSR}( G, r,
\overrightarrow{\rho}, u)$ is said to be of  {\it essentially  infinite type}. Otherwise,
it is said to be   of  {\it non-essentially  infinite type}.
For example, non $-1$-type {\rm RSR} over real group is of  essentially
infinite type. However, it is an  open problem  whether  ${\rm RSR}( G, r,
\overrightarrow{\rho}, u)$ is of  finite type when
it is of  non-essentially infinite type, although paper \cite {AHS08} gave
a partial solution to this problem.

\section {Generalized  quantum linear spaces  } \label {s3}
In this section it is shown that every central quantum linear space
 is  finite dimensional with an  arrow {\rm PBW} basis.

Let $\sigma$ denote the braiding of  the   braided tensor category $({\mathcal
C}, \sigma )$. If $A$ and $ B$ are two objects of ${\mathcal C}$ and $\sigma _{A, B} \sigma _{B, A} =
{\rm id} _{B\otimes A}$ and $\sigma _{B, A} \sigma _{A, B} ={\rm  id} _{A\otimes B}$
then $\sigma $ is said to be
{\it symmetric } on pair $(A, B)$. Furthermore, if $A=B$,
then $\sigma$ is said to be  symmetric  on  object $A$, in short,  or $A$ is said to be
 quantum symmetric.


Every arrow {\rm YD} module $(kQ_1^1, {\rm ad } (G, r, \overrightarrow {\rho}, u))$  has
a decomposition of simple {\rm YD} modules: \begin {eqnarray} \label {e3.1}
(kQ_1^1, {\rm ad } (G, r, \overrightarrow {\rho}, u)) &=&
\oplus _{C\in {\mathcal K}_r(G), i \in I_C(r, u)} kQ_1^1 (G, {\mathcal O}_{u(C)},
\rho _C^{(i)})
\end {eqnarray}
 and $\sigma _{C^{(i)}, D^{(j)}}$ is a map
from $kQ_1^1(G, {\mathcal O}_{u(C)}, \rho _C^{(i)}) \otimes
kQ_1^1(G, {\mathcal O}_{u(D)}, \rho _D^{(j)})$ to $kQ_1^1(G, {\mathcal O}_{u(D)}, \rho _D^{(j)}) \otimes
kQ^1_1(G, {\mathcal O}_{u(C)}, \rho _C^{(i)})$,
where  $kQ_1^1(G, {\mathcal O} _{u(C)}, \rho _C^{(i)}) :=
k\{ a_{x, 1}^{(i, j)} \mid x\in {\mathcal O}_{u(C)}, i \in I_C(r, u), j\in J_C(i) \}$ and
$\sigma _{C^{(i)}, D^{(j)}}$
denotes   $\sigma_
{kQ_1^1(G, {\mathcal O}_{u(C)}, \rho _C^{(i)}),
kQ_1^1(G, {\mathcal O}_{u(D)}, \rho _D^{(j)})}$  for any $i\in I_C(r, u)$, $j\in I_D(r, u)$.

Every {\rm YD} module over $kG$ has a decomposition (\ref {e3.1}) since every
{\rm YD} module is isomorphic to an arrow {\rm YD} module by \cite [Pro. 2.4]
{ZCZ08}, which shows every {\rm YD} module over $kG$ is completely reducible (see
\cite [Section 1.2] {AZ07}).

\begin {Definition}\label {3.1}
An $ {\rm RSR}(G, r, \overrightarrow {\rho}, u) $ is said to be quantum
symmetric if $\sigma _{C^{(i)}, D^{(j)}} = \sigma _{D^{(j)}, C^{(i)}})^{-1}$,
 i.e.
$\sigma $ is symmetric on pair  $(kQ_1^1(G, {\mathcal O}_{u(C)}, \rho _C^{(i)}),
kQ_1^1(G, {\mathcal O}_{u(D)}, \rho _D^{(j)})$,
for any $C, D \in {\mathcal K}_r(G)$,
$i \in I_C(r, u)$ and $j \in I_D(r, u)$.

An $ {\rm RSR}(G, r, \overrightarrow {\rho}, u) $ is said to be quantum weakly
symmetric  if $\sigma _{C^{(i)}, D^{(j)}} = (\sigma _{D^{(j)}, C^{(i)}})^{-1}$ for any $C, D \in {\mathcal K}_r(G)$,
$i \in I_C(r, u)$ and $j \in I_D(r, u)$ with  $(C, i) \not= (D, j)$
(i.e. either $C\not= D $ or $i \not= j$).
\end {Definition}

\begin {Proposition}\label {3.2}
If a non-essentially infinite
${\rm RSR} (G, r,  \overrightarrow \rho,u)$  is quantum weakly symmetric,  then ${\rm RSR} (G, r,  \overrightarrow \rho,u)$  is a finite type.
\end {Proposition}
{\bf Proof.} It follows from \cite [Theorem 2.2]{Gr00}. $\Box$

\begin {Lemma} \label {3.3}
{\rm (i)} Assume that $H$ is a Hopf algebra with an invertible antipode  and
$M$ is a {\rm YD} $H$-module, Then   the braiding $\sigma $ of $^H_H {\mathcal YD}$
is symmetric on $M$ if and only if $\sigma $ is symmetric on $\mathfrak {B}(M)$.

{\rm (ii)} The following conditions are equivalent.

{\rm (1)} ${\rm RSR} (G, r, \overrightarrow {\rho}, u)$
is quantum symmetric.

{\rm (2)}  The braiding  $\sigma$ of $^{kG}_{kG} {\mathcal YD}$
on the arrow {\rm YD} module
$(kQ_1^1, {\rm ad} (G, r, \overrightarrow {\rho}, u))$  is symmetric.


{\rm (3)} The braiding $\sigma $ is symmetric on ${\mathfrak B   } (G, r, \overrightarrow {\rho}, u)$.

{\rm (4)} $ \sigma ^2 (a_{x, 1}^{(i, j)} \otimes a_{y,1}^{(i', j')})= a_{x, 1}^{(i, j)} \otimes
a_{y,1}^{(i', j')}  $
for any $C:= x^G, D:=y^G \in {\mathcal K}_r (G)$, $i\in I_C(r, u)$, $i'\in I_D(r, u)$,
$j \in J_C(i)$, $ j'\in
J_D(i')$.

{\rm (5)}     $a_{x, 1}^{(i, j)} \otimes
a_{y,1}^{(i', j')}= (xyx^{-1} \rhd a_{x, 1}^{(i, j)}
) \otimes ( x \rhd a_{y,1}^{(i', j') })  $
for any $C:= x^G, D:=y^G \in {\mathcal K}_r (G)$, $i\in I_C(r, u)$, $i'\in I_D(r, u)$,
$j \in J_C(i)$, $ j'\in
J_D(i')$.

{\rm (iii)}  The following conditions are equivalent.

{\rm (1)} ${\rm RSR} (G, r, \overrightarrow {\rho}, u)$
is quantum weakly symmetric.


{\rm (2)} $ \sigma ^2 (a_{x, 1}^{(i, j)} \otimes a_{y,1}^{(i', j')})= a_{x, 1}^{(i, j)} \otimes
a_{y,1}^{(i', j')}  $
for any $C:= x^G, D:=y^G \in {\mathcal K}_r (G)$, $i\in I_C(r, u)$, $i'\in I_D(r, u)$,
$j \in J_C(i)$, $ j'\in
J_D(i')$ with $(C, i) \not= (D, i')$.

{\rm (3)}     $a_{x, 1}^{(i, j)} \otimes
a_{y,1}^{(i', j')}= (xyx^{-1} \rhd a_{x, 1}^{(i, j)}
) \otimes ( x \rhd a_{y,1}^{(i', j') })  $
for any $C:= x^G, D:=y^G \in {\mathcal K}_r (G)$, $i\in I_C(r, u)$, $i'\in I_D(r, u)$,
$j \in J_C(i)$, $ j'\in
J_D(i')$  with $(C, i) \not= (D, i')$.

 \end {Lemma}

{\bf Proof.} {\rm (i)} It is clear since $M$ generates $\mathfrak {B} (M)$ as algebras.

{\rm (ii)} It follows from  Definition \ref {3.1} that {\rm (1)} and {\rm (2)} are equivalent. Part (i)
 implies that {\rm (2)} and {\rm (3)} are equivalent. Obviously, {\rm (4)} and {\rm (2)} are equivalent.
Since $(kQ_1^1, {\rm ad} (G, r, \overrightarrow \rho, u))$ is a {\rm YD} module, we have
\begin {eqnarray} \label {e3.2} \sigma ^2 (a_{x, 1}^{(i, j)} \otimes a_{y,1}^{(i', j')})= (xyx^{-1} \rhd a_{x, 1}^{(i, j)}
) \otimes ( x \rhd a_{y,1}^{(i', j') }).
\end {eqnarray}
Therefore, {\rm (4)} and {\rm (5)} are equivalent.

{\rm (iii)} It follows from  (\ref {e3.2}) that {\rm (2)} and {\rm (3)} are equivalent. Obviously
 {\rm (1)} and {\rm (2)}  according to the definition.
$\Box$






\begin {Lemma}\label {3.4}

For  $C:= x^G, D:=y^G \in {\mathcal K}_r (G)$, $i\in I_C(r, u)$, $i'\in I_D(r, u)$,
$j \in J_C(i)$, $ j'\in
J_D(i')$,
 assume that $\rho _C^{(i)}$ and $\rho _D^{(i')}$ are one dimensional representations;
The coset decomposition of $G^{u(C)}$  and $G^{u(D)}$ in $G$ are
\begin {eqnarray*}
G &=&\bigcup_{\theta\in\Theta_C}G^{u(C)}g_{\theta}, \hbox {and }
G = \bigcup_{\eta \in\Theta_D}G^{u(D)}h_{\eta},
\end {eqnarray*}  respectively;
$x = g_\theta ^{-1}u(C)g_\theta $ and $y = h_\eta  ^{-1}u(D)h_\eta $;
$g_\theta y^{-1} = \zeta _{\theta} (y^{-1}) g_{\theta'}$ and $h_\eta x^{-1} =
\zeta _\eta (x^{-1}) h_{\eta'}$
with $\zeta _{\theta} (y^{-1})\in G^{u(C)}$ and  $\zeta _\eta (x^{-1}) \in G^{u(D)}$.

Then  \begin {eqnarray} \label {e3.3}
a_{x, 1}^{(i, j)} \otimes
a_{y,1}^{(i', j')}= (xyx^{-1} \rhd a_{x, 1}^{(i, j)}
) \otimes ( x \rhd a_{y,1}^{(i', j') })
\end {eqnarray} if and only if
\begin {eqnarray} \label {e3.4} xy=yx  \ \  \hbox { and } \ \
\rho _C^{(i)} ( \zeta _\theta (y^{-1})) \rho _D^{(i')} ( \zeta _\eta (x^{-1})) &=&1
\end {eqnarray}

\end {Lemma}
{\bf Proof.} By \cite [Pro. 1.2] {ZCZ08} or \cite [Pro. 1.9] {ZZC04}, $(xyx^{-1} \rhd a_{x, 1}^{(i, j)}
) \otimes ( x \rhd a_{y,1}^{(i', j') }) = \alpha  a ^{(i, j)} _{xy xy^{-1}x^{-1}, 1   }
\otimes
a ^{(i',j')}_{ x yx^{-1}, 1} $, where $\alpha \in k$. Thus
(\ref {e3.3}) holds if and only if $xy=yx$ and $\alpha = 1$. By  \cite [Pro. 1.2] {ZCZ08},
$\alpha =1$ if and only if (\ref {e3.4}) holds. $\Box$


\begin {Proposition}\label {3.5} If
${\rm RSC}(G, r, \overrightarrow \chi, u)$ is  non-essentially infinite and
 (\ref {e3.4}) holds for any  $C:= x^G,
D:=y^G \in {\mathcal K}_r (G)$, $i\in I_C(r, u)$, $i'\in I_D(r, u)$ with $(C, i) \not= (D, i')$,
then
${\rm RSR} (G, r,  \overrightarrow \rho,u)$  is quantum weakly symmetric. Therefore
${\rm RSR} (G, r,  \overrightarrow \rho,u)$  is a finite type.
\end {Proposition}
{\bf Proof.} It follows from Lemma \ref {3.4}, Lemma \ref {3.3} and Proposition \ref {3.2}. $\Box$

If $0 \not= q \in
k$ and $0 \le i \le n < ord (q) $ (the order of $q$), we set
$(0)_{q}! =1$,
$$ \left ( \begin {array} {c} n\\
i
\end {array} \right )_q
 =  \frac{(n)_{q}!}{(i)_{q}!(n - i)_{q}!}, \quad \hbox {where
}(n)_{q}! = \prod_{1 \le i \le n} (i)_{q}, \quad (n)_{q} =
\frac{q^{n} - 1}{q - 1}.$$ In particular, $(n)_q = n$ when $q=1.$

\begin {Lemma}\label{3.6} In $kQ^c(G, r, \overrightarrow \rho, u)$, we have the following results.

{\rm (i)} If $ C= \{g \} \in {\mathcal K}_r(G)$ with $i\in I_C (r, u)$, then there exists
 $0\not =q\in k $ such that $\rho_C^{(i)} (g) = q \ {\rm id} $
and $a_{y,x} ^{(i, j)} \cdot h = q a_{yh, xh} ^{(i, j)} $ for any $x^{-1}y \in C$, $h\in G$,
 $j \in J_C(i)$.

{\rm (ii)} If $a_{w_0,v_0} ^{(i, j)} \cdot h = q a_{w_0h, v_0h} ^{(i, j)}$ for some $v_0, w_0,
 \in G$, $h \in G^{u(C)}$,  $q\in k$ with $v_0^{-1}w_0 \in C \in {\mathcal K}_r(G),$
then  $a_{w,v} ^{(i, j)} \cdot h  = q a_{w_0h, v_0h} ^{(i, j)}$ for any $v, w \in G$ with
$v^{-1}w \in C$
\end {Lemma}

{\bf Proof.} {\rm (i)} It follows from \cite [Pro. 1.2]{ZCZ08}.

{\rm (ii)}  Let $X_C^{(i)}$ be a representation space of $\rho _C^{(i)}$ and
$\{{x_C^{(i,j)}\mid j \in J_C(i)}\}$  a k-basis of $X_C^{(i)}$.
By the proof of \cite [Pro. 1.2] {ZCZ08},
 $a_{w, v}^{(i,j)} \cdot
h =  q   a_{wh, vh}^{(i,j)}$ since
$x_C^{(i,j)} \cdot \zeta_{\theta}(h) = q x_C^{(i, j)}$ by assumption. $\Box$

\begin {Lemma}\label{3.7} In co-path Hopf algebra $kQ^c(G, r, \overrightarrow \rho, u)$, assume $C:= g^G \in
{\mathcal K}_r(G)$,
$i\in I_C(r, u)$, $j \in J_C(i)$ and $ a_{g, 1}^{(i, j)} \cdot g = q  a_{g^2, g}^{(i, j)} $.
 If $i_1, i_2, \cdots, i_m$ are non-negative integers, then
$$\begin{array}{rcl}
a^{(i, j)}_{g^{i_m +1},g^{i_m }}\cdot a^{(i, j)}_{g^{i_{m-1}
+1},g^{i_{m-1}}}\cdot\cdots\cdot a^{(i, j)}_{g^{i_1 +1},g^{i_1 }}
&=&q^{\beta_m}(m)_q! P^{(i, j)}_{g^{\alpha_m}}(g,m)\\
\end{array}$$
where $\alpha _m = i_1 + i_2 + \cdots + i_m $, $P^{(i, j)}_h(g,m) :=$ \
$ a^{(i, j)}_{g^mh,g^{m-1}h}a^{(i, j)}_{g^{m-1}h, g^{m-2}h}\cdots
a^{(i, j)}_{gh,h}$, $\beta_1=0$ and $\beta_m=\sum_{j
=1}^{m-1}(i_1+i_2+\cdots+i_j )$ if $m>1$.
\end{Lemma}
{\bf Proof.} We prove the  equality by induction on $m$. For $m=1$,
it is easy to see that the equality holds. Now suppose $m>1$. We
have

$$\begin{array}{rl}
&a^{(i, j)}_{g^{i_m +1},g^{i_m }}\cdot a^{(i, j)}_{g^{i_{m-1}
+1},g^{i_{m-1}}}\cdot\cdots\cdot a^{(i, j)}_{g^{i_1 +1},g^{i_1 }}\\
=&a^{(i, j)}_{g^{i_m +1},g^{i_m }}\cdot(a^{(i, j)}_{g^{i_{m-1}
+1},g^{i_{m-1}}}\cdot\cdots\cdot a^{(i, j)}_{g^{i_1 +1},g^{i_1 }})\\
=&q^{\beta_{m-1}}(m-1)_q! a^{(i, j)}_{g^{i_m +1},g^{i_m }}\cdot
P^{(i, j)}_{g^{\alpha_{m-1}}}(g,m-1) \ \ \ ( \hbox {by inductive
assumption })\\
=&q^{\beta_{m-1}}(m-1)_q! a^{(i, j)}_{g^{i_m +1},g^{i_m}}\cdot
(a^{(i, j)}_{g^{\alpha_{m-1}+m-1},g^{\alpha_{m-1}+m-2}}\cdots
a^{(i, j)}_{g^{\alpha_{m-1}+1},g^{\alpha_{m-1}}}) \\
 =&q^{\beta_{m-1}}(m-1)_q! \sum_{l=1}^m[(g^{i_m+1}\cdot
a^{(i, j)}_{g^{\alpha_{m-1}+m-1},g^{\alpha_{m-1}+m-2}})
\cdots(g^{i_m+1}\cdot a^{(i, j)}_{g^{\alpha_{m-1}+l},g^{\alpha_{m-1}+l-1}})\\
&(a^{(i, j)}_{g^{i_m+1},g^{i_m}}\cdot g^{\alpha_{m-1}+l-1})
(g^{i_m}\cdot
a^{(i, j)}_{g^{\alpha_{m-1}+l-1},g^{\alpha_{m-1}+l-2}})\cdots
(g^{i_m}\cdot a^{(i, j)}_{g^{\alpha_{m-1}+1},g^{\alpha_{m-1}}})]\\
&\ \ \ ( \hbox {by \cite[Theorem 3.8]{CR02} })
\\
=&q^{\beta_{m-1}}(m-1)_q!
\sum_{l=1}^m[a^{(i, j)}_{g^{\alpha_m+m},g^{\alpha_m+m-1}}
\cdots a^{(i, j)}_{g^{\alpha_m+l+1},g^{\alpha_m+l}}\\
& q^{\alpha_{m-1}+l-1}a^{(i, j)}_{g^{\alpha_m+l},g^{\alpha_m+l-1}}
a^{(i, j)}_{g^{\alpha_m+l-1},g^{\alpha_m+l-2}}\cdots
a^{(i, j)}_{g^{\alpha_m+1},g^{\alpha_m}}]    \ \  ( \hbox  { by     lemma \ref {3.6} } )
 \\
=&q^{\beta_{m-1}}(m-1)_q! \sum_{l=1}^m
q^{\alpha_{m-1}+l-1}P^{(i, j)}_{g^{\alpha_m }}(g,m)\\
=&q^{\beta_{m-1}+\alpha_{m-1}}(m)_q!P^{(i, j)}_{g^{\alpha_m }}(g,m)\\
=&q^{\beta_m}(m)_q! P^{(i, j)}_{g^{\alpha_m }}(g,m). \ \ \Box
\end{array}$$

Recall that a braided algebra $A$ in braided tensor category
$({\mathcal C}, \sigma )$ with braiding $\sigma$ is said to be
braided commutative or quantum commutative, if $ab = \mu \sigma( a
\otimes b)$ for any $a, b\in A$, where $\mu$ is the multiplication
of $A$.

By \cite  [Example 3.11]{CR02},  the multiplication of
any two arrows $a^{(i, j)}_{y,x}$ and $ a^{( m, n)}_{{w,v, }}$ in co-path Hopf algebra
$kQ^c (G, r, \overrightarrow \rho, u)$
 is
\begin{eqnarray}\label {multiplication}
a^{(i, j)}_{y,x} \cdot a^{( m, n)}_{{w,v, }} &=&
(y\cdot a^{(m, n)}_{{w, v}}) (a_{y, x}^{(i, j)} \cdot v)
+ ( a_{y,x} \cdot w) (x \cdot a_{w, v}^{(m, n)}).
\end{eqnarray}

\begin {Lemma} \label {3.8} Let $C:= x^G$, $D:= y^G$ $\in {\mathcal K}_r(G)$,
$i\in I_C(r, u)$, $j\in J_C(i)$, $i'\in I_D(r, u)$, $j'\in J_D(i')$, $\alpha, \beta
\in k $ with $a_{y,1}^{(i', j')} \cdot x = \alpha a_{yx, x}^{(i', j')}$ and $
 a_{x, 1}^{(i, j)} \cdot y = \beta a_{xy, y}^{(i, j)}$ in co-path Hopf
algebra $kQ^c (G, r, \overrightarrow \rho, u)$. If $xy = yx$ then
 $ \alpha \beta = 1$  if and only if
 \begin {eqnarray}\label {e3.6}  a_{x, 1}^{(i, j)}\cdot a_{y, 1} ^{(i', j')} =
\alpha ^{-1}a_{y,1}^{(i', j')}\cdot a_{x,1}^{(i, j)} \end {eqnarray}

\end {Lemma}

{\bf Proof.}
 By (\ref {multiplication}) and  \cite [Pro. 1.2] {ZCZ08}, we have
 \begin{eqnarray} \label {e3.10}
a_{x,1}^{(i, j)}\cdot a_{y,1}^{(i', j')}&=& a^{(i', j')}_{xy,x}a^{(i,j)}_{x,1}+
 \beta  a^{(i,j)}_{xy,y}a^{(i', j')}_{y,1}, \nonumber\\
a_{y,1}^{(i', j')}\cdot a_{x,1}^{(i, j)}& =&
 \alpha  a^{(i',j')}_{yx,x}a^{(i,j)}_{x,1}+
a^{(i,j)}_{yx,y}a^{(i',j')}_{y,1}.
\end{eqnarray} Applying this we can complete the proof. $\Box$

\begin {Lemma} \label {3.9} Assume that $ {\rm RSR}(G, r, \overrightarrow \rho, u) $ satisfies
$C:= \{g_C\}\subseteq Z(G)$ for any $C\in {\mathcal K}_r(G)$. Let  $\rho _C^{(i)} (g_D) =
q_{C, D}^{(i)} \ {\rm id}$ for any $C, D\in {\mathcal K}_r(G)$, $i \in I_C(r, u)$.

{\rm (i)} The following conditions are equivalent:

{\rm (1)} $ {\rm RSR}(G, r, \overrightarrow \rho, u) $ is   quantum
symmetric

{\rm (2)} $  q_{C, D}^{(i)} q_{D, C}^{(i')}=1$ for any $C, D\in {\mathcal K}_r(G)$, $i\in I_C(r, u)$,
$i' \in I_D(r, u)$.

{\rm (3)} $ a_{g_C, 1}^{(i, j)}\cdot a_{g_D, 1} ^{i', j'} =
(q _{D,  C}^{(i')})^{-1}  a_{g_D,1}^{(i', j')}\cdot a_{g_C,1}^{(i, j)}$
for any $C, D\in {\mathcal K}_r(G)$, $i\in I_C(r, u)$, $i' \in I_D(r, u)$,  $j \in I_C(i)$,
$j' \in J_C(i')$.

{\rm (4)} $\mathfrak {B} (G, r, \overrightarrow \rho, u)$ is quantum
commutative in $^{kG}_{kG} {\mathcal YD}$.

{\rm (5)} $\mathfrak {B} (G, r, \overrightarrow \rho, u)$ is quatntum symmetric.

{\rm (6)} $(kQ_1^1, {\rm ad}  (G, r, \overrightarrow \rho, u)  )$ is quantum symmetric.

{\rm (ii)} The following conditions are equivalent:

{\rm (1)} $ {\rm RSR}(G, r, \overrightarrow \rho, u) $ is   quantum weakly
symmetric

{\rm (2)} $  q_{C, D}^{(i)} q_{D, C}^{(i')}=1$ for any $C, D\in {\mathcal K}_r(G)$,
$i \in I_C(r, u)$, $i' \in I_D(r, u)$ with $ (C, i) \not = (D, i')$.

{\rm (3)} $ a_{g_C, 1}^{(i, j)}\cdot a_{g_D, 1} ^{(i', j')} =
(q _{D,  C}^{(i')} )^{-1} a_{g_D,1}^{(i', j')}\cdot a_{g_C,1}^{(i, j)} $
for any $C, D\in {\mathcal K}_r(G)$, $i\in I_C(r, u)$, $i'\in I_D(r, u)$,  $j \in I_C(i)$,
$j' \in J_C(i')$ with $ (C, i) \not = (D, i')$.

\end {Lemma}

{\bf Proof.}   By \cite [Lemma 2,2] {ZCZ08} , $ diag (kG [kQ_1^c, r, \overrightarrow \rho, u
])$ is the Nichols algebra $\mathfrak {B} (G, r,
 \overrightarrow \rho, u)$
in $^{kG}_{kG}
{\mathcal YD}$.  By \cite [Pro. 1.2] {ZCZ08},
\begin{eqnarray} \label {e3.11}
\sigma ^2(a_{g_C,1}^{(i, j)}\otimes  a_{g_D,1}^{(i', j')})&=&
(q_{C, D}^{(i)}q_{D, C}^{(i')} )^{-1} a_{g_C,1}^{(i, j)}\otimes  a_{g_D,1}^{(i', j')}.
\end{eqnarray}

{\rm (i)} By Lemma \ref {3.3} {\rm (ii)}, {\rm (1), (5)} and {\rm (6)} are equivalent. It follows from (\ref {e3.11})
that  {\rm (6)} and {\rm (2)} are equivalent. By Lemma \ref {3.8}, {\rm (3)} and {\rm (2)} are equivalent.
Obviously {\rm (3)} and {\rm (6)} are equivalent. {\rm (3)} and {\rm (4)} are equivalent since
$\mathfrak {B} (G, r, \overrightarrow \rho, u)$ is generated by $kQ_1^1$.

{\rm (ii)} It follows from (\ref {e3.11}) that {\rm (1)} and {(2)} are equivalent.
{\rm (2)} and {(3)} are equivalent according to (\ref {e3.6}).
 \  $\Box$



\begin {Lemma} \label {3.10} (See \cite [Lemma 3.3]{AS98}) Let $B$ be a Hopf algebra and $R$  a braided Hopf algebra in
${}_B^B {\mathcal YD}$ with a linearly independent set $ \{ x_{1} \dots,
x_{t} \}$ $\subseteq $ $P(R)$, the set of all primitive elements in $R$.  Assume that there exist $g_{j} \in
G(B)$ (the set of all group-like elements in $B$) and  $0\not= k_{j, i} \in k$ such that
 $$\delta (x_{i}) = g_{i}\otimes x_{i}, \ g_i\cdot  x_{j} = k_{ij}x_{j},
 \hbox { for all } i, j =1, 2, \cdots, t . $$
Then \begin{eqnarray*}\{x_1^{m_1} x_2^{m_2}\cdots x_t^{m_t} \mid 0
\leq m_j<N_j, 1\le j \le t \}. \end{eqnarray*} is linearly
independent, where $N_i$ is the order of $q_i := k_{ii}$  \ (
$N_i = \infty $ when $q_i$ is not a root of unit, or $q_i =1$ ) for
$1 \le i \le t.$
\end {Lemma}

{\bf Proof.} By the quantum binomial formula, if $ 1 \le n_j < N_j$,
then
$$\Delta (x_{j}^{n_{j}}) = \sum_{0 \le i_{j} \le n_{j}}
 \left( \begin {array}{c} n_j\\
i_j
\end {array} \right )
_{q_{j}}x_{j}^{i_{j}}\otimes x_{j}^{n_{j}-i_{j}}.$$

 We use
the  following notation:  $${\bf n } = (n_{1}, \cdots, n_{j}, \cdots
, n_{t}), \quad  x^{\bf n}= x_{1}^{n_{1}} \cdots x_{j}^{n_{j}}
\cdots x_{t}^{n_{t}}, \quad \vert {\bf n }\vert = n_{1} + \cdots +
n_{j} + \cdots  + n_{t};$$ accordingly, ${\bf N}= (N_{1}, \cdots,
N_{t})$, ${\bf  1 }= (1,\cdots, 1)$. Also, we set
$${\bf i }\le {\bf n}\quad \hbox {if } i_{j} \le n_{j}, \, j = 1,
\cdots, t.$$ Then, for ${\bf n } <{\bf N }$, we deduce from the
quantum binomial formula that \begin {eqnarray} \label
{e3.511}\Delta (x^{{\bf n}}) = x^{{\bf n}}\otimes 1 + 1\otimes
x^{\bf  n} + \sum_{ 0 \le {\bf i} \le {\bf n} , \  0 \ne {\bf i }\ne
{\bf n}} c_{\bf {n, i}}x^{\bf i}\otimes x^{\bf  n- i},\end
{eqnarray} where $c_{\bf n,  i} \ne 0$ for all ${\bf i}$.

We shall prove   by induction on $r$ that the set
$$\{x^{\bf n} \mid  \quad \vert {\bf n} \vert \le r, \quad  {\bf n}
< {\bf N} \}$$ is linearly independent.

Let $r=1$ and let $a_{0} + \sum_{i=1}^{t} a_{i}x_{i} = 0$, with
$a_{j}\in k$, $0\le j \le t$. Applying $\epsilon$, we see that
$a_{0} = 0$; by hypothesis we conclude that the other $a_{j}$'s are
also 0.

Now let $r> 1$ and  suppose that $z = \sum_{   \vert {\bf n}\vert
\le r, {\bf n} < {\bf N} } a_{\bf  n} x^{\bf n} = 0$. Applying
$\epsilon$, we see that $a_{0} = 0$. Then
\begin {eqnarray*} 0
&=& \Delta (z) =z\otimes 1 + 1\otimes z + \sum_{ 1\le \vert{\bf n
}\vert \le r,  {\bf n} < {\bf N}} a_{\bf n}\sum_{ 0 \le {\bf i } \le
{\bf n}, \ 0 \ne {\bf i} \ne { \bf
n } }  c_{\bf n, i}x^{\bf i}\otimes x^{\bf n-  i} \\
&=& \sum_{1\le \vert {\bf n }\vert \le r,  {\bf n} < {\bf N}}\ \
\sum_{ 0 \le {\bf i} \le {\bf n}, \  0 \ne {\bf i }\ne { \bf n}}
a_{\bf n} c_{\bf n, \bf i}x^{\bf i}\otimes x^{\bf n-\bf i}. \end
{eqnarray*} Now, if $\vert{\bf n}\vert \le r$, $0 \le {\bf i} \le
{\bf n}$, and $0 \ne {\bf i }\ne {\bf n}$, then $\vert{\bf i}\vert <
r$ and $\vert{\bf n }- {\bf i}\vert < r$. By inductive hypothesis,
the elements $x^{\bf i}\otimes x^{\bf n-\bf i}$ are linearly
independent. Hence $a_{\bf n} c_{\bf n, \bf i} = 0$ and $a_{\bf n} =
0$ for all ${\bf n}$, $\vert{ \bf n }\vert \ge 1$.  Thus $a_{\bf n}
= 0$ for all ${\bf n}$. $\Box$

The quantum linear space was defined in \cite [Lemma 3.4] {AS98} and now is generalized as
follows.

\begin {Definition}\label {3.11} Let $0 \not= k_{i, j} \in k$ and $1 <N_i : =
{\rm ord } (k_{k_{i, i}})< \infty$  for any $i, j \in \Omega$, where $\Omega$ is a  finite set.
If $R$ is the algebra generated by
set $ \{ x_j \mid j\in \Omega  \}$  with relations
\begin {eqnarray} \label {qlse1} x_l ^{N_l}=0,\  x_{i}x_{j} = k_{i, j}
x_{j}x_{i}  \ \ \ \hbox { for any } i, j  \in \Omega  \hbox { with }
  i \not= j,\end {eqnarray}
 then $R$ is called the {\it  generalized quatum linear space } generated by
$ \{ x_j \mid j\in \Omega  \}$.
\end {Definition}


\begin {Definition}\label {3.12}

{\rm (i)}  ${\rm RSR}(G, r, \overrightarrow {\rho}, u)$ is said to be  a {\it generalized quantum linear
type } if
 the following conditions  are satisfied:

{\rm (GQL1)} $xy = yx$ for any $C:=x^G, D:=y^G \in {\mathcal K}_r(G)$.

{\rm (GQL2)} there exists $k _{x, y} ^{(i, j)} \in k$ such that
$a_{x, 1} ^{(i, j)} \cdot y
=k _{x, y} ^{(i, j)} a_{xy, y} ^{(i, j)} $ for any $C:=x^G, D:=y^G \in {\mathcal K}_r(G)$,
$i\in I_C(r, u)$, $j\in J_C(i)$.

{\rm (GQL3)}  $k _{x, y} ^{(i, j)} k _{y, x} ^{(i', j')} =1 $ for any
$C:=x^G, D:=y^G \in {\mathcal K}_r(G)$, $i\in I_C(r, u)$, $j\in J_C(i)$,
$i' \in  I_D(r, u)$, $j'\in J_D(i')$ with $(x, i, j) \not= (y, i', j') $.

{\rm (GQL4)}  $1 < N_{x} ^{(i, j)} :={\rm ord }(k_{x, x}^{(i, j)}) < \infty$
 for any
$C:=x^G \in {\mathcal K}_r(G)$, $i\in I_C(r, u)$, $j\in J_C(i)$.

{\rm (ii)} ${\rm RSR}(G, r, \overrightarrow {\rho}, u)$ is said to be  a {\it  central
 quantum linear type}  if it is
 quantum  symmetric
  and  of the non-essentially infinite type with
$C \subseteq Z(G)$ for any $C\in {\mathcal K}_r(G)$. In this case, $\mathfrak {B}
(G, r, \overrightarrow \rho, u)$ is called a {\it central  quantum linear space} over $G$.


\end {Definition}

Assume that  $A$ is an algebra with $\{b_\nu \mid \nu \in \Omega \}\subseteq A$ and $\prec$ is a total order
of $\Omega$, $N_\nu \in {\mathbb N}$ or $\infty$ for any $\nu \in \Omega$.
If \begin{eqnarray}\label{basis}
\{ b_{\nu _1} ^{m_{1}} b_{\nu _2}^{m_{2}} \cdots b_{\nu_n} ^{m_{n}} &
&\mid  \nu _1 \prec \nu _2 , \cdots \prec \nu _n; 0 \leq m_{s} <N_{\nu_s};
1 \leq s \leq n;
 \ \ n\in  {\mathbb N} \}
\end{eqnarray} is a basis of $A$, then the basis ( \ref {basis}) is called a   { \rm PBW } basis
generated by $\{b_\nu \mid \nu \in \Omega \}$. If $\{b_\nu \mid \nu \in \Omega \} \subseteq
Q_1$, then it is called an arrow { \rm PBW } basis.

It is well-known that every quantum linear space  is a braided Hopf algebra and
 has a {\rm BPW} basis
(see \cite [Lemma 3.4] {AS98}). Of course, every  generalized quantum linear space is finite dimensional.
However,
 it is  not known  whether every  generalized quantum linear space has  an {\rm PBW} basis.


\begin {Proposition}\label {3.13}
 If ${\rm RSR}(G, r, \overrightarrow {\rho}, u)$ is of the  generalized quantum linear type, then
 $\mathfrak {B}
(G, $ $ r, $ $ \overrightarrow \rho, u)$ is a generalized quantum linear space with
the arraw
{\rm PBW} basis
\begin{eqnarray}\label{basis1}
\{ b_{\nu _1} ^{m_{1}} b_{\nu _2}^{m_{2}} \cdots b_{\nu_n} ^{m_{n}} &
&\mid  \nu _1 \prec \nu _2 , \cdots \prec \nu _n; 0 \leq m_{s} <N_{\nu_s};
1 \leq s \leq n;
 \ \ n\in  {\mathbb N} \} \end{eqnarray} and
\begin {eqnarray} \label {dim}{\rm dim } (\mathfrak {B}
(G, r, \overrightarrow \rho, u)) &=& \prod
 _{C:= x^G \in {\mathcal K}_r(G),
 i \in I_C(r, u), j\in J_C(i)} N_x^{(i, j)},
\end {eqnarray} where  $\{b_\nu \mid \nu \in \Omega \}
: = Q_1^1$ with total order $\prec$ and  $N_{\nu_s} = N_{x}^{(i, j)}$ $ := {\rm ord }
(k_{x, x}^{(i, j)})$ if $b_{\nu_s} =
a_{x, 1}^{(i, j)}$.
\end {Proposition}

{\bf Proof.}  Since any two different arrows in  $\mathfrak {B}
(G, r, \overrightarrow \rho, u)$ are quantum commutative (see
Lemma \ref {3.8}) and
$(b_{\nu _s} ) ^{N_{\nu _s}} =0$ (see Lemma \ref {3.7} ), we have
 $\mathfrak {B}
(G, r, \overrightarrow \rho, u)$ is generated by (\ref {basis1}).

For any $\nu, \nu' \in \Omega$,  $b_\nu = a_{x, 1}^{(i, j)}$
and $b_{\nu'} = a_{y, 1} ^ {(i', j')}$ with
$C:=x^G, D:=y^G \in {\mathcal K}_r(G)$, $i\in I_C(r, u)$, $j\in J_C(i)$,
$i' \in  I_D(r, u)$, $j'\in J_D(i')$, let $g_\nu = x$ and $k_{\nu, \nu' } = (k_{y, x}^{(i', j')})^{-1}$.
By \cite [Pro. 1.2] {ZCZ08} we have
\begin {eqnarray*}
\delta ^-(b_\nu) &=& \delta ^- (a_{x, 1}^{(i, j)})  =
 x \otimes a_{x, 1} ^{(i, j)} = g_\nu \otimes b_\nu \hbox {\ \ and } \\
g _\nu \rhd b_{\nu' } &= & x \cdot a_{y, 1}^{(i', j')} \cdot x^{-1}\\
&=& x \cdot  (k_{y, x}^{(i', j')})^{-1} a _{yx^{-1}, x^{-1}} ^{(i', j')} \ \ (\hbox
{by ({\rm GQL2})})\\
&=& k _{\nu , \nu'} b_{\nu'} \ \ ( \hbox {by ({\rm GQL1})}).
\end {eqnarray*}
Therefore, by Lemma \ref {3.10}, (\ref {basis1}) is  linearly
independent. Thus (\ref {basis1}) is a
basis of $\mathfrak {B}
(G, r, \overrightarrow \rho, u)$.

Let $R$ is the generalized quantum linear space generated by $\{b_\nu \mid \nu \in
\Omega \} := kQ_1^1$.
It is clear that there exists an algebra map $\psi $ from $R$ to $\mathfrak {B}
(G, r, \overrightarrow \rho, u)$ by sending $b_\nu $ to $b_\nu$
 for any $\nu \in \Omega$. Since $\mathfrak {B}
(G, r, \overrightarrow \rho, u)$ has an arrow  ${\rm PBW}$ basis (\ref {basis2}), $\psi$ is
isomorphic. $\Box$

\begin {Proposition}\label {3.15}  Assume that  $C= \{g _C\} \subseteq Z(G)$ and $\rho _C^{(i)} (g_D) =
q_{C, D}^{(i)} \ {\rm id}$ for any $ C, D \in {\mathcal K}_r(G)$, $i\in I_C(r, u)$.
Then

{\rm (i)} ${\rm RSC}(G, r, \overrightarrow {\chi}, u)$  is of the  central
 quantum linear type
if and only if $q_{C,D}^{(i)} q_{D, C}^{(j)}= 1$ and $1 < {\rm ord } (q_{C, C}^{(i)}) < \infty$
for any $C, D\in {\mathcal K}_r (G)$, $i \in I_C(r, u)$, $j\in I_D(r, u)$.

{\rm (ii)} ${\rm RSC}(G, r, \overrightarrow {\chi}, u)$  is quantum weakly symmetric
 with non-essetially infinite type
if and only if $q_{C,D}^{(i)} q_{D, C}^{(j)}= 1$ and $1 < {\rm ord } (q_{C, C}^{(i)}) < \infty$
for any $C, D\in {\mathcal K}_r (G)$, $i \in I_C(r, u)$, $j\in I_D(r, u)$ with $(C, i)
\not= (D, j).$

\end {Proposition}

{\bf Proof.} {\rm (i)}
If ${\rm RSC}(G, r, \overrightarrow {\chi}, u)$  is of the central  quantum
linear type, then
  ${\rm dim} {\mathfrak B}$ $ ( G, {\mathcal O}_{u(C)}, $ $ \rho _C^{(i)}) $ $< \infty$ for any
$C\in {\mathcal K}_r(G)$, $i \in I_C(r, u)$. Let $N_C^{(i)} := {\rm ord } (q_{C, C}^{(i)})$
(
$N_C^{(i)} = \infty$ when $q_{C, C}^{(i)}$ is not a root of unit or $q_{C, C}^{(i)}=1$ ).
By Lemma \ref {3.10}, $\{
(a_{g_C, 1} ^{(i, j)}) ^{m} \mid  0 \leq m  < N_{C}^{(i)}\}$ is linearly independent.
Thus $1 < {\rm ord } (q_{C, C}^{(i)}) < \infty$. Since
${\rm RSR}(G, r, \overrightarrow {\rho}, u)$  is quantum symmetric,
$q_{C,D}^{(i)} q_{D, C}^{(j)}= 1$ by Lemma \ref {3.9}.

Conversely, by Lemma \ref {3.9}, ${\rm RSR}(G, r, \overrightarrow {\rho}, u)$  is quantum
 symmetric. It is clear that ${\rm RSR}(G, r, \overrightarrow {\rho}, u)$
 is of the  generalized quantum
linear type. Thus it is of the  non-essentially  infinite type by Proposition \ref {3.13}.

{\rm (ii)} It is similar to {\rm (i)}.  $\Box$

 The following is the
consequence of Proposition \ref {3.13} and Proposition {3.15}.

\begin {Proposition}\label {3.14}
If
 ${\rm RSR}(G, r, \overrightarrow \rho, u)$ is of the  central   quantum linear
type,  then  $\mathfrak {B}
(G, r, \overrightarrow \rho, u)$ is a generalized quantum linear space with the arrow
{\rm BPW} basis
\begin{eqnarray}\label{basis2}
\{ b_{\nu _1} ^{m_{1}} b_{\nu _2}^{m_{2}} \cdots b_{\nu_n} ^{m_{n}} &
&\mid  \nu _1 \prec \nu _2 , \cdots \prec \nu _n; 0 \leq m_{s} <N_{\nu_s};
1 \leq s \leq n;
 \ \ n\in  {\mathbb N} \}
\end{eqnarray} and
  \begin {eqnarray} \label {dim}{\rm dim } (\mathfrak {B}
(G, r, \overrightarrow \rho, u)) &=& \prod
 _{C\in {\mathcal K}_r(G),
 i \in I_C(r, u)} (N_C^{(i)}) ^ {{\rm deg} (\rho _C^{(i)}) \mid  C  \mid },
\end {eqnarray} where $\{b_\nu \mid \nu \in \Omega \}
: = Q_1^1$ with total order $\prec$ and $N_{\nu_s} = {\rm ord } (q_{C, C}^{(i)})$ if $b_{\nu_s} =
a_{g_C, 1}^{(i, j)}$.

In particular, if
 ${\rm RSR}(G, r, \overrightarrow \rho, u)$ is quantum weakly commutative and  of  $-1$-type
 with
$C \subseteq Z(G)$ for any $C\in {\mathcal K}_r(G)$, then
it is of the central quantum linear type with  $N_C^{(i)} =2$ and
 \begin {eqnarray} \label {dim}{\rm dim } (\mathfrak {B}
(G, r, \overrightarrow \rho, u)) &=& 2 ^{ \sum _{C\in {\mathcal K}_r(G),
 i \in I_C(r, u)}  {{\rm deg} (\rho _C^{(i)})} \mid C \mid}.
\end {eqnarray}
\end {Proposition}

\begin {Remark} \label {3.16}
 ${\rm RSR}(G, r, \overrightarrow {\rho}, u)$  is called a  {\it central ramification system
with irreducible representations}  (or {\rm CRSR } in short ) if  $C \subseteq Z(G)$ for any
$C\in {\mathcal K}_r(G)$. If $G$ is a real group and $r= r_CC$, then
${\rm CRSR}(G, r, \overrightarrow {\rho}, u)$ is of finite type if and only only if
${\rm CRSR}(G, r, \overrightarrow {\rho}, u)$ is  $-1$-type. Indeed, The necessity  follows from
Proposition \ref {1.7}. the sufficiency follows from Proposition \ref {3.15}{\rm (i)}
 since $q_{C, C}^{(i)} =-1$ for any $i\in I_C(u, r)$.
\end {Remark}



\section { Program}\label {s4}
In this section the programs
to  compute the representatives of
conjugacy classes, centralizers of these representatives  and character tables
of these centralizers
in Weyl groups
of exceptional type are given.

 By using the
programs in GAP,  papers \cite {ZWCY08a, ZWCY08b}  obtained  the
representatives of conjugacy classes of Weyl groups of exceptional
type  and all character tables of centralizers of these
representatives. We use the results in \cite {ZWCY08a, ZWCY08b} and
the following program in GAP for Weyl group $W(E_6).$

gap$>$ L:=SimpleLieAlgebra("E",6,Rationals);;

gap$>$ R:=RootSystem(L);;

 gap$>$ W:=WeylGroup(R);Display(Order(W));

gap $>$ ccl:=ConjugacyClasses(W);;

gap$>$ q:=NrConjugacyClasses(W);; Display (q);

gap$>$ for i in [1..q] do

$>$ r:=Order(Representative(ccl[i]));Display(r);;

$>$ od; gap

$>$ s1:=Representative(ccl[1]);;cen1:=Centralizer(W,s1);;

gap$>$ cl1:=ConjugacyClasses(cen1);

gap$>$ s1:=Representative(ccl[2]);;cen1:=Centralizer(W,s1);;

 gap$>$
cl2:=ConjugacyClasses(cen1);

gap$>$ s1:=Representative(ccl[3]);;cen1:=Centralizer(W,s1);;

gap$>$ cl3:=ConjugacyClasses(cen1);

 gap$>$
s1:=Representative(ccl[4]);;cen1:=Centralizer(W,s1);;

 gap$>$
cl4:=ConjugacyClasses(cen1);

 gap$>$
s1:=Representative(ccl[5]);;cen1:=Centralizer(W,s1);;

gap$>$ cl5:=ConjugacyClasses(cen1);

 gap$>$
s1:=Representative(ccl[6]);;cen1:=Centralizer(W,s1);;

gap$>$ cl6:=ConjugacyClasses(cen1);

gap$>$ s1:=Representative(ccl[7]);;cen1:=Centralizer(W,s1);;

gap$>$ cl7:=ConjugacyClasses(cen1);

gap$>$ s1:=Representative(ccl[8]);;cen1:=Centralizer(W,s1);;

gap$>$ cl8:=ConjugacyClasses(cen1);

gap$>$ s1:=Representative(ccl[9]);;cen1:=Centralizer(W,s1);;

gap$>$ cl9:=ConjugacyClasses(cen1);

gap$>$ s1:=Representative(ccl[10]);;cen1:=Centralizer(W,s1);;

 gap$>$
cl10:=ConjugacyClasses(cen1);

gap$>$ s1:=Representative(ccl[11]);;cen1:=Centralizer(W,s1);;

gap$>$ cl11:=ConjugacyClasses(cen1);

gap$>$ s1:=Representative(ccl[12]);;cen1:=Centralizer(W,s1);;

gap$>$ cl2:=ConjugacyClasses(cen1);

gap$>$ s1:=Representative(ccl[13]);;cen1:=Centralizer(W,s1);;

 gap$>$
cl13:=ConjugacyClasses(cen1);

 gap$>$
s1:=Representative(ccl[14]);;cen1:=Centralizer(W,s1);;

gap$>$ cl14:=ConjugacyClasses(cen1);

gap$>$ s1:=Representative(ccl[15]);;cen1:=Centralizer(W,s1);;

gap$>$ cl15:=ConjugacyClasses(cen1);

gap$>$ s1:=Representative(ccl[16]);;cen1:=Centralizer(W,s1);;

gap$>$ cl16:=ConjugacyClasses(cen1);

gap$>$ s1:=Representative(ccl[17]);;cen1:=Centralizer(W,s1);;

gap$>$ cl17:=ConjugacyClasses(cen1);

gap$>$ s1:=Representative(ccl[18]);;cen1:=Centralizer(W,s1);;

gap$>$ cl18:=ConjugacyClasses(cen1);

gap$>$ s1:=Representative(ccl[19]);;cen1:=Centralizer(W,s1);;

gap$>$ cl19:=ConjugacyClasses(cen1);

gap$>$ s1:=Representative(ccl[20]);;cen1:=Centralizer(W,s1);;

gap$>$ cl20:=ConjugacyClasses(cen1);

gap$>$ s1:=Representative(ccl[21]);;cen1:=Centralizer(W,s1);;

 gap$>$
cl21:=ConjugacyClasses(cen1);

gap$>$ s1:=Representative(ccl[22]);;cen1:=Centralizer(W,s1);;

 gap$>$
cl22:=ConjugacyClasses(cen1);

gap$>$ s1:=Representative(ccl[23]);;cen1:=Centralizer(W,s1);;

gap$>$ cl23:=ConjugacyClasses(cen1);

gap$>$ s1:=Representative(ccl[24]);;cen1:=Centralizer(W,s1);;

$>$ cl24:=ConjugacyClasses(cen1);

gap$>$ s1:=Representative(ccl[25]);;cen1:=Centralizer(W,s1);

gap$>$ cl25:=ConjugacyClasses(cen1);

gap$>$ for i in [1..q] do

$>$ s:=Representative(ccl[i]);;cen:=Centralizer(W,s);;

$>$ char:=CharacterTable(cen);;Display (cen);Display(char);

 $>$ od;
gap$>$ for i in [1..q] do

$>$ s:=Representative(ccl[i]);;cen:=Centralizer(W,s);;

$>$ cl:=ConjugacyClasses(cen);;t:=NrConjugacyClasses(cen);;

$>$ for j in [1..t] do

 $>$ if s=Representative(cl[j]) then

 $>$
Display(j);break; $>$ fi;od;

$>$ od;

The programs for Weyl groups of $E_7$, $E_8$, $F_4$ and $ G_2$ are similar. It is possible that the order
of representatives of
conjugacy classes of $G$ changes when one uses the program.

\section {Tables about  $-1$- type} \label {s6}
In this section all $-1$- type bi-one Nichols algebras over Weyl groups of  exceptional type
 up to  graded pull-push
{\rm YD} Hopf algebra isomorphisms,  are listed in table 1--12.

Table 1 is about Weyl group $W(E_6)$; Tables 2--4 are about Weyl
group $W(E_7)$; Tables 5--10 are about Weyl group $W(E_8)$; Table 11
is about Weyl group $W(F_4)$; Table 12 is about Weyl group
$W(G_2)$.

\begin{tabular}{|l|l|l|l|l|l|}
  \hline
  $E_{6}$ &  &  & & &\\\hline
  $s_{i}$ & ${\rm cl}_{i}[p]$ & {\rm Order}$(s_{i})$ &
   the $j$ such that $\mathfrak{B}({\mathcal O}_{s_i},\chi _i^{(j)})$ is of
   $-1$-type
    & $\nu _i ^{(1)}$ & $\nu_i^{(2)}$\\\hline
  $s_{1}$ & ${\rm cl}_{1}[1]$ & 1 &  & 25 & 25 \\\hline
  $s_{2}$ & ${\rm cl}_{2}[3]$ & 4 &  4,5,6,7,17 & 20 & 15 \\\hline
  $s_{3}$ & ${\rm cl}_{3}[24]$ & 2 &  13,14,23,24,25 & 25 & 20 \\\hline
  $s_{4}$ & ${\rm cl}_{4}[17]$ & 4 &  2,4,10,15,16 & 20 & 15 \\\hline
  $s_{5}$ & ${\rm cl}_{5}[7]$ & 4 &  3,4,7,8 & 16 & 12 \\\hline
  $s_{6}$ & ${\rm cl}_{6}[2]$ & 2 &  9,10,11,12,16,17,18,19,25 & 25 & 16 \\\hline
  $s_{7}$ & ${\rm cl}_{7}[19]$ & 2 &  2,3,6,7,10,12,15,16,19,20 & 20 & 10 \\\hline
  $s_{8}$ & ${\rm cl}_{8}[26]$ & 3 &  & 27 & 27 \\\hline
  $s_{9}$ & ${\rm cl}_{9}[2]$ & 6 &  2,4,13 & 18 & 15 \\\hline
  $s_{10}$ & ${\rm cl}_{10}[2]$ & 2 &  2,3,7,8,11,12,15,16,19,20,22 & 22 & 11 \\\hline
  $s_{11}$ & ${\rm cl}_{11}[14]$ & 6 &  3,4,13 & 18 & 15 \\\hline
  $s_{12}$ & ${\rm cl}_{12}[27]$ & 3 &  & 27 & 27 \\\hline
  $s_{13}$ & ${\rm cl}_{13}[4]$ & 10 &  2 & 10 & 9 \\\hline
  $s_{14}$ & ${\rm cl}_{14}[9]$ & 5 &  & 10 & 10 \\\hline
  $s_{15}$ & ${\rm cl}_{15}[13]$ & 4 &  3,4,6,13,14 & 16 & 11 \\\hline
  $s_{16}$ & ${\rm cl}_{16}[3]$ & 8 &  2 & 8 & 7 \\\hline
  $s_{17}$ & ${\rm cl}_{17}[3]$ & 6 &  13 & 15 & 14 \\\hline
  $s_{18}$ & ${\rm cl}_{18}[9]$ & 12 &  2 & 12 & 11 \\\hline
  $s_{19}$ & ${\rm cl}_{19}[11]$ & 6 &  3,4 & 12 & 10 \\\hline
  $s_{20}$ & ${\rm cl}_{20}[2]$ & 9 &  & 9 & 9 \\\hline
  $s_{21}$ & ${\rm cl}_{21}[2]$ & 3 &  & 24 & 24 \\\hline
  $s_{22}$ & ${\rm cl}_{22}[13]$ & 6 &  2,4,13 & 18 & 15 \\\hline
  $s_{23}$ & ${\rm cl}_{23}[3]$ & 12 &  2 & 12 & 11 \\\hline
  $s_{24}$ & ${\rm cl}_{24}[19]$ & 6 &  10,11,12 & 21 & 18 \\\hline
  $s_{25}$ & ${\rm cl}_{25}[4]$ & 6 &  3,4,13 & 18 & 15 \\\hline
\end{tabular}
$$\hbox {Table } 1$$

\begin{tabular}{|l|l|l|l|l|l|}
  \hline
  $E_{7}$ &  &  &  &  &\\\hline
  $s_{i}$ & ${\rm cl}_{i}[p]$ & $Order(s_{i})$ & the $j$ such that $\mathfrak{B}({\mathcal O}_{s_i},\chi _i^{(j)})$ is of
   $-1$-type
    & $\nu _i^{(1)}$ & $\nu_i^{(2)}$ \\\hline
  $s_{1}$ & ${\rm cl}_{1}[1]$ & 1 &  & 60 & 60 \\\hline
  $s_{2}$ & ${\rm cl}_{2}[16]$ & 18 &  2 & 18 & 17 \\\hline
  $s_{3}$ & ${\rm cl}_{3}[15]$ & 9 &  & 18 & 18 \\\hline
  $s_{4}$ & ${\rm cl}_{4}[2]$ & 3 &  & 48 & 48 \\\hline
  $s_{5}$ & ${\rm cl}_{5}[4]$ & 6 &  2,3,4,8,11,12,14,40,43,44 & 48 & 38 \\\hline
  $s_{6}$ & ${\rm cl}_{6}[2]$ & 2 &  2,4,6,8,10,12,15,16,18,20,22,26,27,28,30,& 60  & 30 \\&&&32,35,36,38,41,42,44,46,48,50,52,54,56,58,60 &   &   \\\hline
  $s_{7}$ & ${\rm cl}_{7}[23]$ & 6 &  2,3,6,7,25,26 & 36 & 30 \\\hline
  $s_{8}$ & ${\rm cl}_{8}[23]$ & 3 &  & 54 & 54 \\\hline
  $s_{9}$ & ${\rm cl}_{9}[2]$ & 2 &  2,3,4,5,9,10,13,14,17,18,19,20,25,26,27, & 90 & 45 \\&&&28,33,35,36,39,40,43,44,45,46,51,52,53,54,55,& & \\&&&56,57,58,67,68,69,70,75,77,79,80,83,84,87,88 &  & \\\hline
  $s_{10}$ & ${\rm cl}_{10}[21]$ & 6 &  2,3,4,5,26,28 & 36 & 30 \\\hline
  $s_{11}$ & ${\rm cl}_{11}[2]$ & 2 &  2,3,7,8,11,12,15,16,19,20,27,28,29,30,31,& 74 & 37\\&&&32,37,38,39,40,42,44,47,48,51,52,55,56,58,61,& & \\&&&62,65,66,69,70,73,74 &  &  \\\hline
  $s_{12}$ & ${\rm cl}_{12}[24]$ & 6 &  3,4,7,8,26,28 & 36 & 30 \\\hline
  $s_{13}$ & ${\rm cl}_{13}[4]$ & 2 &  2,3,7,8,11,12,13,14,19,20,23,24,25,26, & 74 & 42
 \\&&& 27,28,37,38,39,40,42,43,45,46,49,50,
       & & \\&&& 55,56,57,59,60,63,64,69,70,73,74  &  &  \\\hline
  $s_{14}$ & ${\rm cl}_{14}[14]$ & 6 &  2,4,6,8,26,27,39,40,41,42 & 60 & 50 \\\hline
  $s_{15}$ & ${\rm cl}_{15}[3]$ & 3 &  & 66 & 66 \\\hline
  $s_{16}$ & ${\rm cl}_{16}[35]$ & 4 &  3,4,7,8,11,12,15,16,34,36,38,40,51,52,& 80 & 60 \\&&&55,56,59,60,63,64 &  &  \\\hline
  $s_{17}$ & ${\rm cl}_{17}[2]$ & 2 &  25,26,27,28,61,62,71,72,73,74,75,76,77,& 106 & 80 \\&&&78,79,80,81,82,99,100,101,102,103,104,105,106 &  &  \\\hline
  $s_{18}$ & ${\rm cl}_{18}[72]$ & 4 &  17,18,19,20,21,22,23,24,25,26,27,28,29,& 76 & 48\\&&&30,31,32,33,34,35,36,45,46,47,48,49,50,51,52 &  &  \\\hline
  $s_{19}$ & ${\rm cl}_{19}[2]$ & 2 &  25,26,27,28,29,30,31,32,43,44,45,46,47,48,& 90 & 60 \\&&&49,50,67,68,69,70,71,72,73,74,77,78,87,88,89,90 &  &  \\\hline
  $s_{20}$ & ${\rm cl}_{20}[4]$ & 8 &  2,4,6,8 & 32 & 28 \\\hline
  $s_{21}$ & ${\rm cl}_{21}[3]$ & 4 &  5,6,7,8,10,12,33,34,35,36,37,38,39,40,41,& 60 & 40 \\&&&42,43,44,50,52 & &  \\\hline
  $s_{22}$ & ${\rm cl}_{22}[8]$ & 12 &  2,5,7,8 & 48 & 44 \\\hline
  $s_{23}$ & ${\rm cl}_{23}[34]$ & 6 &  49,50,51,52 & 60 & 56 \\\hline
\end{tabular}
$$\hbox {Table } 2$$

\begin{tabular}{|l|l|l|l|l|l|}
  \hline
  $E_{7}$ &  &  &  &  &\\\hline
  $s_{i}$ & ${\rm cl}_{i}[p]$ & $Order(s_{i})$ & the $j$ such that $\mathfrak{B}({\mathcal O}_{s_i},\chi _i^{(j)})$ is of
   $-1$-type & $\nu _i^{(1)}$ & $\nu_i^{(2)}$ \\\hline
  $s_{24}$ & ${\rm cl}_{24}[47]$ & 4 &  2,5,7,8,10,13,15,16,34,35,38,39,51,52,53, & 80 & 60\\ &&&54,59,60,61,62 &  &  \\\hline
  $s_{25}$ & ${\rm cl}_{25}[15]$ & 8 &  3,4,7,8 & 32 & 28 \\\hline
  $s_{26}$ & ${\rm cl}_{26}[2]$ & 2 &  2,3,4,5,6,7,8,9,21,22,23,24,27,28,& 106 & 53 \\&&&37,38,39,40,41,42,43,44,53,54,55,56,57,58,& & \\&&&59,60,62,67,68,69,70,73,74,79,80,81,82,87,& & \\&&&88,89,90,95,96,97,98,101,102,104,106 &  &  \\\hline
  $s_{27}$ & ${\rm cl}_{27}[51]$ & 6 &  2,3,6,7,27,28,37,38,39,40 & 60 & 50 \\\hline
  $s_{28}$ & ${\rm cl}_{28}[36]$ & 4 &  2,3,6,7,10,11,14,15,34,36,37,39,49,& 80 & 60 \\&&&50,53,54,59,60,63,64 &  & \\\hline
  $s_{29}$ & ${\rm cl}_{29}[2]$ & 2 &  2,3,4,5,6,7,8,9,19,20,23,24,29,30,31,& 80 & 40 \\&&&32,37,38,39,40,45,46,47,48,53,54,55,56,& & \\&&&61,62,63,64,69,70,71,72,75,76,79,80 &  &  \\\hline
  $s_{30}$ & ${\rm cl}_{30}[2]$ & 2 &  2,3,4,5,6,7,8,9,21,22,23,24,33,34,35,& 80 & 40 \\&&&36,37,38,39,40,49,50,51,52,53,54,55,56,& & \\&&&65,66,67,68,69,70,71,72,77,78,79,80 &  &  \\\hline
  $s_{31}$ & ${\rm cl}_{31}[22]$ & 4 &  2,3,4,5,10,11,12,13,21,22,23,24,29,30,& 76 & 48 \\&&&31,32,34,36,38,40,41,42,45,46,49,50,53,54 &  &  \\\hline
  $s_{32}$ & ${\rm cl}_{32}[30]$ & 6 &  2,7,8,20,29,30,38 & 42 & 35 \\\hline
  $s_{33}$ & ${\rm cl}_{33}[30]$ & 6 &  3,4,7,8,26,28 & 36 & 30 \\\hline
  $s_{34}$ & ${\rm cl}_{34}[80]$ & 4 &  2,3,4,5,10,11,12,13,34,35,38,39,51,52, & 80 & 60 \\&&&53,54,59,60,61,62 &  &  \\\hline
  $s_{35}$ & ${\rm cl}_{35}[14]$ & 12 &  2,4,6,8 & 48 & 44 \\\hline
  $s_{36}$ & ${\rm cl}_{36}[40]$ & 6 &  3,5,6,8,11,13,14,16,49,52 & 60 & 50 \\\hline
  $s_{37}$ & ${\rm cl}_{37}[49]$ & 6 &  2,3,4,5,26,27,37,38,43,44 & 60 & 50 \\\hline
  $s_{38}$ & ${\rm cl}_{38}[50]$ & 6 &  2,5,7,8,27,28,31,32,50 & 54 & 45 \\\hline
  $s_{39}$ & ${\rm cl}_{39}[7]$ & 6 &  2,3,6,7 & 24 & 20 \\\hline
  $s_{40}$ & ${\rm cl}_{40}[2]$ & 10 &  2,3 & 20 & 18 \\\hline
  $s_{41}$ & ${\rm cl}_{41}[21]$ & 5 &  & 30 & 30 \\\hline
  $s_{42}$ & ${\rm cl}_{42}[33]$ & 12 &  3,4,7,8 & 48 & 44 \\\hline
  $s_{43}$ & ${\rm cl}_{43}[39]$ & 6 &  19,20,21,22,27,28 & 42 & 36 \\\hline
  $s_{44}$ & ${\rm cl}_{44}[5]$ & 4 &  3,5,7,9,10,12,14,16,19,21,23,25,26,28,30,32 & 64 & 48 \\\hline
  $s_{45}$ & ${\rm cl}_{45}[6]$ & 6 &  2,3,15,16,19,20,39,40,51,52,62 & 66 & 55 \\\hline
  $s_{46}$ & ${\rm cl}_{46}[6]$ & 6 &  2,3,6,7 & 24 & 20  \\\hline
  $s_{47}$ & ${\rm cl}_{47}[5]$ & 10 &  2,4 & 20 & 18 \\\hline
  $s_{48}$ & ${\rm cl}_{48}[12]$ & 10 &  2,4,21 & 30 & 27 \\\hline
\end{tabular}
$$\hbox {Table } 3$$

\begin{tabular}{|l|l|l|l|l|l|}
  \hline
  $E_{7}$ &  &  & & &\\\hline
  $s_{i}$ & ${\rm cl}_{i}[p]$ & $Order(s_{i})$ & the $j$ such that
   $\mathfrak{B}({\mathcal O}_{s_i},\chi _i^{(j)})$ is of
   $-1$-type& $\nu _i^{(1)}$ & $\nu_i^{(2)}$ \\\hline
  $s_{49}$ & ${\rm cl}_{49}[5]$ & 30 &  2 & 30 & 29 \\\hline
  $s_{50}$ & ${\rm cl}_{50}[15]$ & 15 &  & 30 & 30 \\\hline
  $s_{51}$ & ${\rm cl}_{51}[8]$ & 7 &  & 14 & 14 \\\hline
  $s_{52}$ & ${\rm cl}_{52}[2]$ & 14 &  2 & 14 & 13 \\\hline
  $s_{53}$ & ${\rm cl}_{53}[53]$ & 6 &  3,4,7,8,26,28,39,40,43,44 & 60 & 50 \\\hline
  $s_{54}$ & ${\rm cl}_{54}[5]$ & 8 &  3,5,6,8 & 32 & 28 \\\hline
  $s_{55}$ & ${\rm cl}_{55}[10]$ & 12 &  3,4 & 24 & 22 \\\hline
  $s_{56}$ & ${\rm cl}_{56}[14]$ & 8 &  3,5,6,8 & 32 & 28 \\\hline
  $s_{57}$ & ${\rm cl}_{57}[3]$ & 4 &  2,3,5,6,10,12,21,22,23,24,31,32,37,38,39,& 60 & 40 \\&&&40,41,42,50,52 &  &  \\\hline
  $s_{58}$ & ${\rm cl}_{58}[15]$ & 12 &  2,4 & 24 & 22 \\\hline
  $s_{59}$ & ${\rm cl}_{59}[38]$ & 12 &  3,5,6,8 & 48 & 44 \\\hline
  $s_{60}$ & ${\rm cl}_{60}[9]$ & 4 &  2,4,6,8,10,12,14,16,18,20,22,24,26,28,30,32 & 64 & 48 \\\hline
\end{tabular}
$$\hbox {Table } 4$$

\begin{tabular}{|l|l|l|l|l|l|}
  \hline
  $E_{8}$ &  &  & & &\\\hline
  $s_{i}$ & ${\rm cl}_{i}[p]$ & $Order(s_{i})$ & the $j$ such that $\mathfrak{B}
  ({\mathcal O}_{s_i},\chi _i^{(j)})$ is of
   $-1$-type& $\nu _i^{(1)}$ & $\nu_i^{(2)}$ \\\hline
  $s_{1}$ & ${\rm cl}_{1}[1]$ & 1 &  & 112 & 112 \\\hline
  $s_{2}$ & ${\rm cl}_{2}[29]$ & 30 &  2 & 30 & 29 \\\hline
  $s_{3}$ & ${\rm cl}_{3}[23]$ & 15 &  & 30 & 30 \\\hline
  $s_{4}$ & ${\rm cl}_{4}[2]$ & 5 &  & 45 & 45 \\\hline
  $s_{5}$ & ${\rm cl}_{5}[3]$ & 3 &  & 102 & 102 \\\hline
  $s_{6}$ & ${\rm cl}_{6}[6]$ & 10 &  6,7,27,41 & 45 & 41 \\\hline
  $s_{7}$ & ${\rm cl}_{7}[2]$ & 2 &  3,4,11,12,16,17,18,19,29,30,32,33,34,
& 112 & 67 \\&&&37,38,45,46,51,52,56,57,60,63,64,65,66,71,79, &  & \\&&&80,82,83,89,90,91,92,95,96,99,100,103,104, &  & \\&&&106,107,108,112 &  & \\\hline
  $s_{8}$ & ${\rm cl}_{8}[4]$ & 6 &  4,5,31,33,34,35,36,61,62,79,81,82,92,97 & 102 & 88 \\\hline
  $s_{9}$ & ${\rm cl}_{9}[4]$ & 30 &  2,4 & 60 & 58 \\\hline
  $s_{10}$ & ${\rm cl}_{10}[13]$ & 15 &  & 60 & 60 \\\hline
  $s_{11}$ & ${\rm cl}_{11}[5]$ & 5 &  & 70 & 70 \\\hline
  $s_{12}$ & ${\rm cl}_{12}[2]$ & 3 &  & 150 & 150 \\\hline
  $s_{13}$ & ${\rm cl}_{13}[46]$ & 10 &  2,4,6,8,41,44 & 60 & 54 \\\hline
  $s_{14}$ & ${\rm cl}_{14}[3]$ & 2 &  2,4,6,8,10,12,14,16,18,20,22,24,27,28,31,32, & 120 & 60 \\&&&34,36,38,40,42,44,48,49,50,54,55,56,58,60,62,64, &  & \\&&&67,68,71,72,74,76,79,80,83,84,86,88,90,92,94,96, &  & \\&&&98,100,102,104,106,108,110,112,114,116,118,120 &  &  \\\hline
  $s_{15}$ & ${\rm cl}_{15}[8]$ & 6 &  2,3,6,7,29,30,31,32,37,38,39,40,77,78, & 132 & 110 \\&&&79,80,101,102,103,104,123,124 &  &  \\\hline
  $s_{16}$ & ${\rm cl}_{16}[16]$ & 30 &  2,3 & 60 & 58 \\\hline
  $s_{17}$ & ${\rm cl}_{17}[9]$ & 10 &  2,3,23,24,43,44,62 & 70 & 63 \\\hline
  $s_{18}$ & ${\rm cl}_{18}[4]$ & 6 &  2,3,15,16,26,35,36,37,38,57,58,60,75,76, & 150 & 125 \\&&&87,88,99,100,102,117,118,128,135,136,146 &  &  \\\hline
  $s_{19}$ & ${\rm cl}_{19}[23]$ & 20 &  3,4 & 40 & 38 \\\hline
  $s_{20}$ & ${\rm cl}_{20}[40]$ & 10 &  41,42 & 50 & 48 \\\hline
  $s_{21}$ & ${\rm cl}_{21}[8]$ & 2 &  9,10,43,44,53,54,55,56,77,78,87,88,89,90, & 167 & 130 \\&&&91,92,109,110,111,112,121,134,151,152,153,154, &  & \\&&&155,156,157,158,159,162,163,164,165,166,167 &  &  \\\hline
  $s_{22}$ & ${\rm cl}_{22}[3]$ & 4 &  2,3,4,5,19,20,23,24,35,36,39,40,43,44,47, & 144 & 108 \\&&&48,66,68,77,78,79,80,85,86,87,88,90,92,115, &  & \\&&&116,119,120,131,132,135,136 &  &  \\\hline
  $s_{23}$ & ${\rm cl}_{23}[39]$ & 10 &  3,4,7,8,42,44 & 60 & 54 \\\hline
\end{tabular}
$$\hbox {Table } 5$$

\begin{tabular}{|l|l|l|l|l|l|}
  \hline
  $E_{8}$ &  &  & & &\\\hline
  $s_{i}$ & ${\rm cl}_{i}[p]$ & $Order(s_{i})$ & the $j$ such that $\mathfrak{B}
  ({\mathcal O}_{s_i},\chi _i^{(j)})$ is of
   $-1$-type& $\nu _i^{(1)}$ & $\nu_i^{(2)}$ \\\hline
  $s_{24}$ & ${\rm cl}_{24}[2]$ & 2 &  2,3,7,8,11,12,15,16,19,20,23,24,29,30,31,32, & 120 & 60 \\&&&35,36,39,40,43,44,51,52,53,54,55,56,59,60,63,64,69, &  & \\&&&70,71,72,75,76,81,82,83,84,87,88,91,92,95,96,99, &  & \\&&&100,103,104,107,108,111,112,115,116,119,120 &  &  \\\hline
  $s_{25}$ & ${\rm cl}_{25}[44]$ & 10 &  3,5,6,8,41 & 50 & 45 \\\hline
  $s_{26}$ & ${\rm cl}_{26}[2]$ & 2 &  27,28,29,30,31,32,33,34,55,56,79,80,81,82, & 167 & 105 \\&&&83,84,85,86,89,90,91,92,93,94,95,96,97,98,99,100, &  & \\&&&101,102,103,104,105,106,107,108,122,123,124,125, &  & \\&&&126,127,128,129,130,131,132,133,134,151,152,153, &  & \\&&&154,157,158,159,160,161,162,163 &  &  \\\hline
  $s_{27}$ & ${\rm cl}_{27}[28]$ & 20 &  2,3 & 40 & 38 \\\hline
  $s_{28}$ & ${\rm cl}_{28}[4]$ & 4 &  2,3,4,5,19,20,23,24,37,38,39,40,45,46,47, & 144 & 108 \\&&&48,66,68,78,79,80,81,82,88,89,90,91,92,115,116, &  & \\&&&119,120,131,132,135,136 &  &  \\\hline
  $s_{29}$ & ${\rm cl}_{29}[3]$ & 4 &  2,3,4,5,6,7,8,9,10,11,12,13,14,15,16,17,66, & 160 & 120 \\&&&68,70,72,73,75,77,79,97,98,101,102,105,106,109, &  & \\&&&110,115,116,119,120,123,124,127,128 &  &  \\\hline
  $s_{30}$ & ${\rm cl}_{30}[3]$ & 4 &  2,3,4,5,18,19,27,28,29,30,33,34,35,36,57, & 80 & 60 \\&&&58,59,60,73,76  &  & \\\hline
  $s_{31}$ & ${\rm cl}_{31}[72]$ & 6 &  25,26,31,32,39,40,45,46,68  & 72 & 63 \\\hline
  $s_{32}$ & ${\rm cl}_{32}[2]$ & 3 &  & 135 & 135 \\\hline
  $s_{33}$ & ${\rm cl}_{33}[58]$ & 8 &  3,4,7,8,34,35,49,50,55,56 & 80 & 70 \\\hline
  $s_{34}$ & ${\rm cl}_{34}[3]$ & 4 &  5,6,7,8,11,12,15,16,29,30,31,32,38,39,40,41, & 140 & 95 \\&&&42,55,56,57,58,63,64,65,66,85,86,87,88,89,90,91,92, &  & \\&&&98,99,100,101,102,108,119,120,131,132,135,136 &  &  \\\hline
  $s_{35}$ & ${\rm cl}_{35}[2]$ & 2 &  31,32,33,34,35,36,37,38,63,64,65,66,67,68, & 215 & 140 \\&&&69,70,71,72,73,74,104,105,106,107,108,109,110,111, &  & \\&&&112,113,114,115,118,119,130,131,132,133,146,147, &  & \\&&&148,149,150,151,152,153,154,155,156,157,159,168, &  & \\&&&169,170,171,172,173,174,175,178,179,182,183,198, &  & \\&&&199,200,201,202,203,208,209,210,211,213,215 &  &  \\\hline
  $s_{36}$ & ${\rm cl}_{36}[64]$ & 8 &  2,4,6,8,34,36,51,52,55,56  & 80 & 70 \\\hline
  $s_{37}$ & ${\rm cl}_{37}[44]$ & 4 &  9,10,11,12,13,14,15,16,17,18,21,22,43,44  & 44 & 30 \\\hline
\end{tabular}
$$\hbox {Table } 6$$

\begin{tabular}{|l|l|l|l|l|l|}
  \hline
  $E_{8}$ &  &  & & &\\\hline
  $s_{i}$ & ${\rm cl}_{i}[p]$ & $Order(s_{i})$ & the $j$ such that $\mathfrak{B}
  ({\mathcal O}_{s_i},\chi _i^{(j)})$ is of
   $-1$-type & $\nu _i^{(1)}$ & $\nu_i^{(2)}$ \\\hline
  $s_{38}$ & ${\rm cl}_{38}[2]$ & 2 &  21,22,23,24,25,26,27,28,33,34,35,36,41, & 105 & 65 \\&&&42,43,44,69,70,71,72,75,76,79,80,89,90,91,92, &  & \\&&&93,94,95,96,97,98,99,100,101,102,103,104 &  &  \\\hline
  $s_{39}$ & ${\rm cl}_{39}[73]$ & 4 &  4,5,8,9,10,11,14,15,20,21,24,25,26,27, & 112 & 84 \\&&&30,31,66,67,70,71,75,76,77,78,81,82,87,88 &  &  \\\hline
  $s_{40}$ & ${\rm cl}_{40}[13]$ & 24 &  2,3 & 48 & 46 \\\hline
  $s_{41}$ & ${\rm cl}_{41}[7]$ & 12 &  5,6,7,8,27,28,65,66,71,72 & 96 & 86 \\\hline
  $s_{42}$ & ${\rm cl}_{42}[4]$ & 6 &  49,50,51,52,53,54,55,56,91,92,93,94,95, & 150 & 132 \\&&&96,97,98,145,146 &  &  \\\hline
  $s_{43}$ & ${\rm cl}_{43}[98]$ & 12 &  4,5,6,7,12,13,14,15,98,100 & 120 & 110 \\\hline
  $s_{44}$ & ${\rm cl}_{44}[4]$ & 6 &  57,58,59,60,133,134,139,140,141,142 & 150 & 140 \\\hline
  $s_{45}$ & ${\rm cl}_{45}[3]$ & 4 &  5,6,7,8,11,12,15,16,29,30,31,32,33,34,35, & 140 & 95 \\&&&36,38,55,56,57,58,61,62,63,64,85,86,87,88,89,90, &  & \\&&&91,92,93,94,95,96,98,108,119,120,131,132,133,134 &  &  \\\hline
  $s_{46}$ & ${\rm cl}_{46}[19]$ & 8 &  2,4,6,8,9, 10,12,14,16 & 64 & 56 \\\hline
  $s_{47}$ & ${\rm cl}_{47}[3]$ & 4 &  33,34,35,36,37,38,39,40,41,42,43,44,45,46, & 178 & 114 \\&&&47,48,49,50,51,52,53,54,55,56,57,58,59,60,61,62, &  & \\&&&63,64,69,70,71,72,77,78,79,80,81,82,83,84,85,86, &  & \\&&&87,88,137,138,139,140,141,142,143,144,146,148, &  & \\&&&151,152,155,156,159,160 &  &  \\\hline
  $s_{48}$ & ${\rm cl}_{48}[4]$ & 6 &  2,3,4,5,6,7,8,9,51,52,55,56,59,60,89,90,91, & 150 & 125 \\&&&92,97,98,99,100,134,141,142 &  &  \\\hline
  $s_{49}$ & ${\rm cl}_{49}[4]$ & 8 &  4,5,6,7,33,34,35,36,66,67 & 80 & 70 \\\hline
  $s_{50}$ & ${\rm cl}_{50}[3]$ & 4 &  2,3,4,5,9,10,11,12,18,20,22,24,41,42,43,44, & 120 & 80 \\&&&49,50,57,58,59,60,63,64,69,70,71,72,75,76,77,78, &  & \\&&&79,80,85,86,98,100,102,104 &  &  \\\hline
  $s_{51}$ & ${\rm cl}_{51}[3]$ & 4 &  2,3,4,5,6,7,8,9,10,11,12,13,14,15,16,17,66, & 160 & 120 \\&&&67,70,71,73,76,77,80,97,98,103,104,105,106,111, &  & \\&&&112,115,116,117,118,123,124,125,126 &  &  \\\hline
  $s_{52}$ & ${\rm cl}_{52}[9]$ & 6 &  3,4,6,9,10,13,15,16,49,52,54,55 & 72 & 60 \\\hline
\end{tabular}
$$\hbox {Table } 7$$

\begin{tabular}{|l|l|l|l|l|l|}
  \hline
  $E_{8}$ &  &  & & &\\\hline
  $s_{i}$ & ${\rm cl}_{i}[p]$ & $Order(s_{i})$ & the $j$ such that $\mathfrak{B}
  ({\mathcal O}_{s_i},\chi _i^{(j)})$ is of
   $-1$-type& $\nu _i^{(1)}$ & $\nu_i^{(2)}$ \\\hline
  $s_{53}$ & ${\rm cl}_{53}[2]$ & 2 &  2,3,4,5,6,7,8,9,19,20,23,24,27,28,31,32, & 180 & 90 \\&&&37,38,39,40,45,46,47,48,53,54,55,56,61,62,63,64, &  & \\&&&66,68,71,72,75,76,79,80,83,84,89,90,91,92,97,98, &  & \\&&&99,100,105,106,107,108,113,114,115,116,121, &  & \\&&&122,123,124,129,130,131,132,137,138,139,140, &  & \\&&&145,146,147,148,150,152,154,156,159,160,163, &  & \\&&&164,167,168,171,172,175,176,179,180 &  &  \\\hline
  $s_{54}$ & ${\rm cl}_{54}[19]$ & 12 &  2,4,6,8 & 48 & 44 \\\hline
  $s_{55}$ & ${\rm cl}_{55}[4]$ & 6 &  37,38,39,40,43,44,45,46,59,60,61,62,103, & 126 & 108 \\&&&104,105,106,113,114 &  &  \\\hline
  $s_{56}$ & ${\rm cl}_{56}[2]$ & 3 &  & 144 & 144 \\\hline
  $s_{57}$ & ${\rm cl}_{57}[4]$ & 6 &  2,3,5,6,7,8,15,16,18,23,24,25,26,29,30,33,34, & 144 & 114 \\&&&60,63,64,66,105,106,111,112,113,114,136,139,140 &  &  \\\hline
  $s_{58}$ & ${\rm cl}_{58}[32]$ & 6 &  2,4,7,9,11,13,14,16,53,54,55,56,61,62, & 108 & 90 \\&&&63,64,99,100 &  &  \\\hline
  $s_{59}$ & ${\rm cl}_{59}[42]$ & 8 &  3,5,6,8,11,13,14,16 & 64 & 56 \\\hline
  $s_{60}$ & ${\rm cl}_{60}[10]$ & 6 &  6,7,8,9,10,11,12,13 & 48 & 40 \\\hline
  $s_{61}$ & ${\rm cl}_{61}[70]$ & 6 &  2,3,4,5,10,11,12,13,50,52,54,56 & 72 & 60 \\\hline
  $s_{62}$ & ${\rm cl}_{62}[2]$ & 2 &  2,3,4,5,6,7,8,9,21,22,23,24,29,30,31,32,41, & 180 & 90 \\&&&42,43,44,45,46,47,48,57,58,59,60,61,62,63,64,67, &  & \\&&&68,73,74,75,76,81,82,83,84,93,94,95,96,97,98,99, &  & \\&&&100,117,118,119,120,121,122,123,124,125,126, &  & \\&&&127,128,129,130,131,132,141,142,143,144,145, &  & \\&&&146,147,148,151,152,155,156,161,162,163, &  & \\&&&164,169,170,171,172,177,178,179,180 &  &  \\\hline
  $s_{63}$ & ${\rm cl}_{63}[4]$ & 6 &  25,26,33,34,51,52,53,54,68,77,78,79,80, & 135 & 117 \\&&&109,112,113,130,131 &  &  \\\hline
  $s_{64}$ & ${\rm cl}_{64}[36]$ & 6 &  2,3,15,16,21,22,38,39,47,48,53,54,74,75 & 84 & 70 \\\hline
  $s_{65}$ & ${\rm cl}_{65}[28]$ & 18 &  2,4,37 & 54 & 51 \\\hline
  $s_{66}$ & ${\rm cl}_{66}[26]$ & 9 &  & 54 & 54 \\\hline
  $s_{67}$ & ${\rm cl}_{67}[18]$ & 18 &  3,4 & 36 & 34 \\\hline
  $s_{68}$ & ${\rm cl}_{68}[4]$ & 6 &  2,3,5,6,7,8,15,16,21,22,23,24,27,28,79, & 96 & 76 \\&&&80,85,86,87,88 &  &  \\\hline
  $s_{69}$ & ${\rm cl}_{69}[12]$ & 18 &  2,4 & 36 & 34 \\\hline
  $s_{70}$ & ${\rm cl}_{70}[4]$ & 6 &  2,3,5,6,7,8,14,15,19,20,21,22,26,27,78, & 96 & 76 \\&&&79,83,84,85,86 &  &  \\\hline
\end{tabular}
$$\hbox {Table } 8$$

\begin{tabular}{|l|l|l|l|l|l|}
  \hline
  $E_{8}$ &  &  & & &\\\hline
  $s_{i}$ & ${\rm cl}_{i}[p]$ & $Order(s_{i})$ & the $j$ such that
  $\mathfrak{B}({\mathcal O}_{s_i},\chi _i^{(j)})$ is of
   $-1$-type& $\nu _i^{(1)}$ & $\nu_i^{(2)}$ \\\hline
  $s_{71}$ & ${\rm cl}_{71}[18]$ & 9 &  & 54 & 54 \\\hline
  $s_{72}$ & ${\rm cl}_{72}[24]$ & 18 &  2,5,6 & 54 & 51 \\\hline
  $s_{73}$ & ${\rm cl}_{73}[40]$ & 12 &   4,5,6,7 & 48 & 44 \\\hline
  $s_{74}$ & ${\rm cl}_{74}[3]$ & 4 &  2,3,4,5,9,10,11,12,17,20,22,23,41,42, & 120 & 80 \\&&&43,44,49,50,57,58,59,60,63,64,69,70,71,72, &  & \\&&&75,76,77,78,79,80,85,86,97,100,102,103 &  &  \\\hline
  $s_{75}$ & ${\rm cl}_{75}[7]$ & 12 &  2,3,5,6,27,28,61,62,65,66 & 96 & 86 \\\hline
  $s_{76}$ & ${\rm cl}_{76}[12]$ & 12 &  2,3,6,7,49,52 & 72 & 66 \\\hline
  $s_{77}$ & ${\rm cl}_{77}[36]$ & 12 &  2,9,10,38,47,48,74 & 84 & 77 \\\hline
  $s_{78}$ & ${\rm cl}_{78}[4]$ & 6 &  29,56,57,58,59,60,61,62,63,100,101, & 147 & 126 \\&&&102,103,104,105,106,107,108,127,129,130 &  &  \\\hline
  $s_{79}$ & ${\rm cl}_{79}[43]$ & 12 &  7,8,14,37,38 & 48 & 43 \\\hline
  $s_{80}$ & ${\rm cl}_{80}[4]$ & 4 &  11,12,17,18,19,20,26,35,36,39,40,41, & 59 & 43 \\&&&46,47,48,49 &  &  \\\hline
  $s_{81}$ & ${\rm cl}_{81}[15]$ & 20 &  2 & 20 & 19 \\\hline
  $s_{82}$ & ${\rm cl}_{82}[7]$ & 12 &  2,3,4,5,50,52,75,76,79,80 & 120 & 110 \\\hline
  $s_{83}$ & ${\rm cl}_{83}[3]$ & 4 &  2,3,4,5,6,7,8,9,35,36,39,40,43,44,47, & 200 & 150 \\&&&48,69,70,71,72,77,78,79,80,85,86,87,88,93, &  & \\&&&94,95,96,130,132,139,140,143,144,147,148, &  & \\&&&151,152,173,174,175,176,181,182,183,184 &  &  \\\hline
  $s_{84}$ & ${\rm cl}_{84}[107]$ & 12 &  4,5,6,7,12,13,14,15,98,100 & 120 & 110 \\\hline
  $s_{85}$ & ${\rm cl}_{85}[11]$ & 14 &  2,3 & 28 & 26 \\\hline
  $s_{86}$ & ${\rm cl}_{86}[26]$ & 7 &  & 28 & 28 \\\hline
  $s_{87}$ & ${\rm cl}_{87}[3]$ & 14 &  3,4 & 28 & 26 \\\hline
  $s_{88}$ & ${\rm cl}_{88}[18]$ & 14 &  2,4 & 28 & 26 \\\hline
  $s_{89}$ & ${\rm cl}_{89}[100]$ & 6 &  3,5,7,9,10,12,14,16,53,54,55,56,61,62, & 108 & 90 \\&&&63,64,99,100 &  &  \\\hline
  $s_{90}$ & ${\rm cl}_{90}[56]$ & 6 &  3,4,23,24,29,30,38,39,43,44,57,58,73,76 & 84 & 70 \\\hline
\end{tabular}
$$\hbox {Table } 9$$

\begin{tabular}{|l|l|l|l|l|l|}
  \hline
  $E_{8}$ &  &  & & &\\\hline
  $s_{i}$ & ${\rm cl}_{i}[p]$ & $Order(s_{i})$ & the $j$ such that
  $\mathfrak{B}({\mathcal O}_{s_i},\chi _i^{(j)})$ is of
   $-1$-type  & $\nu _i^{(1)}$ & $\nu_i^{(2)}$ \\\hline
  $s_{91}$ & ${\rm cl}_{91}[29]$ & 8 &  2,3,10,23,24 & 32 & 27 \\\hline
  $s_{92}$ & ${\rm cl}_{92}[8]$ & 24 &  2 & 24 & 23 \\\hline
  $s_{93}$ & ${\rm cl}_{93}[9]$ & 12 &  2,4,6,8,50,52 & 72 & 66 \\\hline
  $s_{94}$ & ${\rm cl}_{94}[4]$ & 6 &  2,3,5,6,7,8,38,41,42,51,52,75,76,77,78, & 126 & 105 \\&&&93,94,104,113,114,122 &  &  \\\hline
  $s_{95}$ & ${\rm cl}_{95}[35]$ & 12 &  2,4,17,18,29,30 & 72 & 66 \\\hline
  $s_{96}$ & ${\rm cl}_{96}[60]$ & 6 &  31,32,33,34,67,68 & 72 & 66 \\\hline
  $s_{97}$ & ${\rm cl}_{97}[7]$ & 12 &  2,3,4,5,50,52,75,76,79,80 & 120 & 110 \\\hline
  $s_{98}$ & ${\rm cl}_{98}[19]$ & 12 &  2,4 & 24 & 22 \\\hline
  $s_{99}$ & ${\rm cl}_{99}[59]$ & 12 &  6,7,8,9,10,11,12,13 & 96 & 88 \\\hline
  $s_{100}$ & ${\rm cl}_{100}[4]$ & 6 &  2,3,4,5,6,7,8,9,50,52,54,56,75,76,79,80, & 120 & 100 \\&&&83,84,87,88 &  &  \\\hline
  $s_{101}$ & ${\rm cl}_{101}[80]$ & 12 &  4,5,8,9,10,11,14,15 & 96 & 88 \\\hline
  $s_{102}$ & ${\rm cl}_{102}[84]$ & 6 &  2,3,6,7,10,11,14,15,51,52,53,54,59,60, & 108 & 90 \\&&&61,62,98,99 &  &   \\\hline
  $s_{103}$ & ${\rm cl}_{103}[82]$ & 12 &  2,9,10,37,39,40,74 & 84 & 77 \\\hline
  $s_{104}$ & ${\rm cl}_{104}[54]$ & 8 &  3,5,6,8,17,18,25,26,65,68  & 80 & 70 \\\hline
  $s_{105}$ & ${\rm cl}_{105}[35]$ & 30 &  2,4 & 60 & 58 \\\hline
  $s_{106}$ & ${\rm cl}_{106}[6]$ & 6 &  2,3,4,5,27,28,31,32,35,36,39,40,75,76, & 132 & 110 \\&&&79,80,99,100,103,104,122,124 &  &  \\\hline
  $s_{107}$ & ${\rm cl}_{107}[11]$ & 8 &   3,4 & 16 & 14 \\\hline
  $s_{108}$ & ${\rm cl}_{108}[4]$ & 6 &  2,3,4,5,27,28,31,32,53,54,55,56,58,81, & 150 & 125 \\&&&82,95,96,97,98,119,120,121,122,140,146 &  &  \\\hline
  $s_{109}$ & ${\rm cl}_{109}[4]$ & 6 &  2,3,4,5,6,7,8,9,49,50,51,52,73,74,75,76, & 120 & 100 \\&&&77,78,79,80 &  &  \\\hline
  $s_{110}$ & ${\rm cl}_{110}[2]$ & 24 &   3,4 & 48 & 46 \\\hline
  $s_{111}$ & ${\rm cl}_{111}[23]$ & 12 &  2,4,26 & 36 & 33 \\\hline
  $s_{112}$ & ${\rm cl}_{112}[92]$ & 6 &  3,5,6,8,23,24,35,36,49,50,61,62,73,76, & 108 & 90 \\&&&77,78,101,102 &  &  \\\hline
\end{tabular}
$$\hbox {Table } 10$$

\begin{tabular}{|l|l|l|l|l|l|}
  \hline
  $F_{4}$ &  &  & & &\\\hline
  $s_{i}$ & ${\rm cl}_{i}[p]$ & $Order(s_{i})$ & the $j$
  such that $\mathfrak{B}({\mathcal O}_{s_i},\chi _i^{(j)})$ is of
   $-1$-type & $\nu _i^{(1)}$ & $\nu_i^{(2)}$ \\\hline
  $s_{1}$ & ${\rm cl}_{1}[1]$ & 1 &  & 25 & 25\\\hline
  $s_{2}$ & ${\rm cl}_{2}[2]$ & 2 &  9,10,11,12,16,17,18,19,25 & 25 & 16 \\\hline
  $s_{3}$ & ${\rm cl}_{3}[25]$ & 2 &  17,18,19,20,25 & 25 & 20 \\\hline
  $s_{4}$ & ${\rm cl}_{4}[16]$ & 4 &  3,4,6,13,14 & 16 & 11 \\\hline
  $s_{5}$ & ${\rm cl}_{5}[5]$ & 3 &  & 18 & 18 \\\hline
  $s_{6}$ & ${\rm cl}_{6}[15]$ & 6 &  3,4,13 & 18 & 15 \\\hline
  $s_{7}$ & ${\rm cl}_{7}[10]$ & 2 &  2,3,4,5,10,12,15,16,19,20 & 20 & 10 \\\hline
  $s_{8}$ & ${\rm cl}_{8}[19]$ & 2 &  2,4,6,8,9,10,13,14,15,16 & 20 & 10 \\\hline
  $s_{9}$ & ${\rm cl}_{9}[16]$ & 4 &  3,5,6,8 & 16 & 12 \\\hline
  $s_{10}$ & ${\rm cl}_{10}[16]$ & 3 &  & 18 & 18 \\\hline
  $s_{11}$ & ${\rm cl}_{11}[12]$ & 6 &  3,4,13 & 18 & 15 \\\hline
  $s_{12}$ & ${\rm cl}_{12}[9]$ & 3 &  & 21 & 21 \\\hline
  $s_{13}$ & ${\rm cl}_{13}[21]$ & 6 &  10,11,12 & 21 & 18 \\\hline
  $s_{14}$ & ${\rm cl}_{14}[11]$ & 12 &  2 & 12 & 11 \\\hline
  $s_{15}$ & ${\rm cl}_{15}[10]$ & 6 &  2,3 & 12 & 10 \\\hline
  $s_{16}$ & ${\rm cl}_{16}[12]$ & 6 &  2,4 & 12 & 10 \\\hline
  $s_{17}$ & ${\rm cl}_{17}[19]$ & 2 &  2,3,4,6,7,10,12,15,16,19,20 & 20 & 10 \\\hline
  $s_{18}$ & ${\rm cl}_{18}[18]$ & 2 &  2,4,6,8,9,12,13,14,19,20 & 20 & 10 \\\hline
  $s_{19}$ & ${\rm cl}_{19}[6]$ & 4 &  2,4,6,8 & 16 & 12 \\\hline
  $s_{20}$ & ${\rm cl}_{20}[11]$ & 6 &  2,3 & 12 & 10 \\\hline
  $s_{21}$ & ${\rm cl}_{21}[12]$ & 6 &  2,4 & 12 & 10 \\\hline
  $s_{22}$ & ${\rm cl}_{22}[2]$ & 2 &  6,7,8,9,10,11,12,13 & 16 & 8 \\\hline
  $s_{23}$ & ${\rm cl}_{23}[4]$ & 8 &  2 & 8 & 7 \\\hline
  $s_{24}$ & ${\rm cl}_{24}[17]$ & 4 &  2,4,6,8,18 & 20 & 15 \\\hline
  $s_{25}$ & ${\rm cl}_{25}[9]$ & 4 &  2,4,6,8,17 & 20 & 15 \\\hline
   \end{tabular}
$$\hbox {Table } 11$$

\vskip 1cm

\begin{tabular}{|l|l|l|l|l|l|}
  \hline
  $G_{2}$ &  &  & & &\\\hline
  $s_{i}$ & ${\rm cl}_{i}[p]$ & $Order(s_{i})$ & the $j$ such that
  $\mathfrak{B}({\mathcal O}_{s_i},\chi _i^{(j)})$ is of
   $-1$-type   & $\nu _i^{(1)}$ & $\nu_i^{(2)}$ \\\hline
  $s_{1}$ & ${\rm cl}_{1}[1]$ & 1 &  & 6 & 6 \\\hline
  $s_{2}$ & ${\rm cl}_{2}[3]$ & 2 &  2,4 & 4 & 2 \\\hline
  $s_{3}$ & ${\rm cl}_{3}[3]$ & 2 &  2,4 & 4 & 2 \\\hline
  $s_{4}$ & ${\rm cl}_{4}[4]$ & 2 &  3,4,5 & 6 & 3 \\\hline
  $s_{5}$ & ${\rm cl}_{5}[3]$ & 6 &  2 & 6 & 5 \\\hline
  $s_{6}$ & ${\rm cl}_{6}[5]$ & 3 &  & 6 & 6 \\\hline
\end{tabular}
$$\hbox {Table } 12$$

\section{ Bi-one Nichols algebras over  Weyl groups of exceptional type}\label {s6}
In this section all  $-1$-type bi-one Nichols algebra over
 Weyl groups $G$ of exceptional type up to  graded pull-push
{\rm YD} Hopf algebra isomorphisms are  given.

In  Table 1--12, we use the following notations. $s_i$ denotes
the representative of $i$-th conjugacy class of $G$ ($G$ is the Weyl
group of exceptional type); $\chi _i ^{(j)}$  denotes the $j$-th character of
 $G^{s_i}$ for any $i$;
  $\nu_i ^{(1)}$ denotes    the number of
conjugacy classes of  the centralizer  $G^{s_i}$;
$\nu_i ^{(2)}$ denote the number of character $\chi _i^{(j)}$
 of
$G^{s_i}$ with non  $-1$-type
$\mathfrak{B}(\mathcal{O}_{s_i},\chi _i ^{(j)})$;  ${\rm cl}_i [j]$
denote that  $s_i$ is in  $j$-th conjugacy class of $G^{s_i}$.

We give one of the main results.

\begin {Theorem} \label {1} Let $G$ be a Weyl group of exceptional
type. Then

 {\rm (i)} For any bi-one Nichols algebra
 $\mathfrak{B}(\mathcal{O}_{s},\chi)$ over Weyl group $G$, there exist $s_{i}$
in the first column of the table
of $G$ and $j $ with $1 \le j \le \nu _{i}^{(1)}$ such that
 $(kG, \mathfrak{B}(\mathcal{O}_s,\chi)) \cong (kG, \mathfrak{B}(\mathcal{O}_{s_{i}},\chi _{i}
 ^{(j)}))$
 as graded pull-push {\rm YD} Hopf
algebras;

{\rm (ii)} $\mathfrak{B} ({\mathcal O}_{s_i}, \chi _i^{(j)})$  is of
$-1$-type if and only if $j$ appears in the  fourth column of the table of
$G$;

{\rm (iii)} ${\rm dim } (\mathfrak{B} ({\mathcal O}_{s_i}, \chi _i^{(j)}))
= \infty $ if $j$  does not appears in the fourth column of the table of
$G$.

\end {Theorem}

{\bf Proof.} {\rm (i)} We assume that $G$ is the Weyl group of $E_6$
without  loss of generality. There exists $s_i$ such that $s_i$ and $s$
are in  the same conjugacy class since $s_1, s_2, \cdots, s_{25}$ are
the representatives of all conjugacy classes of $G$. Lemma \ref
{1.1} and \cite [The remark of Pro. 1.5] {ZCZ08} or Proposition \ref {1.9}
yield that
 there
exists $j$ such that  $(kG, \mathfrak{B}(\mathcal{O}_s,\chi)) \cong
(kG, \mathfrak{B}(\mathcal{O}_{s_i},\chi _i ^{(j)}))$
 as  graded pull-push {\rm YD} Hopf
algebras,  since $\chi _i ^{(1)}, \chi _i^{(2)}, \cdots,
\chi _i ^{\nu _i^{(1)}}$ are all characters of all irreducible
representations of $G^{s_i}$.

{\rm (ii)} It follows from the  program.

{\rm (iii)} It follows from Lemma \ref {1.2}. $\Box$

By \cite {Ca72}, $W(G_2)$ is isomorphic to dihedral group $D_{6}$.
Set $y= s_5$ and $x= s_3$. It is clear that $xyx = y^{-1}$ with
${\rm ord } (y) = 6$ and  ${\rm ord } (x) = 2$. Thus  it follows
from \cite [Table 2] {AF07} that  ${\rm dim } (\mathfrak{B}
({\mathcal O}_{s_5}, \chi _5^{(2)})) =4 < \infty $.

It is clear that if there exists $\phi\in {\rm Aut} (G)$ such that
$\phi (s_i)= s_j$ then ${\rm ord } (s_i) = {\rm ord } (s_j)$, $\nu
_i^{(1)}= \nu_j ^{(1)}$, $\nu _i^{(2)}= \nu_j ^{(2)}$ for Weyl group
$G$ of exceptional type. Consequently, the representative system of
iso-conjugacy classes of $W(E_6)$ is $\{ s_i \mid 1\le i \le 25 \}$.
The representative system of iso-conjugacy classes of $W(F_4)$ is
$\{ s_i \mid 1\le i \le 25, i\not=  8, 10,  11, 16, 17, 18, 19, 20,
21, 25 \}$. The representative system of iso-conjugacy classes of
$W(G_2)$ is $\{ s_1, s_2, s_4, s_5, s_6 \}$.

\section { Pointed Hopf algebras over  Weyl groups of exceptional type}\label {s7}

In this section all central quantum linear spaces over Weyl groups of exceptional type are found.


\begin {Lemma}\label {7.1} $Z(W(E_6))= \{1\}$;
$Z(W(E_7))= \{1, s_6\}$;  $Z(W(E_8))= \{1, s_7\}$;
$Z(W(F_4))= \{1,s_2 \}$;  $Z(W(G_2))= \{1, s_4\}$.

\end {Lemma}
{\bf Proof.} If $s_i \in Z(G)$, then $G^{s_i} = G$.

{\rm (i)} Let  $G= W(W_6)$. The number of conjugacy classes of  $G$ is 25 by table 1. The numbers
of conjugacy classes of both  $G^{s_3}$ and  $G^{s_6}$ also
are  25.  $G$, $G^{s_3}$ and $G^{s_6}$ have 16, 8 and 4 one dimensional
 representations,   respectively, according to the character tables in \cite {ZWCY08a}.
Thus $s_3$ and $s_6$ do not belong to the center of $G$.

{\rm (ii)} Let  $G= W(W_7)$. The number of conjugacy classes of  $G$, $G^{s_6}$,
 $G^{s_{14}}$, $G^{s_{21}}$, $G^{s_{23}}$, $G^{s_{27}}$, $G^{s_{36}}$, $G^{s_{37}}$,
$G^{s_{53}}$, $G^{s_{57}}$ is  60  by table 1 --4. They have 2, 2, 24, 8, 3, 24,
48, 24, 24 and 8
one dimensional
 representations,   respectively, according to the character tables in \cite {ZWCY08a}.
Thus they do not belong to the center of $G$ but $s_6$. Obviously
$s_6 $
 $\in Z(G)$.

{\rm (iii)} Let  $G= W(W_8)$. The number of conjugacy classes of  $G$, $G^{s_7}$,
 and $G^{s_{39}}$
 is  112  by table 5--10. They have 2, 2 and 64
one dimensional
 representations,   respectively, according to the character tables in \cite {ZWCY08b}.
Thus $s_{39}$ does not belong to the center of $G$. Obviously $s_7 $
 $\in Z(G)$.

{\rm (iv)} Let  $G= W(F_4)$. The number of conjugacy classes of  $G$, $G^{s_2}$,
 and $G^{s_{3}}$
 is 25  by table 11. They have 4, 4 and 16
one dimensional
 representations,   respectively, according to the character tables in \cite {ZWCY08a}.
Thus $s_3$ does not belong to the center of $G$. Obviously $s_2 $
$\in Z(G)$.

{\rm (v)} Let  $G= W(G_2)$. The number of conjugacy classes of  $G$, $G^{s_4}$,  $G^{s_5}$,
 and $G^{s_{6}}$
 is 6  by table 12. They have 4, 4, 6 and 6
one dimensional
 representations,   respectively, according to the character tables in \cite {ZWCY08a}.
Thus $s_5$ and $s_6$ do not belong to the center of $G$. Obviously
$s_4 \in Z(G)$. $\Box$

We give the other main result.
\begin {Theorem}\label {2} Every central quantum linear space
$\mathfrak {B}  (G, r, \overrightarrow {\rho}, u)$  over Weyl Groups of exceptional type
is one case in the following:

{\rm (i)}  $G= W(E_7)$, $C= \{  s_6 \}$, $r = r_CC$ and
$\chi  _C^{(i)} $ $\in  \{  \chi _6^{(j)} \mid j $ =  2, 4, 6, 8, 10, 12, 15, 16, 18, 20, 22,
26, 27, 28, 30,
 32, 35, 36, 38, 41, 42, 44, 46, 48, 50, 52, 54, 56, 58, 60  $\}$ for any $i \in I_C(r, u)$.

{\rm (ii)} Let $G= W(E_8)$, $C= \{  s_7 \}$, $r = r_CC$ and $\chi _C^{(i)} $
$\in  \{  \chi _7^{(j)}
 \mid j$ = 3, 4, 11, 12, 16, 17, 18, 19,29, 30, 32, 33, 34, 37, 38, 45, 46,
51, 52, 56, 57, 60, 63, 64, 65, 66, 71, 79,
80, 82, 83, 89, 90, 91, 92, 95, 96, 99, 100, 103, 104, 106, 107, 108, 112
 $\}$ for any $i \in I_C(r, u)$.

{\rm (iii)} Let $G= W(F_4)$, $C= \{  s_2 \}$, $r = r_CC$ and $\chi _C^{(i)} $
$\in  \{  \chi _2^{(j)}
 \mid j$ = 9, 10, 11, 12, 16, 17, 18, 19, 25
 $\}$ for any $i \in I_C(r, u)$.

{\rm (iv)} Let $G= W(G_2)$, $C= \{  s_4 \}$, $r = r_CC$ and $\chi _C^{(i)} $
 $\in  \{  \chi _4^{(3)},
\chi _4^{(4)},  \chi _4^{(5)} \}$ for any $i \in I_C(r, u)$.

\end {Theorem}
{\bf Proof.} Let us  first consider the case of {\rm (i)}. By Theorem  \ref {1} and
Table 2, ${\rm RSR}(G, r, \overrightarrow \rho. u)$ is of
 $-1$-type. Applying Lemma \ref {7.1} we have that $\mathfrak {B}
(G, r, \overrightarrow \rho. u)$
 is  a central quantum linear space. Similarly, $ \mathfrak {B}
(G, r, \overrightarrow \rho. u)$
 is  a central quantum linear space under the other case.

Conversely, if $\mathfrak {B}
(G, r, \overrightarrow \rho. u)$
 is  a central quantum linear space over Weyl Group $G$ of exceptional type, then
for any $C\in {\mathcal K}_r(G)$, $C$ has to be $\{s_6\}$ with $G = W(E_7)$
or $\{s_6\}$ with $G= W(E_8)$ or $\{s_2\}$ with $G=W(F_4)$ or
$\{s_4\}$ with $G= W(G_2)$ by Lemma \ref {7.1}. This implies $r= r_CC$ and $C$ is one  case in this theorem.
 Furthermore, every  bi-one type
 ${\rm RSR} (G, {\mathcal O}_{u(C)},
\rho _C^{(i)})$ for any $i\in I_C(r,u)$ is of $-1$-type by Proposition \ref {1.7}.
Applying Theorem \ref {1} and Table 2, Table 5, Table 11 and Table 12, we have that
$\chi _C^{(i)}$ has to be one   case in this theorem for any $i\in I_C(r,u)$. $\Box$

In other words we have

\begin {Remark}\label {7.2} Let $G$ be a Weyl Group of exceptional type
and  $M = M({\mathcal O}_a, \rho ^{(1)})\oplus M({\mathcal O}_a,
\rho ^{(2)}) \oplus \cdots \oplus M({\mathcal O}_a, \rho ^{(m)})$ is
a  {\rm YD} module over $kG$. Then $\mathfrak B (M)$ is finite
dimensional in the following cases:

{\rm (i)}  $G= W(E_7)$, $a= s_6$ and the characters of $\rho
^{(1)}$, $\rho ^{(2)}$, $ \cdots,$ $\rho ^{(m)}$ are in   $ \{ \chi
_6^{(j)} \mid j $ =  2, 4, 6, 8, 10, 12, 15, 16, 18, 20, 22, 26, 27,
28, 30,
 32, 35, 36, 38, 41, 42, 44, 46, 48, 50, 52, 54, 56, 58, 60  $\}$.

{\rm (ii)}  $G= W(E_8)$, $a=  s_7 $ and  the characters of $\rho
^{(1)}$, $\rho ^{(2)}$, $ \cdots,$ $\rho ^{(m)}$ are in
 $ \{  \chi _7^{(j)}
 \mid j$ = 3, 4, 11, 12, 16, 17, 18, 19,29, 30, 32, 33, 34, 37, 38, 45, 46,
51, 52, 56, 57, 60, 63, 64, 65, 66, 71, 79, 80, 82, 83, 89, 90, 91,
92, 95, 96, 99, 100, 103, 104, 106, 107, 108, 112
 $\}$.

{\rm (iii)}  $G= W(F_4)$, $a=  s_2 $ and the characters of $\rho
^{(1)}$, $\rho ^{(2)}$, $ \cdots,$ $\rho ^{(m)}$ are in $ \{ \chi
_2^{(j)}
 \mid j$ = 9, 10, 11, 12, 16, 17, 18, 19, 25
 $\}$.

{\rm (iv)} $G= W(G_2)$, $a  s_4 $ and the characters of $\rho
^{(1)}$, $\rho ^{(2)}$, $ \cdots,$ $\rho ^{(m)}$ are in $\{ \chi
_4^{(3)}, \chi _4^{(4)},  \chi _4^{(5)} \}$.

  \end {Remark}

\section { Nichols algebras of reducible {\rm YD} modules} \label {s8}

In this section it is proved that except a few cases Nichols algebras of reducible
Yetter-Drinfeld modules over  Weyl groups of exceptional type are infinite dimensional.

 ${\mathcal O}_{s_i} $ and
${\mathcal O}_{s_j}$ are  said to be square-commutative if  $stst =
tsts$ for any $s\in {\mathcal O}_{s_i} $, $t\in {\mathcal O}_{s_j}$.
$a$ and $b$ are  said to be square-commutative if  $abab = baba$.

\begin {Lemma} \label {8.1} Let $G$ be a Weyl group of Exceptional
Type.

{\rm (i)} ${\mathcal O}_{s_i}$ and  ${\mathcal O}_{s_j}$ are not
commutative for any $i$ and $j$
 with $i, j\not = 1$ when
   $G= W(E_6)$.

{\rm (ii)} ${\mathcal O}_{s_i}$ and  ${\mathcal O}_{s_j}$ are not
square-commutative when $G= W(E_7)$ and  $ (i, j) \not= (9, 11),  $
$(9, 13)$, $(11, 19)$, $(13, 19)$ with $i, j \not= 1, 6$.

{\rm (iii)} ${\mathcal O}_{s_i}$ and  ${\mathcal O}_{s_j}$ are not
square-commutative when $G= W(E_8)$ and  $ (i, j) \not= (5, 14),  $
$(5, 24)$, $(8, 14)$, $(8, 24)$, $(14, 35)$, $(14, 80)$, $(24, 35)$,
$(24, 80)$ with $i, j \not= 1, 7$.

{\rm (iv)} ${\mathcal O}_{s_i}$ and  ${\mathcal O}_{s_j}$ are not
square-commutative when $G= W(F_4)$ and  $ (i, j) \not= (3, 3),  $
 $(3,4 )$, $(3, 7)$, $(3,8 )$, $(3, 17)$, $(3, 18 )$, $(3, 24)$,
$(3, 25 )$, $(4, 4 )$, $(4, 7 )$, $(4, 8 )$, $(4, 17 )$, $(4, 18 )$,
$(4, 24 )$, $(4, 25 )$, $(7, 12 )$, $(7, 13 )$, $(7, 17 )$, $(7, 18
)$, $(8, 12)$, $(8, 13)$,  $(8, 17 )$, $(8, 18 )$, $(12, 17)$, $(12,
18 )$, $(13, 17 )$, $(13, 18 )$ with $i, j \not = 1, 2$.

{\rm (v)}  ${\mathcal O}_{s_i}$ and  ${\mathcal O}_{s_j}$ are not
square-commutative when $G= W(G_2)$ and  $ (i, j) \not= (2, 5),  $
$(2, 6)$, $(3, 5)$,  $(3, 6 )$, $(5,5 )$, $(5,6 )$, $(6,6 )$ with
$i, j \not= 1, 4$.

\end {Lemma}
{\bf Proof.} Let $A := \{ (i, j ) \mid $ $ (i, j) = (9, 11),$ $(9,
13)$, $(11, 19)$, $(13, 19)$,  or  $i, j = 1, 6$ $\}$, $B := \{ (i,
j ) \mid $ $ (i, j) =$ $(5, 14),  $ $(5, 24)$, $(8, 14)$, $(8, 24)$,
$(14, 35)$, $(14, 80)$, $(24, 35)$, $(24, 80)$, or $i, j = 1, 7$
$\}$, $C := \{ (i, j ) \mid $ $ (i, j)= (3, 3),  $
 $(3,4 )$, $(3, 7)$, $(3,8 )$, $(3, 17)$, $(3, 18 )$, $(3, 24)$,
$(3, 25 )$, $(4, 4 )$, $(4, 7 )$, $(4, 8 )$, $(4, 17 )$, $(4, 18 )$,
$(4, 24 )$, $(4, 25 )$, $(7, 12 )$, $(7, 13 )$, $(7, 17 )$, $(7, 18
)$, $(8, 12)$, $(8, 13)$,  $(8, 17 )$, $(8, 18 )$, $(12, 17)$, $(12,
18 )$, $(13, 17 )$, $(13, 18 )$,  or  $i, j  = 1, 2$ $\}$.

{\rm (i)} It follows from Table 13.

{\rm (ii)} ${\mathcal O}_{s_{i}}$ and ${\mathcal O}_{s_{j}}$ are
square-commutative in $W(E_7)$ for $(i, j) \in A$. $s_i$ and $s_j $
are not square-commutative if  $(i, j) \not\in A $ and there does
not exist  $t$ such that $s_i$ and $s_ts_js_t^{-1}$  are in table
14--16.

{\rm (iii)} ${\mathcal O}_{s_{i}}$ and ${\mathcal O}_{s_{j}}$ are
square-commutative in $W(E_8)$ for $(i, j) \in B$.
 $s_i$ and $s_{110}s_j s_{110}^{-1}$ in $W(E_8)$ are not
square-commutative if  $(i, j) \not\in B $ and there does not exist
$t$ such that $s_i$ and $s_ts_js_t^{-1}$ are  in table 17.

{\rm (iv)} ${\mathcal O}_{s_{i}}$ and ${\mathcal O}_{s_{j}}$ are
square-commutative in $W(F_4)$ for $(i, j) \in C$.  $s_i$ and
$s_{3}s_j s_{3}^{-1}$ are not square-commutative in $W(F_4)$ if $(i,
j) \not\in C $ and there does not exist  $t$ such that $s_i$ and
$s_ts_js_t^{-1}$ are in table 18.

{\rm (v)} ${\mathcal O}_{s_{i}}$ and ${\mathcal O}_{s_{j}}$ are
square-commutative in $W(G_2)$ for any $(i, j)$ but $(i, j)= (2,3),
$ $(2, 2),$ $(3,3 )$. $s_2$ and $s_{5}s_3 s_{5}^{-1}$, $s_2$ and
$s_{6}s_2 s_{6}^{-1}$, $s_3$ and $s_{5}s_3 s_{5}^{-1}$ are not
square-commutative, respectively. $\Box$

Note that we have proved that ${\mathcal O}_{s_{i}}$ and ${\mathcal
O}_{s_{j}}$ are square-commutative in $G=(W(E_7))$, $G=(W(E_8))$ and
$G=(W(F_2))$ if and only if $(i, j)\in A$, $B$, $C$, respectively.
The programs to prove that ${\mathcal O}_{s_{i}}$ and ${\mathcal
O}_{s_{j}}$ in $W(E_7)$ are square-commutative are the following:

gap$>$ L:=SimpleLieAlgebra("E",7,Rationals);;

 gap$>$ R:=RootSystem(L);;

gap$>$ W:=WeylGroup(R);;

gap$>$ ccl:=ConjugacyClasses(W);

 gap$>$
q:=NrConjugacyClasses(W);;Display (q);

gap$>$ con1:=Elements(ccl[11]);;m:=Size(con1);

gap$>$ for k in [1..m] do

$>$ s:=con1[k];

$>$ con2:=Elements(ccl[19]);n:=Size(con2);

$>$ for l in [1..n] do

$>$ t:=con2[l];

$>$ if $(s*t)\hat {\ }2=(t*s)\hat {\ }2$ then

$>$ Print( " k=",k," AND l=",l, " ${\setminus n}$");

$>$ fi;

$>$ od;

$>$ od;

For any reducible {\rm YD} module $M$ over $kG$, there are at least
two irreducible {\rm YD} sub-modules of $M$. Therefore  we only
consider the direct sum of two irreducible {\rm YD} modules.

We give the final main result.
\begin {Theorem} \label {3} Let $G$ be a Weyl group of  Exceptional
Type. Then ${\rm dim } ({\mathfrak B } ( M({\mathcal O}_{s _i}, \rho
^{(1)})\oplus M({\mathcal O}_{s_j}, \rho^{(2)})) = \infty$ in the
following cases:

{\rm (i)}   $G= W(E_6)$ .

{\rm (ii)}  $G= W(E_7)$ and  $ (i, j) \not= (9, 11),  $ $(9, 13)$,
$(11, 19)$, $(13, 19)$ and $i, j \not=  6$.

{\rm (iii)} $G= W(E_8)$ and  $ (i, j) \not= $ $(8, 14)$, $(8, 24)$,
$(14, 35)$, $(14, 80)$, $(24, 35)$, $(24, 80)$ and $i, j \not= 7$.

{\rm (iv)} $G= W(F_4)$ and  $ (i, j) \not= (3, 3),  $
 $(3,4 )$, $(3, 7)$, $(3,8 )$, $(3, 17)$, $(3, 18 )$, $(3, 24)$,
$(3, 25 )$, $(4, 4 )$, $(4, 7 )$, $(4, 8 )$, $(4, 17 )$, $(4, 18 )$,
$(4, 24 )$, $(4, 25 )$, $(7, 13 )$, $(7, 17 )$, $(7, 18 )$, $(8,
13)$, $(8, 17 )$, $(8, 18 )$, $(13, 17 )$, $(13, 18 )$ and $i, j
\not = 2$.

{\rm (v)}  $G= W(G_2)$ and  $ (i, j) \not= (2, 5)$,   $(3, 5)$,
$(5,5 )$   and $i, j \not=  4$.

\end {Theorem}
{\bf Proof.}  It follows from \cite [Theorem 8.2, Theorem 8.6] {HS}
and Lemma \ref {8.1}. Note that the orders of $s_{12}$ in $W(F_4)$ ,
 $s_5$ in $W(E_8)$  and   $s_6$ in $W(G_2)$  are odd. $\Box$

\begin{tabular}{|l|l|}
  \hline
  $E_{6}$ &   \\\hline

$s_i$ & $s_i$ and $s_t s_js_t^{-1}$ are not  commutative \\\hline

  $s_{2}$ & $s_{7}$$s_{2}$$s_{7}^{-1},$  $s_{7}$$s_{3}$$s_{7}^{-1},$
  $s_{7}$$s_{4}$$s_{7}^{-1},$ $s_{5}$, $s_{6}$,  $s_{7}$, $s_{8}$, $s_{9}$, $s_{10}$, $s_{11}$,
  $s_{12}$, $s_{13}$, $s_{14}$, $s_{5}s_{15}s_{5}^{-1},$ \\
  &$s_{16}$,  $s_{17}$,  $s_{18}$,  $s_{19}$,  $s_{20}$,  $s_{21}$,  $s_{22}$,  $s_{23}$,
    $s_{24}$,  $s_{25}$ \\\hline

  $s_{3}$ & $s_{7}$$s_{3}$$s_{7}^{-1},$  $s_{7}$$s_{4}$$s_{7}^{-1},$  $s_{7}$$s_{5}$$s_{7}^{-1},$   $s_{8}$$s_{6}$$s_{8}^{-1},$   $s_{7}$, $s_{8}$, $s_{9}$, $s_{10}$, $s_{11}$, $s_{12}$, $s_{13}$  \\
  &$s_{14}$,   $s_{7}s_{15}s_{7}^{-1},$  $s_{16}$,  $s_{17}$,  $s_{18}$,  $s_{19}$,
    $s_{20}$,  $s_{21}$,  $s_{22}$,  $s_{23}$,  $s_{24}$,  $s_{25}$ \\\hline

   $s_{4}$ & $s_{7}$$s_{4}$$s_{7}^{-1},$  $s_{7}$$s_{5}$$s_{7}^{-1},$ $s_{8}$$s_{6}$$s_{8}^{-1},$   $s_{7}$, $s_{8}$, $s_{9}$, $s_{10}$, $s_{11}$, $s_{12}$, $s_{13}$  \\
  &$s_{14}$,   $s_{7}$$s_{15}s_{7}^{-1}$,  $s_{16}$,  $s_{17}$,  $s_{18},$  $s_{19}$,
      $s_{20}$,  $s_{21}$, $s_{22}$,  $s_{23}$,  $s_{24}$,  $s_{25}$ \\\hline

   $s_{5}$ & $s_{2}$$s_{5}$$s_{2}^{-1},$  $s_2s_{6}s_2^{-1}$,   $s_2s_{7}s_2^{-1}$,
    $s_2s_{8}s_2^{-1}$,   $s_2s_{9}s_2^{-1}$, $s_2s_{10}s_2^{-1}$,    $s_{11}$, $s_{12}$, $s_{13}$  \\
  &$s_{14}$,   $s_{15}$,   $s_{16}$,  $s_{17}$,  $s_{18}$,  $s_{19}$,  $s_2s_{20}s_2^{-1}$,
   $s_{21},$  $s_{22}$,  $s_{23}$,  $s_{24}$,  $s_{25}$ \\\hline

   $s_{6}$ &    $s_{2}$$s_{6}$$s_{2}^{-1}$,  $s_{2}$$s_{7}$$s_{2}^{-1},$
    $s_2s_{8}s_2^{-1}$,   $s_2s_{9}s_2^{-1}$,  $s_{2}$$s_{10}$$s_{2}^{-1},$ $s_{11}$, $s_{12}$, $s_{12}$  \\
  &$s_{14}$, $s_{15}$, $s_{16}$, $s_{5}$$s_{17}$$s_{5}^{-1},$  $s_{18}$,  $s_{19}$,
    $s_2s_{20}s_2^{-1}$,  $s_{21}$, $s_{5}$$s_{22}$$s_{5}^{-1},$  $s_{23}$, $s_{24}$, $s_{5}$$s_{25}$$s_{5}^{-1},$ \\\hline

    $s_{7}$ &   $s_{2}$$s_{7}$$s_{2}^{-1},$    $s_{2}$$s_{8}$$s_{2}^{-1},$
    $s_{9}$, $s_{10}$, $s_{11}$, $s_{12}$, $s_{13}$  \\
  &$s_{14}$, $s_2s_{15}s_2^{-1}$, $s_{16}$, $s_{17}$, $s_{18}$, $s_{19}$, $s_{20}$, $s_{21}$, $s_{22}$, $s_{23}$, $s_{24}$, $s_{25}$ \\\hline

 $s_{8}$ & $s_{2}$$s_{8}$$s_{2}^{-1},$  $s_{2}$$s_{9}$$s_{2}^{-1},$ $s_{3}$$s_{10}$$s_{3}^{-1},$   $s_{2}$$s_{11}$$s_{2}^{-1},$  $s_{2}$$s_{12}$$s_{2}^{-1},$  $s_{13}$  \\
  &$s_{14}$, $s_{15}$, $s_{16}$, $s_{17}$, $s_{18}$, $s_{19}$,  $s_2s_{20}s_2^{-1}$, $s_{21}$, $s_{22}$, $s_{23}$, $s_{24}$, $s_{25}$ \\\hline

$s_{9}$ &  $s_{2}$$s_{9}$$s_{2}^{-1},$ $s_{2}$$s_{10}$$s_{2}^{-1},$   $s_{2}$$s_{11}$$s_{2}^{-1},$  $s_{2}$$s_{12}$$s_{2}^{-1},$  $s_{13}$  \\
  &$s_{14}$, $s_{15}$, $s_{16}$, $s_{17}$, $s_{18}$, $s_{19}$, $s_2s_{20}s_2^{-1}$, $s_{21}$, $s_{22}$, $s_{23}$, $s_{24}$, $s_{25}$ \\\hline

$s_{10}$ &   $s_{2}$$s_{10}$$s_{2}^{-1},$   $s_{2}$$s_{11}$$s_{2}^{-1},$  $s_{2}$$s_{12}$$s_{2}^{-1},$  $s_{2}$$s_{13}$$s_{2}^{-1},$  $s_{2}$$s_{14}$$s_{2}^{-1},$  \\
  &  $s_{15}$, $s_{16}$, $s_{17}$, $s_{18}$, $s_{19}$, $s_2s_{20}s_2^{-1}$,
    $s_{21}$, $s_{22}$, $s_{23}$, $s_{24}$, $s_2s_{25}s_2^{-1}$ \\\hline

$s_{11}$ &    $s_{2}$$s_{11}$$s_{2}^{-1},$
$s_{2}$$s_{12}$$s_{2}^{-1},$ $s_2s_{13}s_2^{-1}$, $s_{14}$,
$s_{15}$, $s_{16}$, $s_{17}$, $s_{18}$, $s_{19}$, $s_{20}$,
$s_{21}$, $s_{22}$, $s_{23}$, $s_{24}$, $s_{25}$ \\\hline

$s_{12}$ &      $s_{2}$$s_{12}$$s_{2}^{-1},$ $s_8s_{13}s_8^{-1}$,
$s_2s_{14}s_2^{-1}$    $s_{15}$, $s_{16}$, $s_{17}$, $s_{18}$,
$s_{19}$, $s_{20}$, $s_{21}$, $s_{22}$, $s_{23}$, $s_{24}$, $s_{25}$
\\\hline

$s_{13}$ &      $s_{2}$$s_{13}$$s_{2}^{-1},$
$s_{2}$$s_{14}$$s_{2}^{-1},$    $s_{15}$, $s_{16}$, $s_{17}$,
$s_{18}$, $s_{19}$, $s_{20}$, $s_{21}$, $s_{22}$, $s_{23}$,
$s_{24}$, $s_{25}$ \\\hline

$s_{14}$ &      $s_{2}$$s_{14}$$s_{2}^{-1},$     $s_{15}$, $s_{16}$,
$s_{17}$, $s_{18}$, $s_{19}$, $s_{20}$, $s_{21}$, $s_{22}$,
$s_{23}$, $s_{24}$, $s_{25}$ \\\hline

$s_{15}$ &      $s_{5}$$s_{15}$$s_{5}^{-1},$
$s_{5}$$s_{16}$$s_{5}^{-1},$     $s_{17}$, $s_8s_{18}s_8^{-1}$,
$s_8s_{19}s_8^{-1}$, $s_{20}$, $s_{21}$, $s_{22}$, $s_{23}$,
$s_{24}$, $s_{25}$ \\\hline

$s_{16}$ &        $s_{2}$$s_{16}$$s_{2}^{-1},$    $s_{17}$,
$s_{18}$, $s_{19}$, $s_{20}$, $s_{21}$, $s_{22}$, $s_{23}$,
$s_{24}$, $s_{25}$ \\\hline

$s_{17}$ &        $s_{2}$$s_{17}$$s_{2}^{-1},$
$s_{2}$$s_{18}$$s_{2}^{-1},$    $s_{19}$, $s_{20}$, $s_{21}$,
$s_{2}$$s_{22}$$s_{2}^{-1},$ $s_{23}$, $s_{24}$, $s_{25}$ \\\hline

$s_{18}$ &        $s_{2}$$s_{18}$$s_{2}^{-1},$ $s_8s_{19}s_8^{-1}$,
$s_{20}$, $s_{21}$, $s_{22}s_{23}$, $s_{24}$, $s_{25}$ \\\hline

$s_{19}$ &        $s_{2}$$s_{19}$$s_{2}^{-1},$       $s_{20}$,
$s_{21}$, $s_{22}$, $s_{23}$, $s_{24}$, $s_{25}$ \\\hline

$s_{20}$ &        $s_{2}$$s_{20}$$s_{2}^{-1},$
$s_{2}$$s_{21}$$s_{2}^{-1},$    $s_{22}$,  $s_{23}$, $s_{24}$,
$s_{25}$ \\\hline

$s_{21}$ &        $s_{2}$$s_{21}$$s_{2}^{-1},$ $s_8s_{22}s_8^{-1}$,
$s_{2}$$s_{23}$$s_{2}^{-1},$ $s_{2}$$s_{24}$$s_{2}^{-1},$ $s_{25}$
\\\hline

$s_{22}$ & $s_{11}$$s_{22}$$s_{11}^{-1}, $ $s_8s_{23}s_8^{-1}$,
$s_8s_{24}s_8^{-1}$, $s_{25}$

\\\hline

$s_{23}$ & $s_{2}$$s_{23}$$s_{2}^{-1},$ $s_{2}$$s_{24}$$s_{2}^{-1},$
$s_{25}$

\\\hline

$s_{24}$ & $s_{2}$$s_{24}$$s_{2}^{-1},$   $s_{25}$\\\hline

$s_{25}$ & $s_{2}$$s_{25}$$s_{2}^{-1}$\\\hline

\end{tabular}

$$ \hbox {Table } 13$$

\begin{tabular}{|l|l|}
  \hline
  $E_{7}$ &   \\\hline

 $s_i$ & $s_i$ and $s_t s_js_t^{-1}$ are not square-commutative \\\hline
  $s_{2}$ & $s_{60}$$s_{2}$$s_{60}^{-1},$  $s_{60}$$s_{3}$$s_{60}^{-1},$
   $s_{60}$$s_{4}$$s_{60}^{-1},$ $s_{60}$$s_{5}$$s_{60}^{-1},$ \\\hline

 $s_{3}$ & $s_{60}$$s_{3}$$s_{60}^{-1},$   $s_{60}$$s_{4}$$s_{60}^{-1}$,
  $s_{60}$$s_{5}$$s_{60}^{-1},$ \\\hline

$s_{4}$ & $s_{60}$$s_{4}$$s_{60}^{-1},$
$s_{59}$$s_{5}$$s_{59}^{-1}$,
  $s_{59}$$s_{9}$$s_{59}^{-1}$,   $s_{44}$$s_{11}$$s_{44}^{-1},$
   $s_{44}$$s_{13}$$s_{44}^{-1}$ \\\hline

$s_{5}$ & $s_{60}$$s_{5}$$s_{60}^{-1},$
$s_{44}$$s_{9}$$s_{44}^{-1}$,
 $s_{44}$$s_{11}$$s_{44}^{-1},$   $s_{44}$$s_{13}$$s_{44}^{-1}$,
  $s_{44}$$s_{34}$$s_{44}^{-1}$, $s_{44}$$s_{57}$$s_{44}^{-1}$\\\hline

$s_{7}$ & $s_{60}$$s_{7}$$s_{60}^{-1},$
$s_{44}$$s_{8}$$s_{44}^{-1}$,
 $s_{44}$$s_{9}$$s_{44}^{-1}$,  $s_{44}$$s_{10}$$s_{44}^{-1}$,
  $s_{44}$$s_{11}$$s_{44}^{-1},$   $s_{44}$$s_{12}$$s_{44}^{-1}$,
  $s_{44}$$s_{13}$$s_{44}^{-1}$,  \\ &$s_{44}$$s_{14}$$s_{44}^{-1},$ $s_{44}$$s_{15}$$s_{44}^{-1}$, $s_{44}$$s_{21}$$s_{44}^{-1}$ \\\hline

$s_{8}$ &   $s_{60}$$s_{8}$$s_{60}^{-1}$,  $s_{44}$$s_{9}$$s_{44}^{-1}$,  $s_{44}$$s_{10}$$s_{44}^{-1}$,  $s_{44}$$s_{11}$$s_{44}^{-1},$   $s_{44}$$s_{12}$$s_{44}^{-1}$,  $s_{44}$$s_{13}$$s_{44}^{-1}$,  $s_{44}$$s_{14}$$s_{44}^{-1},$ \\
  &$s_{44}$$s_{15}$$s_{44}^{-1}$,  $s_{44}$$s_{16}$$s_{44}^{-1}$,
  $s_{44}$$s_{18}$$s_{44}^{-1}$\\\hline

$s_{9}$ &   $s_{2}$$s_{9}$$s_{2}^{-1}$,
$s_{44}$$s_{10}$$s_{44}^{-1}$,
 $s_{44}$$s_{12}$$s_{44}^{-1}$,  $s_{2}$$s_{14}$$s_{2}^{-1}$,  $s_{2}$$s_{15}$$s_{2}^{-1}$,  $s_{2}$$s_{16}$$s_{2}^{-1},$\\
  &  $s_{2}$$s_{17}$$s_{2}^{-1}$, $s_{2}$$s_{18}$$s_{2}^{-1}$,
   $s_{2}$$s_{19}$$s_{2}^{-1}$,  $s_{2}$$s_{21}$$s_{2}^{-1}$,
    $s_{2}$$s_{24}$$s_{2}^{-1}$,  $s_{2}$$s_{26}$$s_{2}^{-1}$,
  $s_{2}$$s_{45}$$s_{2}^{-1}$\\\hline

$s_{10}$ &   $s_{2}$$s_{10}$$s_{2}^{-1}$,
$s_{2}$$s_{11}$$s_{2}^{-1},$
 $s_{2}$$s_{12}$$s_{2}^{-1}$,  $s_{2}$$s_{13}$$s_{2}^{-1}$,  $s_{2}$$s_{14}$$s_{2}^{-1}$,
  $s_{2}$$s_{15}$$s_{2}^{-1}$,    $s_{2}$$s_{21}$$s_{2}^{-1}$,
    $s_{2}$$s_{45}$$s_{2}^{-1}$\\\hline

$s_{11}$ &     $s_{3}$$s_{11}$$s_{3}^{-1},$
$s_{2}$$s_{12}$$s_{2}^{-1}$, $s_{3}$$s_{13}$$s_{3}^{-1}$,
$s_{2}$$s_{14}$$s_{2}^{-1}$,  $s_{2}$$s_{15}$$s_{2}^{-1}$,
$s_{2}$$s_{16}$$s_{2}^{-1},$ $s_{3}$$s_{17}$$s_{3}^{-1}$,
$s_{3}$$s_{18}$$s_{3}^{-1}$, \\&
  $s_{2}$$s_{20}$$s_{2}^{-1},$   $s_{3}$$s_{21}$$s_{3}^{-1}$,
  $s_{2}$$s_{25}$$s_{2}^{-1}$,  $s_{3}$$s_{26}$$s_{3}^{-1}$,  $s_{2}$$s_{27}$$s_{2}^{-1},$
  $s_{2}$$s_{28}$$s_{2}^{-1}$,  $s_{2}$$s_{29}$$s_{2}^{-1}$,
  $s_{3}$$s_{30}$$s_{3}^{-1}$,  \\& $s_{3}$$s_{31}$$s_{3}^{-1}$, $s_{2}$$s_{34}$$s_{2}^{-1}$,
   $s_{2}$$s_{36}$$s_{2}^{-1}$,  $s_{2}$$s_{37}$$s_{2}^{-1},$
   $s_{2}$$s_{38}$$s_{2}^{-1}$,  $s_{2}$$s_{39}$$s_{2}^{-1}$,  $s_{2}$$s_{40}$$s_{2}^{-1}$,
   $s_{2}$$s_{41}$$s_{2}^{-1}$, \\ & $s_{2}$$s_{44}$$s_{2}^{-1}$,  $s_{2}$$s_{45}$$s_{2}^{-1}$,
   $s_{2}$$s_{54}$$s_{2}^{-1}$, $s_{2}$$s_{56}$$s_{2}^{-1}$,  $s_{2}$$s_{57}$$s_{2}^{-1}$,
     $s_{2}$$s_{59}$$s_{2}^{-1}$,  $s_{2}$$s_{60}$$s_{2}^{-1}$\\\hline

$s_{12}$ &      $s_{2}$$s_{12}$$s_{2}^{-1}$,
$s_{2}$$s_{13}$$s_{2}^{-1}$, $s_{2}$$s_{14}$$s_{2}^{-1}$,
$s_{2}$$s_{15}$$s_{2}^{-1}$,  $s_{2}$$s_{16}$$s_{2}^{-1},$
$s_{2}$$s_{45}$$s_{2}^{-1}$\\\hline

$s_{13}$ &      $s_{3}$$s_{13}$$s_{3}^{-1}$,
$s_{2}$$s_{14}$$s_{2}^{-1}$, $s_{2}$$s_{15}$$s_{2}^{-1}$,
$s_{2}$$s_{16}$$s_{2}^{-1},$ $s_{3}$$s_{17}$$s_{3}^{-1}$,
$s_{3}$$s_{18}$$s_{3}^{-1}$,
  $s_{2}$$s_{20}$$s_{2}^{-1},$   $s_{3}$$s_{21}$$s_{3}^{-1}$, \\ &  $s_{2}$$s_{25}$$s_{2}^{-1}$,
   $s_{3}$$s_{26}$$s_{3}^{-1}$,  $s_{2}$$s_{27}$$s_{2}^{-1},$
  $s_{2}$$s_{28}$$s_{2}^{-1}$,  $s_{2}$$s_{29}$$s_{2}^{-1}$,
   $s_{3}$$s_{30}$$s_{3}^{-1}$,  $s_{3}$$s_{31}$$s_{3}^{-1}$, \\ &$s_{2}$$s_{34}$$s_{2}^{-1}$,
   $s_{2}$$s_{36}$$s_{2}^{-1}$,  $s_{2}$$s_{37}$$s_{2}^{-1},$
   $s_{2}$$s_{38}$$s_{2}^{-1}$,  $s_{2}$$s_{39}$$s_{2}^{-1}$,  $s_{2}$$s_{40}$$s_{2}^{-1}$,
  $s_{2}$$s_{41}$$s_{2}^{-1}$, \\ & $s_{2}$$s_{44}$$s_{2}^{-1}$,  $s_{2}$$s_{45}$$s_{2}^{-1}$,
   $s_{2}$$s_{54}$$s_{2}^{-1}$,
  $s_{2}$$s_{56}$$s_{2}^{-1}$,  $s_{2}$$s_{57}$$s_{2}^{-1}$,
  $s_{2}$$s_{59}$$s_{2}^{-1}$,  $s_{2}$$s_{60}$$s_{2}^{-1}$\\\hline

$s_{14}$ &       $s_{2}$$s_{14}$$s_{2}^{-1}$,
$s_{2}$$s_{15}$$s_{2}^{-1}$, $s_{2}$$s_{16}$$s_{2}^{-1}$,
$s_{2}$$s_{17}$$s_{2}^{-1}$,     $s_{2}$$s_{21}$$s_{2}^{-1}$,
 $s_{2}$$s_{25}$$s_{2}^{-1}$,   $s_{2}$$s_{29}$$s_{2}^{-1}$,
  $s_{2}$$s_{36}$$s_{2}^{-1}$,  \\ & $s_{2}$$s_{37}$$s_{2}^{-1}$,
   $s_{2}$$s_{45}$$s_{2}^{-1}$\\\hline

$s_{15}$ &     $s_{3}$$s_{15}$$s_{3}^{-1}$,
$s_{2}$$s_{16}$$s_{2}^{-1}$,
 $s_{2}$$s_{17}$$s_{2}^{-1}$,  $s_{2}$$s_{18}$$s_{2}^{-1}$,
  $s_{2}$$s_{27}$$s_{2}^{-1}$, $s_{2}$$s_{29}$$s_{2}^{-1}$,
   $s_{2}$$s_{36}$$s_{2}^{-1}$,  $s_{2}$$s_{37}$$s_{2}^{-1}$, \\ &
    $s_{2}$$s_{42}$$s_{2}^{-1}$,  $s_{2}$$s_{43}$$s_{2}^{-1}$,
   $s_{2}$$s_{56}$$s_{2}^{-1}$\\\hline

$s_{16}$ &    $s_{2}$$s_{16}$$s_{2}^{-1}$,
$s_{2}$$s_{17}$$s_{2}^{-1}$,  $s_{2}$$s_{18}$$s_{2}^{-1}$,
$s_{3}$$s_{19}$$s_{3}^{-1}$,  $s_{2}$$s_{20}$$s_{2}^{-1}$,
$s_{2}$$s_{21}$$s_{2}^{-1},$
 $s_{3}$$s_{26}$$s_{3}^{-1}$, $s_{2}$$s_{28}$$s_{2}^{-1}$, \\ & $s_{2}$$s_{29}$$s_{2}^{-1}$,
$s_{2}$$s_{45}$$s_{2}^{-1}$,  $s_{2}$$s_{55}$$s_{2}^{-1}$,
$s_{2}$$s_{57}$$s_{2}^{-1}$,  $s_{2}$$s_{58}$$s_{2}^{-1}$\\\hline

$s_{17}$ &    $s_{3}$$s_{17}$$s_{3}^{-1}$,
$s_{3}$$s_{18}$$s_{3}^{-1}, $, $s_{23}$$s_{19}$$s_{23}^{-1}$,
$s_{2}$$s_{20}$$s_{2}^{-1}$,
 $s_{3}$$s_{21}$$s_{3}^{-1},$ $s_{2}$$s_{25}$$s_{2}^{-1}$,  $s_{3}$$s_{26}$$s_{3}^{-1}$,
  $s_{2}$$s_{28}$$s_{2}^{-1}$,  \\ & $s_{2}$$s_{29}$$s_{2}^{-1}$,  $s_{3}$$s_{30}$$s_{3}^{-1},$
   $s_{3}$$s_{31}$$s_{3}^{-1}$,  $s_{2}$$s_{36}$$s_{2}^{-1}$,
   $s_{2}$$s_{45}$$s_{2}^{-1}$,  $s_{2}$$s_{56}$$s_{2}^{-1}$,  $s_{2}$$s_{57}$$s_{2}^{-1}$,
    $s_{2}$$s_{60}$$s_{2}^{-1},$\\\hline

$s_{18}$ &    $s_{2}$$s_{18}$$s_{2}^{-1}$,
$s_{3}$$s_{19}$$s_{3}^{-1}, $, $s_{2}$$s_{20}$$s_{2}^{-1}$,
$s_{2}$$s_{21}$$s_{2}^{-1},$ $s_{3}$$s_{26}$$s_{3}^{-1},$
   $s_{2}$$s_{28}$$s_{2}^{-1}$,  $s_{2}$$s_{29}$$s_{2}^{-1}$,  $s_2$$s_{30}$$s_{2}^{-1},$
   \\ &
      $s_2s_{44}s_2^{-1}$,   $s_{2}$$s_{45}$$s_{2}^{-1},
    $, $s_{2}$$s_{57}$$s_{2}^{-1}$
  \\\hline

$s_{19}$ &     $s_{3}$$s_{19}$$s_{3}^{-1}$,
$s_{2}$$s_{20}$$s_{2}^{-1}$,  $s_{2}$$s_{21}$$s_{2}^{-1},$
$s_{3}$$s_{26}$$s_{3}^{-1},$
   $s_3s_{27}s_3^{-1}$,   $s_{2}$$s_{28}$$s_{2}^{-1}$,  \\ & $s_{2}$$s_{29}$$s_{2}^{-1},$
    $s_{2}$$s_{30}$$s_{2}^{-1},$   $s_{2}$$s_{31}$$s_{2}^{-1}$,  $s_2s_{32}s_2^{-1}$,
      $s_2s_{33}s_2^{-1}$,  $s_2s_{44}s_2^{-1}$,  $s_2s_{46}s_2^{-1}$,  $s_2s_{54}s_2^{-1}$,
    \\ &  $s_2s_{55}s_2^{-1}$,
 $s_{3}$$s_{57}$$s_{3}^{-1}$,  $s_2s_{58}s_2^{-1}$,  $s_2s_{59}s_2^{-1}$,
   $s_{2}$$s_{60}$$s_{2}^{-1}$\\\hline

\end{tabular}

$$ \hbox {Table } 14$$

\begin{tabular}{|l|l|}
  \hline
  $E_{7}$ &   \\\hline
$s_{20}$ &    $s_{2}$$s_{20}$$s_{2}^{-1}$,
$s_{2}$$s_{21}$$s_{2}^{-1},$
 $s_{2}$$s_{26}$$s_{2}^{-1}$,  $s_2s_{27}s_2^{-1}$,
  $s_{2}$$s_{29}$$s_{2}^{-1}$,  $s_2s_{55}s_2^{-1}$,
   $s_2s_{58}s_2^{-1}$ \\\hline

$s_{21}$ &     $s_{2}$$s_{21}$$s_{2}^{-1},$
$s_{3}$$s_{26}$$s_{3}^{-1},$
   $s_{2}$$s_{28}$$s_{2}^{-1}$,  $s_{2}$$s_{29}$$s_{2}^{-1},$
    $s_{2}$$s_{30}$$s_{2}^{-1},$   $s_2s_{32}s_2^{-1}$,   $s_2s_{40}s_2^{-1}$,
    $s_2s_{44}s_2^{-1}$, \\ & $s_2s_{55}s_2^{-1}$,
  $s_{2}$$s_{57}$$s_{2}^{-1}$,  $s_2s_{58}s_2^{-1}$\\\hline

$s_{22}$ &        $s_2s_{22}s_2^{-1}$,  $s_2s_{23}s_2^{-1}$,
$s_2s_{24}s_2^{-1}$,  $s_{2}$$s_{45}$$s_{2}^{-1}$\\\hline

$s_{23}$ &        $s_2s_{23}s_2^{-1}$, $s_2s_{24}s_2^{-1}$,
 $s_{2}$$s_{45}$$s_{2}^{-1} $\\\hline

$s_{24}$ &        $s_2s_{24}s_2^{-1}$,  $s_{2}$$s_{45}$$s_{2}^{-1}
$\\\hline

$s_{25}$ &        $s_2s_{25}s_2^{-1}$,  $s_2s_{49}s_2^{-1}$,
$s_2s_{56}s_2^{-1}$, $s_2s_{57}s_2^{-1}$\\\hline

$s_{26}$ &         $s_3s_{26}s_3^{-1}$,
 $s_2s_{27}s_2^{-1}$,  $s_2s_{28}s_2^{-1}$,  $s_2s_{29}s_2^{-1}$,  $s_2s_{30}s_2^{-1}$,
 $s_2s_{31}s_2^{-1}$,  $s_2s_{44}s_2^{-1}$, \\ & $s_2s_{54}s_2^{-1}$,
 $s_2s_{57}s_2^{-1}$,
    $s_2s_{59}s_2^{-1}$,   $s_2s_{60}s_2^{-1}$\\\hline

$s_{27}$ &
 $s_2s_{27}s_2^{-1}$,   $s_2s_{29}s_2^{-1}$,   $s_2s_{38}s_2^{-1}$,  $s_2s_{39}s_2^{-1}$,
   $s_{2}$$s_{45}$$s_{2}^{-1},$   $s_2s_{59}s_2^{-1}$\\\hline

$s_{28}$ &     $s_2s_{28}s_2^{-1}$,  $s_2s_{29}s_2^{-1}$,
$s_2s_{30}s_2^{-1}$,
 $s_2s_{57}s_2^{-1}$\\\hline

$s_{29}$ &        $s_2s_{29}s_2^{-1}$,   $s_2s_{30}s_2^{-1}$,
$s_2s_{31}s_2^{-1}$,
 $s_2s_{32}s_2^{-1}$,  $s_2s_{33}s_2^{-1}$,  $s_2s_{34}s_2^{-1}$,  $s_2s_{36}s_2^{-1}$,
  $s_{2}$$s_{45}$$s_{2}^{-1} $,  $s_2s_{59}s_2^{-1}$\\\hline

$s_{30}$ &         $s_2s_{30}s_2^{-1}$, $s_2s_{31}s_2^{-1}$,
  $s_2s_{34}s_2^{-1}$\\\hline

$s_{31}$ &         $s_2s_{31}s_2^{-1}$\\\hline

$s_{32}$ &
 $s_2s_{32}s_2^{-1}$,  $s_2s_{33}s_2^{-1}$,  $s_2s_{46}s_2^{-1}$,
  $s_2s_{55}s_2^{-1}$,
   $s_2s_{58}s_2^{-1}$\\\hline

$s_{33}$ &         $s_2s_{33}s_2^{-1}$,  $s_2s_{46}s_2^{-1}$\\\hline

 $s_{34}$ &        $s_2s_{34}s_2^{-1}$, $s_2s_{35}s_2^{-1}$,  $s_2s_{53}s_2^{-1}$,
  $s_2s_{59}s_2^{-1}$\\\hline

   $s_{35}$ &         $s_2s_{35}s_2^{-1}$,    $s_2s_{53}s_2^{-1}$\\\hline

   $s_{36}$ &        $s_2s_{36}s_2^{-1}$,
 $s_2s_{37}s_2^{-1}$, $s_2s_{43}s_2^{-1}$\\\hline

   $s_{37}$ &
 $s_2s_{37}s_2^{-1}$, $s_3s_{43}s_3^{-1}$\\\hline

 $s_{38}$ &
 $s_2s_{38}s_2^{-1}$,  $s_2s_{39}s_2^{-1}$\\\hline

 $s_{39}$ &
 $s_2s_{39}s_2^{-1}$\\\hline

 $s_{40}$ &
  $s_2s_{40}s_2^{-1}$, $s_2s_{41}s_2^{-1}$\\\hline

 $s_{41}$ &
 $s_2s_{41}s_2^{-1}$,
     $s_2s_{59}s_2^{-1}$\\\hline

 $s_{42}$ &
     $s_2s_{42}s_2^{-1}$\\\hline

 $s_{43}$ &
 $s_2s_{43}s_2^{-1}$\\\hline

 $s_{44}$ &
 $s_2s_{44}s_2^{-1}$\\\hline

 $s_{45}$ &
 $s_2s_{45}s_2^{-1}$\\\hline

 $s_{46}$ &
 $s_2s_{46}s_2^{-1}$\\\hline

 $s_{47}$ &
 $s_2s_{47}s_2^{-1}$,
      $s_2s_{48}s_2^{-1}$,  $s_2s_{49}s_2^{-1}$,  $s_2s_{50}s_2^{-1}$\\\hline

 $s_{48}$ &
 $s_2s_{48}s_2^{-1}$,    $s_2s_{49}s_2^{-1}$,  $s_2s_{50}s_2^{-1}$\\\hline

 $s_{49}$ &
 $s_2s_{49}s_2^{-1}$,    $s_2s_{50}s_2^{-1}$\\\hline

 $s_{50}$ &
 $s_2s_{50}s_2^{-1}$\\\hline

 $s_{51}$ &
 $s_2s_{51}s_2^{-1}$,  $s_2s_{52}s_2^{-1}$\\\hline

 $s_{52}$ &
 $s_2s_{52}s_2^{-1}$\\\hline

\end{tabular}

$$ \hbox {Table }  15$$

\begin{tabular}{|l|l|}
  \hline
  $E_{7}$ &   \\\hline

 $s_{53}$ &
 $s_2s_{53}s_2^{-1}$\\\hline

 $s_{54}$ &
 $s_2s_{54}s_2^{-1}$\\\hline

 $s_{55}$ &
 $s_2s_{55}s_2^{-1}$,
   $s_2s_{58}s_2^{-1}$\\\hline

 $s_{56}$ &
 $s_2s_{56}s_2^{-1}$,  $s_2s_{57}s_2^{-1}$\\\hline

 $s_{57}$ &
 $s_2s_{57}s_2^{-1}$\\\hline

 $s_{58}$ &
 $s_2s_{58}s_2^{-1}$\\\hline

 $s_{59}$ &
 $s_2s_{59}s_2^{-1}$\\\hline

 $s_{60}$ &
 $s_2s_{60}s_2^{-1}$ \\\hline

\end{tabular}

$$ \hbox {Table }  16$$

\vskip 1.0cm

\begin{tabular}{|l|l|}
  \hline
  $E_{8}$ &   \\\hline

 $s_i$ & $s_i$ and $s_t s_js_t^{-1}$ are not square-commutative \\\hline

$s_{5}$ &       $s_{41}s_{5}s_{41}^{-1}$, $s_{54}s_{15}s_{54}^{-1}$,
   $s_{112}s_{18}s_{112}^{-1}$, $s_{9}s_{26}s_{9}^{-1}$, $s_{2}s_{38}s_{2}^{-1}$,
   $s_{2}s_{106}s_{2}^{-1}$,
    \\\hline

   $s_{6}$ &         $s_{2}s_{12}s_{2}^{-1}$, $s_{2}s_{44}s_{2}^{-1}$  \\\hline

   $s_{8}$ &         $s_{41}s_{8}s_{41}^{-1}$, $s_{112}s_{12}s_{112}^{-1}$,
   $s_{2}s_{22}s_{2}^{-1},$ $s_{9}s_{26}s_{9}^{-1},$  $s_{2}s_{38}s_{2}^{-1}
   $\\\hline

   $s_{12}$ &   $s_{2}s_{24}s_{2}^{-1}$, $s_{2}s_{26}s_{2}^{-1}$,
    $s_{2}s_{50}s_{2}^{-1}$,  $s_{2}s_{51}s_{2}^{-1}$,
           $s_{2}s_{62}s_{2}^{-1}$,       $s_{41}s_{80}s_{41}^{-1}$\\\hline

   $s_{14}$ &         $s_{2}s_{14}s_{2}^{-1}$, $s_{41}s_{21}s_{41}^{-1}$,
   $s_{9}s_{26}s_{9}^{-1}$, $s_{2}s_{32}s_{2}^{-1}$,
   $s_{2}s_{38}s_{2}^{-1},$ $s_{2}s_{39}s_{2}^{-1},$
   $s_{2}s_{53}s_{2}^{-1}$,     $s_{9}s_{56}s_{9}^{-1}$, \\ & $s_{9}s_{57}s_{9}^{-1}$,
 $s_{2}s_{58}s_{2}^{-1}$, $s_{9}s_{68}s_{9}^{-1}$, $s_{70}$, $s_{9}s_{108}s_{9}^{-1}$
   \\\hline

$s_{15}$ &         $s_{2}s_{26}s_{2}^{-1}$\\\hline

   $s_{18}$ &    $s_{2}s_{24}s_{2}^{-1}$, $s_{2}s_{26}s_{2}^{-1}$,
    $s_{2}s_{62}s_{2}^{-1}$,  $s_{2}s_{74}s_{2}^{-1}$      $s_{41}s_{80}s_{41}^{-1}$\\\hline

$s_{21}$ &         $s_{41}s_{21}s_{41}^{-1}$,
$s_{2}s_{53}s_{2}^{-1}$\\\hline

$s_{22}$ &         $s_{2}s_{56}s_{2}^{-1}$, $s_{2}s_{57}s_{2}^{-1}$,
   $s_{2}s_{70}s_{2}^{-1}$\\\hline

$s_{24}$ &         $s_{41}s_{38}s_{41}^{-1}$,
$s_{2}s_{39}s_{2}^{-1}$,
   $s_{41}s_{42}s_{41}^{-1}$,      $s_{2}s_{51}s_{2}^{-1}$,
$s_{35}s_{53}s_{35}^{-1}$,
    $s_{2}s_{105}s_{2}^{-1}$,
   $s_{9}s_{106}s_{9}^{-1}$\\\hline

$s_{26}$ &         $s_{9}s_{35}s_{9}^{-1}$, $s_{2}s_{42}s_{2}^{-1}$,
$s_{2}s_{51}s_{2}^{-1}$, $s_{2}s_{106}s_{2}^{-1}$
   \\\hline

$s_{35}$ &        $s_{2}s_{48}s_{2}^{-1}$
   \\\hline

   $s_{45}$ &        $s_{2}s_{75}s_{2}^{-1}$,
$s_{2}s_{80}s_{2}^{-1}$, $s_{2}s_{108}s_{2}^{-1}$
   \\\hline
$s_{80}$ &         $s_{2}s_{106}s_{2}^{-1}$
   \\\hline

$s_{110}$ &         $s_{2}s_{110}s_{2}^{-1}$
\\\hline
\end{tabular}

$$ \hbox {Table }  17$$

\vskip 0.3cm
\begin{tabular}{|l|l|}
  \hline
  $F_4$ &   \\\hline

   $s_i$ & $s_i$ and $s_t s_js_t^{-1}$ are not square-commutative \\\hline
$s_3$  &   $ s_{23}s_9 s_{23}^{-1}$, $ s_{23}s_{19} s_{23}^{-1}$, $
s_{5}s_{23} s_{5}^{-1}$,  \\\hline

$s_4$ & $ s_{14}s_{22} s_{14}^{-1}$,  $ s_{20}s_{23}
s_{20}^{-1}$\\\hline

$s_5$ &   $ s_{23}s_{12} s_{23}^{-1}$,  $ s_{23}s_{13}
s_{23}^{-1}$\\\hline

$s_6$ &   $ s_{23}s_{12} s_{23}^{-1}$,  $ s_{23}s_{13}
s_{23}^{-1}$\\\hline

$s_7$ &   $ s_{14}s_{7} s_{14}^{-1}$,  $ s_{14}s_{8} s_{14}^{-1}$, $
s_{14}s_{9} s_{14}^{-1}$,
 $ s_{14}s_{22} s_{14}^{-1}$, $ s_{14}s_{24} s_{14}^{-1}$,
$ s_{14}s_{25} s_{14}^{-1}$  \\\hline

$s_8$ &  $ s_{14}s_{8} s_{14}^{-1}$, $ s_{14}s_{9} s_{14}^{-1}$, $
s_{14}s_{22} s_{14}^{-1}$, $ s_{14}s_{24} s_{14}^{-1}$, $
s_{14}s_{25} s_{14}^{-1}$
\\\hline

   $s_{9}$ &  $ s_{14}s_{9} s_{14}^{-1}$, $ s_{14}s_{19}
   s_{14}^{-1}$, $ s_{23}s_{23} s_{23}^{-1}$
         \\\hline

   $s_{10}$ &     $ s_{23}s_{12} s_{23}^{-1}$, $ s_{23}s_{13} s_{23}^{-1}$
   \\\hline

   $s_{11}$ &     $ s_{23}s_{12} s_{23}^{-1}$, $ s_{23}s_{13} s_{23}^{-1}$
   \\\hline

   $s_{12}$ &     $ s_{23}s_{24} s_{23}^{-1}$, $ s_{23}s_{25} s_{23}^{-1}$
   \\\hline

     $s_{13}$ &     $ s_{23}s_{24} s_{23}^{-1}$, $ s_{23}s_{25} s_{23}^{-1}$
   \\\hline

     $s_{14}$ &     $ s_{14}s_{22} s_{14}^{-1}$, $ s_{23}s_{23} s_{23}^{-1}$
   \\\hline

  $s_{15}$ &     $ s_{23}s_{20} s_{23}^{-1}$, $ s_{23}s_{21} s_{23}^{-1}$
   \\\hline

  $s_{16}$ &     $ s_{23}s_{20} s_{23}^{-1}$, $ s_{23}s_{21} s_{23}^{-1}$
   \\\hline
  $s_{17}$ &    $ s_{14}s_{17} s_{14}^{-1}$, $ s_{14}s_{18} s_{14}^{-1}$,
  $ s_{14}s_{19} s_{14}^{-1}$, $ s_{14}s_{22} s_{14}^{-1}$,
  $ s_{14}s_{24} s_{14}^{-1}$,  $ s_{14}s_{25} s_{14}^{-1}$
   \\\hline
$s_{18}$ &    $ s_{14}s_{18} s_{14}^{-1}$,
  $ s_{14}s_{19} s_{14}^{-1}$, $ s_{14}s_{22} s_{14}^{-1}$,
  $ s_{14}s_{24} s_{14}^{-1}$,  $ s_{14}s_{25} s_{14}^{-1}$
   \\\hline

   $s_{19}$ &
  $ s_{14}s_{19} s_{14}^{-1}$, $ s_{14}s_{23} s_{14}^{-1}$
   \\\hline

$s_{22}$ & $ s_{14}s_{22} s_{14}^{-1}$,
  $ s_{14}s_{23} s_{14}^{-1}$, $ s_{14}s_{24} s_{14}^{-1}$,
  $ s_{14}s_{25} s_{14}^{-1}$
   \\\hline

$s_{23}$ &
  $ s_{14}s_{23} s_{14}^{-1}$, $ s_{14}s_{24} s_{14}^{-1}$,
  $ s_{14}s_{25} s_{14}^{-1}$
   \\\hline
$s_{24}$ &    $ s_{14}s_{24} s_{14}^{-1}$,
  $ s_{14}s_{25} s_{14}^{-1}$
   \\\hline
$s_{25}$ &
  $ s_{14}s_{25} s_{14}^{-1}$
   \\\hline

\end{tabular}
$$ \hbox {Table }  18$$

\vskip 1.0cm

\noindent{\large\bf Acknowledgement}: We would like to thank Prof.
N. Andruskiewitsch and Dr. F. Fantino for suggestions and help. The first  author and the
second author were financially supported by the Australian Research
Council. S.C.Zhang thanks the School of Mathematics and Physics,
The University of Queensland for hospitality.

\end {document}